\documentclass[11pt]{book}
\usepackage{amssymb}
\usepackage{amsmath}
\usepackage{amsthm}
\usepackage{graphics}
\usepackage{epsfig,amssymb,latexsym,verbatim}
\usepackage{graphicx}
\usepackage{multirow}
\usepackage{multicol}
\usepackage{hypernat}
\usepackage{float}
\usepackage{hyperref}
\floatstyle{ruled}
\newfloat{algorithm}{tbp}{loa}
\floatname{algorithm}{Algorithm}
\usepackage{graphicx, amssymb}

\newtheorem{theorem}{{\bf Theorem}}
\newtheorem{lemma}{{\bf Lemma}}

\newtheorem{proposition}{{\bf Proposition}}

\begin{document}
\renewcommand{\baselinestretch}{1.2}
\markboth{\hfill{\footnotesize\rm LI WANG AND LIJIAN YANG}\hfill}
{\hfill {\footnotesize\rm SINGLE-INDEX PREDICTION MODEL} \hfill}
\renewcommand{\thefootnote}{}
$\ $\par \fontsize{10.95}{14pt plus.8pt minus .6pt}\selectfont
\vspace{0.8pc} \centerline{\large\bf SPLINE SINGLE-INDEX
PREDICTION MODEL} %\vspace{2pt}
%\centerline{\large\bf IF A SECOND LINE IS NEEDED}
\vspace{.4cm} \centerline{Li Wang and Lijian Yang
\footnote{\emph{Address for correspondence}: Lijian Yang,
Department of Statistics and Probability, Michigan State
University, East Lansing, MI 48824, USA. E-mail:
yang@stt.msu.edu}} \vspace{.4cm} \centerline{\it University of
Georgia and Michigan State University} \vspace{.55cm}
\fontsize{9}{11.5pt plus.8pt minus .6pt}\selectfont

\begin{quotation}
\noindent \textit{Abstract:} For the past two decades,
single-index model, a special case of projection pursuit
regression, has proven to be an efficient way of coping with the
high dimensional problem in nonparametric regression. In this
paper, based on weakly dependent sample, we investigate the
single-index prediction (SIP) model which is robust against
deviation from the single-index model. The single-index is
identified by the best approximation to the multivariate
prediction function of the response variable, regardless of
whether the prediction function is a genuine single-index
function. A polynomial spline estimator is proposed for the
single-index prediction coefficients, and is shown to be root-n
consistent and asymptotically normal. An iterative optimization
routine is used which is sufficiently fast for the user to analyze
large data of high dimension within seconds. Simulation
experiments have provided strong evidence that corroborates with
the asymptotic theory. Application of the proposed procedure to
the rive flow data of Iceland has yielded superior out-of-sample
rolling forecasts.

\vspace{9pt} \noindent \textit{Key words and phrases:} B-spline, geometric
mixing, knots, nonparametric regression, root-n rate, strong consistency.
\end{quotation}

\fontsize{10.95}{14pt plus.8pt minus .6pt}\selectfont

\thispagestyle{empty}

\setcounter{chapter}{1} \label{SEC:introduction}
\setcounter{equation}{0}
%-1
\noindent \textbf{1. Introduction} \vskip 0.1in

Let $\left\{ \mathbf{X}_{i}^{T},Y_{i}\right\} _{i=1}^{n}=\left\{
X_{i,1},...,X_{i,d},Y_{i}\right\} _{i=1}^{n}$ be a length $n$ realization of
a $\left( d+1\right) $-dimensional strictly stationary process following the
heteroscedastic model
\begin{equation}
Y_{i}=m\left( \mathbf{X}_{i}\right) +\sigma \left( \mathbf{X}_{i}\right)
\varepsilon _{i},m\left( \mathbf{X}_{i}\right) =E\left( Y_{i}|\mathbf{X}%
_{i}\right) ,  \label{sindmodel}
\end{equation}
in which $E\left( \varepsilon _{i}\left| \mathbf{X}_{i}\right. \right) =0$, $%
E\left( \varepsilon _{i}^{2}\left| \mathbf{X}_{i}\right. \right)
=1$, $1\leq i\leq n$. The $d$-variate functions $m$, $\sigma $ are
the unknown mean and standard deviation of the response $Y_{i}$
conditional on the predictor vector $\mathbf{X}_{i}$, often
estimated nonparametrically. In what follows, we let $\left(
\mathbf{X}^{T},Y,\varepsilon \right) $ have the stationary
distribution of $\left( \mathbf{X}_{i}^{T},Y_{i},\varepsilon
_{i}\right) $. When the dimension of $\mathbf{X}$ is high, one
unavoidable issue is the ``curse of dimensionality'', which refers
to the poor convergence rate of nonparametric estimation of
general multivariate function. Much effort has been devoted to the
circumventing of this difficulty. In the words of Xia, Tong, Li
and Zhu (2002), there are essentially two approaches: function
approximation and dimension reduction. A favorite function
approximation technique is the generalized additive model
advocated by Hastie and Tibshirani (1990), see also, for example,
Mammen, Linton and Nielsen (1999), Huang and Yang (2004), Xue and
Yang (2006 a, b), Wang and Yang (2007). An attractive dimension
reduction method is the single-index model, similar to the first
step of projection pursuit regression, see Friedman and Stuetzle
(1981), Hall (1989), Huber (1985), Chen (1991). The basic appeal
of single-index model is its simplicity: the $d $-variate function
$m\left( \mathbf{x}\right) =m\left(
x_{1},...,x_{d}\right) $ is expressed as a univariate function of $\mathbf{x}%
^{T}\mathbf{\theta }_{0}=\sum_{p=1}^{d}x_{p}\theta _{0,p}$. Over
the last two decades, many authors had devised various intelligent
estimators of the single-index coefficient vector $\mathbf{\theta
}_{0}=\left( \theta _{0,1},...,\theta _{0,d}\right) ^{T}$, for
instance, Powell, Stock and Stoker (1989), H\"{a}rdle and Stoker
(1989), Ichimura (1993), Klein and Spady (1993), H\"{a}rdle, Hall
and Ichimura (1993), Horowitz and H\"{a}rdle (1996), Carroll, Fan,
Gijbels and Wand (1997), Xia and Li (1999), Hristache, Juditski
and Spokoiny (2001). More recently, Xia, Tong, Li and Zhu (2002)
proposed the minimum average variance estimation (MAVE) for
several index vectors.

All the aforementioned methods assume that the $d$-variate
regression function $m\left( \mathbf{x}\right) $ is exactly a
univariate function of some $\mathbf{x}^{T}\mathbf{\theta }_{0}$
and obtain a root-$n$ consistent estimator of $\mathbf{\theta
}_{0}$. If this model is misspecified ($m$ is not a genuine
single-index function), however, a goodness-of-fit test then
becomes necessary and the estimation of $\mathbf{\theta }_{0}$
must be redefined, see Xia, Li, Tong and Zhang (2004). In this
paper, instead of presuming that underlying true function $m$ is a
single-index function, we estimate a univariate function $g$ that
optimally approximates the multivariate function $m$ in the sense
of
\begin{equation}
g\left( \nu \right) =E\left[ \left. m\left( \mathbf{X}\right)
\right| \mathbf{X}^{T}\mathbf{\theta }_{0}=\nu \right] ,
\label{DEF:g}
\end{equation}
where the unknown parameter $\mathbf{\theta }_{0}$ is called the SIP
coefficient, used for simple interpretation once estimated; $\mathbf{X}^{T}%
\mathbf{\theta }_{0}$ is the latent SIP variable; and $g$ is a
smooth but unknown function used for further data summary, called
the link prediction function. Our method therefore is clearly
interpretable regardless of the goodness-of-fit of the
single-index model, making it much more relevant in applications.

We propose estimators of $\mathbf{\theta }_{0}$ and $g$ based on
weakly dependent sample, which includes many existing
nonparametric time series models, that are (i) computationally
expedient and (ii) theoretically reliable. Estimation of both
$\mathbf{\theta }_{0}$ and $g$ has been done via the kernel
smoothing techniques in existing literature, while we use
polynomial spline smoothing. The greatest advantages of spline
smoothing, as pointed out in Huang and Yang (2004), Xue and Yang
(2006 b) are its simplicity and fast computation. Our proposed
procedure
involves two stages: estimation of $\mathbf{\theta }_{0}$ by some $\sqrt{n}$%
-consistent $\hat{\mathbf{\theta }}$, minimizing an empirical version of the
mean squared error, $R(\mathbf{\theta })=E\{Y-E(\left. Y\right| \mathbf{X}%
^{T}\mathbf{\theta })\}^{2}$; spline smoothing of $Y$ on $\mathbf{X}^{T}\hat{%
\mathbf{\theta }}$ to obtain a cubic spline estimator $\hat{g}$ of
$g$. The best single-index approximation to $m(\mathbf{x})$ is
then $\hat{m}(\mathbf{x})=\hat{g}\left(
\mathbf{x}^{T}\hat{\mathbf{\theta }}\right)$.

Under geometrically strong mixing condition, strong consistency
and $\sqrt{n}$-rate asymptotic normality of the estimator
$\hat{\mathbf{\theta }}$ of the SIP coefficient $\mathbf{\theta
}_{0}$ in (\ref{DEF:g}) are obtained. Proposition \ref{PROP:unif}
is the key in understanding the efficiency of the proposed
estimator. It shows that the derivatives of the risk function up
to order 2 are uniformly almost surely approximated by their
empirical versions.

Practical performance of the SIP estimators is examined via Monte
Carlo examples. The estimator of the SIP coefficient performs very
well for data of both moderate and high dimension $d$, of sample
size $n$ from small to large, see Tables \ref{TAB:meansdxia} and
\ref{TAB:mse-time}, Figures \ref {FIG:Xia2004} and
\ref{FIG:estimation}. By taking advantages of the spline smoothing
and the iterative optimization routines, one reduces the
computation burden immensely for massive data sets. Table
\ref{TAB:mse-time} reports the computing time of one simulation
example on an ordinary PC, which shows that for massive data sets,
 the SIP method is much faster than the
MAVE method. For instance, the SIP estimation of a
$200$-dimensional $\mathbf{\theta }_{0}$ from a data of size
$1000$ takes on average mere $2.84$ seconds, while the MAVE method
needs to spend $2432.56$ seconds on average to obtain a comparable
estimates. Hence on account of criteria (i) and (ii), our method
is indeed appealing. Applying the proposed SIP procedure to the
rive flow data of Iceland, we have obtained superior forecasts,
based on a $9$-dimensional index selected by BIC, see Figure
\ref{FIG:riverflow_fitted}.

The rest of the paper is organized as follows. Section 2 gives
details of the model specification, proposed methods of estimation
and main results. Section 3 describes the actual procedure to
implement the estimation method. Section 4 reports our findings in
an extensive simulation study. The proposed SIP model and the
estimation procedure are applied in Section 5 to the rive flow
data of Iceland. Most of the technical proofs are contained in the
Appendix.

\setcounter{chapter}{2} \renewcommand{\theproposition}{{2.%
\arabic{proposition}}} \setcounter{equation}{0} \setcounter{lemma}{0} %
\setcounter{theorem}{0} \setcounter{proposition}{0} \setcounter{corollary}{0}
\vskip .12in \noindent \textbf{2. The Method and Main Results} \label%
{SEC:method}

\vskip .10in \noindent \textbf{2.1. Identifiability and definition
of the index coefficient} \vskip .10in

It is obvious that without constraints, the SIP coefficient vector $\mathbf{%
\theta }_{0}=\left( \theta _{0,1},...,\theta _{0,d}\right) ^{T}$
is identified only up to a constant factor. Typically, one
requires that $\left\| \mathbf{\theta }_{0}\right\| =1$ which
entails that at least one of the coordinates $\theta
_{0,1},...,\theta _{0,d} $ is nonzero. One could assume without
loss of generality that $\theta _{0,d}>0$, and the candidate
$\mathbf{\theta }_{0}$ would then belong to the upper unit
hemisphere $S_{+}^{d-1}=\left\{ \left( \theta _{1},...,\theta
_{d}\right) |\sum_{p=1}^{d}\theta _{p}^{2}=1,\theta _{d}>
0\right\} $.

For a fixed $\mathbf{\theta }=\left( \theta _{1},...,\theta _{d}\right) ^{T}$%
, denote $X_{\mathbf{\theta }}=\mathbf{X}^{T}\mathbf{\theta }$, $X_{\mathbf{%
\theta },i}=\mathbf{X}_{i}^{T}\mathbf{\theta }$, $1\leq i\leq n$. Let
\begin{equation}
m_{\mathbf{\theta }}\left( X_{\mathbf{\theta }}\right) =E\left( Y|X_{\mathbf{%
\theta }}\right) =E\left\{ m\left( \mathbf{X}\right) |X_{\mathbf{\theta }%
}\right\} .  \label{DEF:mtheta}
\end{equation}
Define the risk function of $\mathbf{\theta }$ as
\begin{equation}
R\left( \mathbf{\theta }\right) =E\left[ \left\{ Y-m_{\mathbf{\theta }%
}\left( X_{\mathbf{\theta }}\right) \right\} ^{2}\right] =E\left\{ m\left(
\mathbf{X}\right) -m_{\mathbf{\theta }}\left( X_{\mathbf{\theta }}\right)
\right\} ^{2}+E\sigma ^{2}\left( \mathbf{X}\right) ,  \label{DEF:Rtheta}
\end{equation}
which is uniquely minimized at $\mathbf{\theta }_{0}$ $\in
S_{+}^{d-1}$, i.e.
\[
\mathbf{\theta }_{0}=\arg \min_{\mathbf{\theta \in
}S_{+}^{d-1}}R\left( \mathbf{\theta }\right).
\]
\vskip 0.05in \noindent \textbf{Remark 2.1.} Note that
$S_{+}^{d-1}$ is not a compact set, so we introduce a cap shape
subset of $S_{+}^{d-1}$
\[
S_{c}^{d-1}=\left\{ \left( \theta _{1},...,\theta _{d}\right)
|\sum_{p=1}^{d}\theta _{p}^{2}=1,\theta _{d}\geq
\sqrt{1-c^{2}}\right\}, c\in \left( 0,1\right)
\]
Clearly, for an appropriate choice of $c$, $\mathbf{\theta }_{0}\in
S_{c}^{d-1}$, which we assume in the rest of the paper.

Denote $\mathbf{\theta }_{-d}=\left( \theta _{1},...,\theta _{d-1}\right)
^{T}$, since for fixed $\mathbf{\theta }\in S_{+}^{d-1}$, the risk function $%
R\left( \mathbf{\theta }\right) $ depends only on the first $d-1$ values in $%
\mathbf{\theta }$, so $R\left( \mathbf{\theta }\right) $ is a function of $%
\mathbf{\theta }_{-d}$
\[
R^{*}\left( \mathbf{\theta }_{-d}\right) =R\left( \theta _{1},\theta
_{2},...,\theta _{d-1},\sqrt{1-\left\| \mathbf{\theta }_{-d}\right\| _{2}^{2}%
}\right) ,
\]
with well-defined score and Hessian matrices
\begin{equation}
S^{*}\left( \mathbf{\theta }_{-d}\right) =\frac{\partial }{\partial \mathbf{%
\theta }_{-d}}R^{*}\left( \mathbf{\theta }_{-d}\right) \text{, }H^{*}\left(
\mathbf{\theta }_{-d}\right) =\frac{\partial ^{2}}{\partial \mathbf{\theta }%
_{-d}\partial \mathbf{\theta }_{-d}^{T}}R^{*}\left( \mathbf{\theta }%
_{-d}\right). \label{DEF:SHstarmatrices}
\end{equation}

\noindent \textbf{Assumption A1:} \textit{The Hessian matrix} $H^{*}\left(
\mathbf{\theta }_{0,-d}\right) $ \textit{is positive definite and the risk
function }$R^{*}$ \textit{is locally convex at} $\mathbf{\theta }_{0,-d}$%
\textit{, i.e., for any }$\varepsilon >0$\textit{, there exists} $\delta >0$
\textit{such that} $R^{*}\left( \mathbf{\theta }_{-d}\right) -R^{*}\left(
\mathbf{\theta }_{0,-d}\right) <\delta $ \textit{implies} $\left\| \mathbf{%
\theta }_{-d}-\mathbf{\theta }_{0,-d}\right\| _{2}<\varepsilon $.

\vskip .12in \noindent \textbf{2.2. Variable transformation}
\vskip .10in

Throughout this paper, we denote by $B_{a}^{d}=\left\{ \mathbf{x}\in
R^{d}\left| \left\| \mathbf{x}\right\| \leq a\right. \right\} $ the $d$%
-dimensional ball with radius $a$ and center $\mathbf{0}$ and
\[
C^{(k)}\left( B_{a}^{d}\right) =\left\{ m\left| \text{the\ }k\text{th order
partial derivatives of }m\text{ are continuous on }B_{a}^{d}\right. \right\}
\]
the space of $k$-th order smooth functions.

\noindent \textbf{Assumption A2:} \textit{The density function of }$\mathbf{X%
}$,\textit{\ }$f\left( \mathbf{x}\right) \in C^{(4)}\left( B_{a}^{d}\right) $%
\textit{, and there are constants }$0<c_{f}\leq C_{f}$\textit{\ such that }
\[
\left\{
\begin{array}{ll}
c_{f}/\text{Vol}_{d}\left( B_{a}^{d}\right) \leq f\left( \mathbf{x}\right)
\leq C_{f}/\text{Vol}_{d}\left( B_{a}^{d}\right) , & \mathbf{x}\in B_{a}^{d}
\\
f\left( \mathbf{x}\right) \equiv 0, & \mathbf{x}\notin B_{a}^{d}
\end{array}
\right. .
\]

For a fixed $\mathbf{\theta }$, define the transformed variables of the SIP
variable $X_{\mathbf{\theta }}$
\begin{equation}
U_{\mathbf{\theta }}=F_{d}\left( X_{\mathbf{\theta }}\right) ,U_{\mathbf{%
\theta },i}=F_{d}\left( X_{\mathbf{\theta },i}\right) ,1\leq i\leq n,
\label{DEF:Utheta}
\end{equation}
in which $F_{d}$ is the a rescaled centered $\text{Beta}\left\{ \left(
d+1\right) /2,\left( d+1\right) /2\right\} $ cumulative distribution
function, i.e.
\begin{equation}
F_{d}\left( \nu \right) =\int_{-1}^{\nu /a}\frac{\Gamma \left( d+1\right) }{%
\Gamma \left\{ \left( d+1\right) /2\right\} ^{2}2^{d}}\left( 1-t^{2}\right)
^{\left( d-1\right) /2}dt,\nu \in \left[ -a,a\right] .  \label{DEF:Fd}
\end{equation}

\vskip 0.05in \noindent \textbf{Remark 2.2.} For any fixed $\mathbf{\theta }$%
, the transformed variable $U_{\mathbf{\theta }}$ in
(\ref{DEF:Utheta}) has a quasi-uniform $[0,1]$ distribution. Let
$f_{\mathbf{\theta }}\left( u\right) $ be the probability density
function of $U_{\mathbf{\theta }}$, then for any $u\in \left[
0,1\right] $%
\[
f_{\mathbf{\theta }}\left( u\right) =\left\{ F_{d}^{^{\prime
}}\left( v\right) \right\} f_{X_{\theta }}\left( v\right), \
v=F_{d}^{-1}\left( u\right),
\]
in which $f_{X_{\theta}}\left( v\right)=\lim_{\triangle \nu
\rightarrow 0} P\left( \nu \leq X_{\mathbf{\theta }}\leq \nu
+\triangle \nu \right)$. Noting that $x_{\mathbf{\theta}}$ is
exactly the projection of $\mathbf{x}$
on $\mathbf{\theta }$, let $\mathcal{D}_{\nu }=\left\{\mathbf{x}| \nu \leq x_{%
\mathbf{\theta }}\leq \nu +\triangle \nu \right\} \cap B_{a}^{d}$,
then one has
\[
P\left( \nu \leq X_{\mathbf{\theta }}\leq \nu +\triangle \nu \right)
=P\left( \mathbf{X}\in \mathcal{D}_{\nu }\right) =\int_{\mathcal{D}_{\nu
}}f\left( \mathbf{x}\right) d\mathbf{x}.
\]
According to Assumption A2
\[
\frac{c_{f}\text{Vol}_{d}(\mathcal{D}_{\nu })}{\text{Vol}_{d}\left(
B_{a}^{d}\right) }\leq P\left( \nu \leq X_{\mathbf{\theta }}\leq \nu
+\triangle \nu \right) \leq \frac{C_{f}\text{Vol}_{d}(\mathcal{D}_{\nu })}{%
\text{Vol}_{d}\left( B_{a}^{d}\right) }.
\]
On the other hand
\[
\text{Vol}_{d}(\mathcal{D}_{\nu })=\text{Vol}_{d-1}(\mathcal{J}_{\nu
})\triangle \nu +o\left( \triangle \nu \right) ,
\]
where $\mathcal{J}_{\nu
}=\left\{\mathbf{x}|x_{\mathbf{\theta}}=v\right\} \cap
B_{a}^{d}$. Note that the volume of $B_{a}^{d}$ is $%
\pi ^{d/2}a^{d}/\Gamma\left( d/2+1\right)$ and
\[
\text{Vol}_{d-1}\left( \mathcal{J}_{\nu }\right) =\left. \pi
^{\left( d-1\right) /2}\left( a^{2}-\nu ^{2}\right) ^{\left(
d-1\right) /2}\right/ \Gamma\left\{ (d+1)/{2}\right\} ,
\]
thus
\[
\frac{\text{Vol}_{d-1}(\mathcal{J}_{\nu })}{\text{Vol}_{d}\left(
B_{a}^{d}\right) }=\frac{1}{a\sqrt{\pi }}\frac{\Gamma\left(
d+1\right) }{\left\{
\Gamma\left( \frac{d+1}{2}\right) \right\} ^{2}2^{d}}\left\{ 1-\left( \frac{\nu }{%
a}\right) ^{2}\right\} ^{\left( d-1\right) /2}.
\]
Therefore $0<c_{f}\leq f_{\mathbf{\theta }}\left( u\right) \leq
C_{f}<\infty $, for any fixed $\mathbf{\theta }$ and $u\in \left[
0,1\right] $.

In terms of the transformed SIP variable $U_{\mathbf{\theta }}$ in (\ref
{DEF:Utheta}), we can rewrite the regression function $m_{\mathbf{\theta }}$
in (\ref{DEF:mtheta}) for fixed $\mathbf{\theta }$
\begin{equation}
\gamma _{\mathbf{\theta }}\left( U_{\mathbf{\theta }}\right) =E\left\{
m\left( \mathbf{X}\right) |U_{\mathbf{\theta }}\right\} =E\left\{ m\left(
\mathbf{X}\right) |X_{\mathbf{\theta }}\right\} =m_{\mathbf{\theta }}\left(
X_{\mathbf{\theta }}\right) ,  \label{EQ:mthetagtheta}
\end{equation}
then the risk function $R\left( \mathbf{\theta }\right) $ in (\ref
{DEF:Rtheta}) can be expressed as
\begin{equation}
R\left( \mathbf{\theta }\right) =E\left[ \left\{ Y-\gamma _{\mathbf{\theta }%
}\left( U_{\mathbf{\theta }}\right) \right\} ^{2}\right] =E\left\{ m\left(
\mathbf{X}\right) -\gamma _{\mathbf{\theta }}\left( U_{\mathbf{\theta }%
}\right) \right\} ^{2}+E\sigma ^{2}\left( \mathbf{X}\right) .
\label{EQ:Rthetagtheta}
\end{equation}

\vskip 0.10in \noindent \textbf{2.3. Estimation Method} \vskip
0.10in

Estimation of both $\mathbf{\theta }_{0}$ and $g$ requires a degree of
statistical smoothing, and all estimation here is carried out via cubic
spline. In the following, we define the estimator $\hat{\mathbf{\theta }}$
of $\mathbf{\theta }_{0}$ and the estimator $\hat{g}$ of $g$.

To introduce the space of splines, we pre-select an integer $n^{1/6}\ll
N=N_{n}\ll n^{1/5}\left( \log n\right) ^{-2/5}$, see Assumption A6 below.
Divide $\left[ 0,1\right] $ into $\left( N+1\right) $ subintervals $%
J_{j}=\left[ t_{j},t_{j+1}\right) $, $j=0,...,N-1,J_{N}=\left[
t_{N},1\right] $, where $T:=\left\{ t_{j}\right\} _{j=1}^{N}$ is a sequence
of equally-spaced points, called interior knots, given as
\[
t_{1-k}=...=t_{-1}=t_{0}=0<t_{1}<...<t_{N}<1=t_{N+1}=...=t_{N+k},
\]
in which $t_{j}=jh,\,\,j=0,1,...,N+1,h=1/\left( N+1\right) $ is the distance
between neighboring knots. The $j$-th B-spline of order $k$ for the knot
sequence $T$ denoted by $B_{j,k}$ is recursively defined by de Boor (2001).

Denote by $\Gamma ^{\left( k-2\right) }=\Gamma ^{\left( k-2\right) }\left[
0,1\right] $ the space of all $C^{\left( k-2\right) }\left[ 0,1\right] $
functions that are polynomials of degree $k-1$ on each interval. For fixed $%
\mathbf{\theta }$, the cubic spline estimator $\hat{\gamma}_{\mathbf{\theta }%
}$ of $\gamma _{\mathbf{\theta }}$ and the related estimator $\hat{m}_{%
\mathbf{\theta }}$ of $m_{\mathbf{\theta }}$ are defined as
\begin{equation}
\hat{\gamma}_{\mathbf{\theta }}\left( \cdot \right) =\arg \min_{\gamma
\left( \cdot \right) \in \Gamma ^{\left( 2\right) }\left[ 0,1\right]
}\sum_{i=1}^{n}\left\{ Y_{i}-\gamma \left( U_{\mathbf{\theta },i}\right)
\right\} ^{2},\mbox{\ }\hat{m}_{\mathbf{\theta }}\left( \nu \right) =\hat{%
\gamma}_{\mathbf{\theta }}\left\{ F_{d}\left( \nu \right) \right\} .
\label{DEF:mthetahat}
\end{equation}

Define the empirical risk function of $\mathbf{\theta }$
\begin{equation}
\hat{R}\left( \mathbf{\theta }\right) =n^{-1}\sum_{i=1}^{n}\left\{ Y_{i}-%
\hat{\gamma}_{\mathbf{\theta }}\left( U_{\mathbf{\theta },i}\right) \right\}
^{2}=n^{-1}\sum_{i=1}^{n}\left\{ Y_{i}-\hat{m}_{\mathbf{\theta }}\left( X_{%
\mathbf{\theta },i}\right) \right\} ^{2},  \label{DEF:Rhat}
\end{equation}
then the spline estimator of the SIP coefficient $\mathbf{\theta }_{0}$ is
defined as
\[
\hat{\mathbf{\theta }}=\arg \min_{\mathbf{\theta \in }S_{c}^{d-1}}\hat{R}%
\left( \mathbf{\theta }\right),
\]
and the cubic spline estimator of $g$ is $\hat{m}_{\mathbf{\theta
}}$ with $\mathbf{\theta }$ replaced by $\hat{\mathbf{\theta }}$,
i.e.
\begin{equation}
\hat{g}\left( \nu \right) =\left\{ \arg \min_{\gamma \left( \cdot
\right) \in \Gamma ^{\left( 2\right) }\left[ 0,1\right]
}\sum_{i=1}^{n}\left\{ Y_{i}-\gamma \left(
U_{\hat{\mathbf{\theta}},i}\right) \right\} ^{2}\right\} \left\{
F_{d}\left( \nu \right) \right\} .  \label{DEF:ghat}
\end{equation}

\vskip 0.05in \noindent \textbf{2.4. Asymptotic results} \vskip
0.1in

Before giving the main theorems, we state some other assumptions.

\noindent \textbf{Assumption A3:} \textit{The regression function }$m\in
C^{(4)}\left( B_{a}^{d}\right) $\textit{\ for some }$a>0$.

\noindent \textbf{Assumption A4:} \textit{The noise }$\varepsilon $\textit{\
satisfies }$E\left( \varepsilon \left| \mathbf{X}\right. \right) =0$, $%
E\left( \varepsilon ^{2}\left| \mathbf{X}\right. \right) =1$ \textit{and
there exists a positive constant} $M$ \textit{such that} $\sup\limits_{%
\mathbf{x}\in B^{d}}E\left( \left| \varepsilon \right| ^{3}\left| \mathbf{X}=%
\mathbf{x}\right. \right) <M$. \textit{The standard deviation function }$%
\sigma \left( \mathbf{x}\right) $\textit{\ is continuous on }$B_{a}^{d}$,
\[
0<c_{\sigma }\leq \inf_{\mathbf{x}\in B_{a}^{d}}\sigma \left( \mathbf{x}%
\right) \leq \sup_{\mathbf{x}\in B_{a}^{d}}\sigma \left( \mathbf{x}\right)
\leq C_{\sigma }<\infty .
\]

\noindent \textbf{Assumption A5:} \textit{There exist positive constants $%
K_{0}$ and }$\lambda _{0}$\textit{\ such that $\alpha \left( n\right) \leq
K_{0}e^{-\lambda _{0}n} $ holds for all $n$, with the }$\alpha $\textit{%
-mixing coefficient for }$\left\{ \mathbf{Z}_{i}=\left( \mathbf{X}%
_{i}^{T},\varepsilon _{i}\right) \right\} _{i=1}^{n}$\ \textit{\ defined as}
\[
\alpha \left( k\right) =\sup_{B\in \sigma \left\{ \mathbf{Z}_{s},s\leq
t\right\} ,C\in \sigma \left\{ \mathbf{Z}_{s},s\geq t+k\right\} }\left|
P\left( B\cap C\right) -P\left( B\right) P\left( C\right) \right| ,\mbox{\ }%
k\geq 1.
\]

\noindent \textbf{Assumption A6:} \textit{The number of interior knots} $N$%
\textit{\ satisfies:} $n^{1/6}\ll N\ll n^{1/5}\left( \log n\right) ^{-2/5}$%
\textit{.}

\vskip 0.1in \noindent \textbf{Remark 2.3.} Assumptions A3 and A4 are
typical in the nonparametric smoothing literature, see for instance,
H\"{a}rdle (1990), Fan and Gijbels (1996), Xia, Tong Li and Zhu (2002). By
the result of Pham (1986), a geometrically ergodic time series is a strongly
mixing sequence. Therefore, Assumption A5 is suitable for (\ref{sindmodel})
as a time series model under aforementioned assumptions.

We now state our main results in the next two theorems.

\begin{theorem}
\label{THM:strconsistent} Under Assumptions A1-A6, one has
\begin{equation}
\hat{\mathbf{\theta }}_{-d}\mathbf{\longrightarrow \theta }_{0,-d},a.s..
\label{EQ:strongconsist}
\end{equation}
\end{theorem}

\noindent \textbf{Proof.} Denote by $\left( \Omega ,\mathcal{F},\mathcal{P}%
\right) $ the probability space on which all $\left\{ \left( \mathbf{X}%
_{i}^{T},Y_{i}\right) \right\} _{i=1}^{\infty }$ are defined. By Proposition
\ref{PROP:unif}, given at the end of this section
\begin{equation}
\sup_{\left\| \mathbf{\theta }_{-d}\right\| _{2}\leq \sqrt{1-c^{2}}}\left|
\hat{R}^{*}\left( \mathbf{\theta }_{-d}\right) -R^{*}\left( \mathbf{\theta }%
_{-d}\right) \right| \longrightarrow 0,a.s..  \label{EQ:k0}
\end{equation}
So for any $\delta >0$ and $\omega \in \Omega $, there exists an integer $%
n_{0}\left( \omega \right) $, such that when $n>n_{0}\left( \omega \right) $%
, $\hat{R}^{*}\left( \mathbf{\theta }_{0,-d},\omega \right) -R^{*}\left(
\mathbf{\theta }_{0,-d}\right) <$ $\delta /2$. Note that $\hat{\mathbf{\theta%
}}_{-d}=\hat{\mathbf{\theta}}_{-d}\left( \omega \right) $ is the
minimizer
of $\hat{R}^{*}\left( \mathbf{\theta }_{-d},\omega \right) $, so $\hat{R}%
^{*}\left(\hat{\mathbf{\theta}}_{-d}\left( \omega \right) ,\omega
\right)
-R^{*}\left(\mathbf{\theta }_{0,-d}\right) <\delta /2$. Using (\ref{EQ:k0}%
), there exists $n_{1}\left( \omega \right) $, such that when
$n>n_{1}\left( \omega \right) $, $R^{*}\left(
\hat{\mathbf{\theta}}_{-d}\left( \omega \right) ,\omega \right)
-\hat{R}^{*}\left(\hat{\mathbf{\theta}}_{-d}\left( \omega \right)
,\omega \right) <$ $\delta /2$. Thus, when $n>\max \left(
n_{0}\left( \omega \right) ,n_{1}\left( \omega \right) \right)$,
\[
R^{*}\left( \hat{\mathbf{\theta}}_{-d}\left( \omega \right)
,\omega \right) -R^{*}\left( \mathbf{\theta }_{0,-d}\right)
<\delta /2+\hat{R}^{*}\left( \hat{\mathbf{\theta}}_{-d}\left(
\omega \right) ,\omega \right) -R^{*}\left( \mathbf{\theta
}_{0,-d}\right) <\delta /2+\delta /2=\delta .
\]
According to Assumption A1, $R^{*}$ is locally convex at
$\mathbf{\theta }_{0,-d}$, so for any $\varepsilon >0$ and any
$\omega$, if $R^{*}\left( \hat{\mathbf{\theta}}_{-d}\left( \omega
\right) ,\omega \right) -R^{*}\left(
\mathbf{\theta }_{0,-d}\right) <\delta $, then $\left\| \hat{\mathbf{\theta}}%
_{-d}\left( \omega \right) \mathbf{-\theta }_{0,-d}\right\|
<\varepsilon $ for $n$ large enough , which implies the strong
consistency.\hfill

\begin{theorem}
\label{THM:normality} Under Assumptions A1-A6, one has
\[
\sqrt{n}\left( \hat{\mathbf{\theta }}_{-d}\mathbf{-\theta }_{0,-d}\right)
\stackrel{d}{\longrightarrow }N\left\{ \mathbf{0},\Sigma \left( \mathbf{%
\theta }_{0}\right) \right\} ,
\]
where $\Sigma \left( \mathbf{\theta }_{0}\right) =\left\{ H^{*}\left(
\mathbf{\theta }_{0,-d}\right) \right\} ^{-1}\Psi \left( \mathbf{\theta }%
_{0}\right) \left\{ H^{*}\left( \mathbf{\theta }_{0,-d}\right) \right\} ^{-1}
$, $H^{*}\left( \mathbf{\theta }_{0,-d}\right) =\left\{ l_{pq}\right\}
_{p,q=1}^{d-1}$ and $\Psi \left( \mathbf{\theta }_{0}\right) =\left\{ \psi
_{pq}\right\} _{p,q=1}^{d-1}$ with
\begin{eqnarray}
l_{p,q}=-2E\left[ \left\{ \dot{\gamma}_{p}\dot{\gamma}_{q}+\gamma _{\mathbf{%
\theta }_{0}}\ddot{\gamma}_{p,q}\right\} \left( U_{\mathbf{\theta }%
_{0}}\right) \right] +2\theta _{0,q}\theta _{0,d}^{-1}E\left[ \left\{ \dot{%
\gamma}_{p}\dot{\gamma}_{d}\left( U_{\mathbf{\theta }_{0}}\right) +\gamma _{%
\mathbf{\theta }_{0}}\ddot{\gamma}_{p,d}\right\} \left( U_{\mathbf{\theta }%
_{0}}\right) \right]   \nonumber \\
+2\theta _{0,d}^{-3}E\left[ \left( \gamma _{\mathbf{\theta }_{0}}\dot{\gamma}%
_{d}\right) \left( U_{\mathbf{\theta }_{0}}\right) \right] \left\{ \left(
\theta _{0,d}^{2}+\theta _{0,p}^{2}\right) I_{\left\{ p=q\right\} }+\theta
_{0,p}\theta _{0,q}I_{\left\{ p\neq q\right\} }\right\}   \nonumber \\
+2\theta _{0,p}\theta _{0,d}^{-1}E\left[ \left\{ \dot{\gamma}_{p}\dot{\gamma}%
_{q}+\gamma _{\mathbf{\theta }_{0}}\ddot{\gamma}_{p,q}\right\} \left( U_{%
\mathbf{\theta }_{0}}\right) \right] -2\theta _{0,p}\theta _{0,q}\theta
_{0,d}^{-2}E\left[ \left\{ \dot{\gamma}_{d}^{2}+\gamma _{\mathbf{\theta }%
_{0}}\ddot{\gamma}_{d,d}\right\} \left( U_{\mathbf{\theta }_{0}}\right)
\right] ,  \nonumber
\end{eqnarray}
\begin{equation}
\psi _{pq}=4E\left[ \left\{ \left( \dot{\gamma}_{p}-\theta _{0,p}\theta
_{0,d}^{-1}\dot{\gamma}_{d}\right) \left( \dot{\gamma}_{q}-\theta
_{0,q}\theta _{0,d}^{-1}\dot{\gamma}_{d}\right) \right\} \left( U_{\mathbf{%
\theta }_{0}}\right) \left\{ \gamma _{\mathbf{\theta }_{0}}\left( U_{\mathbf{%
\theta }_{0}}\right) -Y\right\} ^{2}\right] ,  \nonumber
\end{equation}
in which $\dot{\gamma}_{p}$ and $\ddot{\gamma}_{p,q}$ are the values of $%
\frac{\partial }{\partial \theta _{p}}\gamma _{\mathbf{\theta }}$, $\frac{%
\partial ^{2}}{\partial \theta _{p}\partial \theta _{q}}\gamma _{\mathbf{%
\theta }}$ taking at $\mathbf{\theta =\theta }_{0}$, for any $p,q=1,2,...,d-1
$ and $\gamma _{\mathbf{\theta }}$ is given in (\ref{EQ:mthetagtheta}).
\end{theorem}

\noindent \textbf{Remark 2.4.} Consider the Generalized Linear Model (GLM): $%
Y=g\left( \mathbf{X}^{T}\mathbf{\theta }_{0}\right) +\sigma \left( \mathbf{X}%
\right) \varepsilon $, where $g$ is a known link function. Let $\tilde{%
\mathbf{\theta }}$ be the nonlinear least squared estimator of $\mathbf{%
\theta }_{0}$ in GLM. Theorem \ref{THM:normality} shows that under the
assumptions A1-A6, the asymptotic distribution of the $\hat{\mathbf{\theta }}%
_{-d}$ is the same as that of $\tilde{\mathbf{\theta }}$. This
implies that our proposed SIP estimator $\hat{\mathbf{\theta
}}_{-d}$ is as efficient as if the true link function $g$ is
known.

The next two propositions play an important role in our proof of
the main results. Proposition \ref{PROP:ghattheta-gtheta}
establishes the uniform convergence rate of the derivatives of
$\hat{\gamma}_{\mathbf{\theta }}$ up to order 2 to those of
$\gamma_{\mathbf{\theta }}$ in $\mathbf{\theta}$. Proposition
\ref{PROP:unif} shows that the derivatives of the risk function up
to order 2 are uniformly almost surely approximated by their
empirical versions.

\begin{proposition}
\label{PROP:ghattheta-gtheta} Under Assumptions
A2-A6, with probability $1$%
\begin{equation}
\sup\limits_{\mathbf{\theta }\in S_{c}^{d-1}}\sup_{u\in \left[
0,1\right]
}\left| \hat{\gamma}_{\mathbf{\theta }}\left( u\right) -\gamma _{\mathbf{%
\theta }}\left( u\right) \right| =O\left\{ \left( nh\right)
^{-1/2}\log n+h^{4}\right\} ,  \label{EQ:ghattheta-gtheta}
\end{equation}
\begin{equation}
\sup\limits_{1\leq p\leq d}\sup\limits_{\mathbf{\theta }\in
S_{c}^{d-1}}\max\limits_{1\leq i\leq n}\left| \frac{\partial
}{\partial
\theta _{p}}\left\{ \hat{\gamma}_{\mathbf{\theta }}\left( U_{\mathbf{\theta }%
,i}\right) -\gamma _{\mathbf{\theta }}\left( U_{\mathbf{\theta
},i}\right) \right\} \right| =O\left( \frac{\log
n}{\sqrt{nh^{3}}}+h^{3}\right) ,
\label{EQ:ghatthetaderiv-gthetaderv}
\end{equation}
\begin{equation}
\sup\limits_{1\leq p,q\leq d}\sup\limits_{\mathbf{\theta }\in
S_{c}^{d-1}}\max\limits_{1\leq i\leq n}\left| \frac{\partial
^{2}}{\partial
\theta _{p}\partial \theta _{q}}\left\{ \hat{\gamma}_{\mathbf{\theta }%
}\left( U_{\mathbf{\theta },i}\right) -\gamma _{\mathbf{\theta }}\left( U_{%
\mathbf{\theta },i}\right) \right\} \right| =O\left( \frac{\log n}{\sqrt{%
nh^{5}}}+h^{2}\right) .  \label{EQ:ghatthetaderiv2-gthetaderv2}
\end{equation}
\end{proposition}
\begin{proposition}
\label{PROP:unif} Under Assumptions A2-A6, one has for $k=0,1,2$
\[
\sup_{\left\| \mathbf{\theta }_{-d}\right\| \leq \sqrt{1-c^{2}}}\left| \frac{%
\partial ^{k}}{\partial ^{k}\mathbf{\theta }_{-d}}\left\{ \hat{R}^{*}\left(
\mathbf{\theta }_{-d}\right) -R^{*}\left( \mathbf{\theta
}_{-d}\right) \right\} \right| =o(1),a.s..
\]
\end{proposition}

Proofs of Theorem \ref{THM:normality}, Propositions
\ref{PROP:ghattheta-gtheta} and \ref{PROP:unif} are given in
Appendix.

\setcounter{chapter}{3} \label{SEC:implementation} \renewcommand{%
\thetheorem}{{\arabic{theorem}}}
\renewcommand{\thelemma}{{3.\arabic{lemma}}}
\renewcommand{\theproposition}{{3.\arabic{proposition}}} %
\setcounter{equation}{0} \setcounter{lemma}{0} \setcounter{theorem}{0}%
\setcounter{proposition}{0}\setcounter{corollary}{0}\vskip .12in
\noindent \textbf{3. Implementation} \vskip 0.10in

In this section, we will describe the actual procedure to implement the
estimation of $\mathbf{\theta }_{0}$ and $g$. We first introduce some new
notation. For fixed $\mathbf{\theta }$, write the B-spline matrix as $%
\mathbf{B}_{\mathbf{\theta }}=\left\{ B_{j,4}\left( U_{\mathbf{\theta }%
,i}\right) \right\} _{i=1,j=-3}^{n,\text{ }N}$ and
\begin{equation}
\mathbf{P}_{\mathbf{\theta }}=\mathbf{B}_{\mathbf{\theta }}\left( \mathbf{B}%
_{\mathbf{\theta }}^{T}\mathbf{B}_{\mathbf{\theta }}\right) ^{-1}\mathbf{B}_{%
\mathbf{\theta }}^{T}  \label{DEF:Ptheta}
\end{equation}
as the projection matrix onto the cubic spline space $\Gamma _{n,\mathbf{%
\theta }}^{\left( 2\right) }$. For any $p=1,...,d$, denote
\[
\mathbf{\dot{B}}_{p}=\frac{\partial }{\partial \theta _{p}}\mathbf{B}_{%
\mathbf{\theta }}, \ \mathbf{\dot{P}}_{p}=\frac{\partial }{\partial \theta _{p}}\mathbf{P}_{\mathbf{\theta }%
}.
\]
as the first order partial derivatives of $\mathbf{B}_{\mathbf{%
\theta }}$ and $\mathbf{P}_{\mathbf{\theta }}$ with respect to $\mathbf{%
\theta }$.

Let $\hat{S}^{*}(\mathbf{\theta }_{-d})$ be the score vector of
$\hat{R}^{*}\left( \mathbf{\theta }_{-d}\right)$, i.e.
\begin{equation}
\hat{S}^{*}(\mathbf{\theta }_{-d})=\frac{\partial }{\partial \mathbf{\theta }%
_{-d}}\hat{R}^{*}\left( \mathbf{\theta }_{-d}\right).
\label{DEF:Shatstarmatrices}
\end{equation}
The next lemma provides the exact forms of $\hat{S}^{*}(\mathbf{\theta }%
_{-d})$.

\begin{lemma}
\label{LEM:Shatstarmatrices} For the score vector of
$\hat{R}^{*}\left( \mathbf{\theta }_{-d}\right) $ defined in
(\ref{DEF:Shatstarmatrices}), one has
\begin{equation}
\hat{S}^{*}\left( \mathbf{\theta }_{-d}\right) =-n^{-1}\left\{ \mathbf{Y}^{T}%
\mathbf{\dot{P}}_{p}\mathbf{Y}-\theta _{p}\mathbf{\theta }_{d}^{-1}\mathbf{Y}%
^{T}\mathbf{\dot{P}}_{d}\mathbf{Y}\right\} _{p=1}^{d-1},
\label{EQ:Shatstarmatrix}
\end{equation}
where for any $p=1,2,...,d$
\begin{equation}
\mathbf{Y}^{T}\mathbf{\dot{P}}_{p}\mathbf{Y}=2\mathbf{Y}^{T}\left( \mathbf{I}%
-\mathbf{P}_{\mathbf{\theta }}\right) \mathbf{\dot{B}}_{p}\left( \mathbf{B}_{%
\mathbf{\theta }}^{T}\mathbf{B}_{\mathbf{\theta }}\right) ^{-1}\mathbf{B}_{%
\mathbf{\theta }}^{T}\mathbf{Y},  \label{EQ:YPdotY}
\end{equation}
where $\mathbf{\dot{B}}_{p}=\left\{ \left\{ B_{j,3}\left(
U_{\mathbf{\theta },i}\right) -B_{j+1,3}\left( U_{\mathbf{\theta
},i}\right) \right\} \dot{F}_{d}\left(
\mathbf{X}_{\mathbf{\theta },i}\right) h^{-1}X_{i,p}\right\} _{i=1,j=-3}^{n,%
\text{ }N}$ with
\[
\dot{F}_{d}\left( x\right)=\frac{d}{dx}F_{d}=\frac{\Gamma \left(
d+1\right) }{a\Gamma \left\{ \left( d+1\right) /2\right\} ^{2}2^{d}}\left( 1-%
\frac{x^{2}}{a^{2}}\right) ^{\frac{d-1}{2}}I\left( \left| x\right|
\leq a\right).
\]
\end{lemma}
\noindent \textbf{Proof}. For any $p=1,2,...,d$, the derivatives
of B-splines in de Boor (2001) implies
\begin{eqnarray*}
\mathbf{\dot{B}}_{p} &=&\left\{ \frac{\partial }{\partial \theta _{p}}%
B_{j,4}\left( U_{\mathbf{\theta },i}\right) \right\} _{i=1,j=-3}^{n,\text{ }%
N}=\left\{ \frac{d}{du}B_{j,4}\left( U_{\mathbf{\theta },i}\right) \frac{d}{%
d\theta _{p}}U_{\mathbf{\theta },i}\right\} _{i=1,j=-3}^{n,\text{ }N} \\
&=&3\left\{ \left\{ \frac{B_{j,3}\left( U_{\mathbf{\theta },i}\right) }{%
t_{j+3}-t_{j}}-\frac{B_{j+1,3}\left( U_{\mathbf{\theta },i}\right) }{%
t_{j+4}-t_{j+1}}\right\} \dot{F}_{d}\left( \mathbf{X}_{\mathbf{\theta }%
,i}\right) X_{i,p}\right\} _{i=1,j=-3}^{n,\text{ }N} \\
&=&\left\{ \left\{ B_{j,3}\left( U_{\mathbf{\theta },i}\right)
-B_{j+1,3}\left( U_{\mathbf{\theta },i}\right) \right\}
\dot{F}_{d}\left(
\mathbf{X}_{\mathbf{\theta },i}\right) h^{-1}X_{i,p}\right\} _{i=1,j=-3}^{n,%
\text{ }N}.
\end{eqnarray*}
Next, note that
\begin{eqnarray*}
\mathbf{\dot{P}}_{p} &=&\mathbf{\dot{B}}_{p}\left( \mathbf{B}_{\mathbf{%
\theta }}^{T}\mathbf{B}_{\mathbf{\theta }}\right) ^{-1}\mathbf{B}_{\mathbf{%
\theta }}^{T}+\mathbf{B}_{\mathbf{\theta }}\left[ \frac{\partial
}{\partial
\theta _{p}}\left\{ \left( \mathbf{B}_{\mathbf{\theta }}^{T}\mathbf{B}_{%
\mathbf{\theta }}\right) ^{-1}\mathbf{B}_{\mathbf{\theta
}}^{T}\right\}
\right] \\
&=&\mathbf{\dot{B}}_{p}\left( \mathbf{B}_{\mathbf{\theta }}^{T}\mathbf{B}_{%
\mathbf{\theta }}\right) ^{-1}\mathbf{B}_{\mathbf{\theta }}^{T}+\mathbf{B}_{%
\mathbf{\theta }}\left\{ \frac{\partial }{\partial \theta
_{p}}\left( \mathbf{B}_{\mathbf{\theta
}}^{T}\mathbf{B}_{\mathbf{\theta }}\right)
^{-1}\right\} \mathbf{B}_{\mathbf{\theta }}^{T}+\mathbf{B}_{\mathbf{\theta }%
}\left( \mathbf{B}_{\mathbf{\theta }}^{T}\mathbf{B}_{\mathbf{\theta }%
}\right) ^{-1}\mathbf{\dot{B}}_{p}^{T}.
\end{eqnarray*}
Since
\[
0\equiv \frac{\partial \left\{ \left( \mathbf{B}_{\mathbf{\theta }}^{T}%
\mathbf{B}_{\mathbf{\theta }}\right) ^{-1}\mathbf{B}_{\mathbf{\theta }}^{T}%
\mathbf{B}_{\mathbf{\theta }}\right\} }{\partial \theta
_{p}}=\frac{\partial \left( \mathbf{B}_{\mathbf{\theta
}}^{T}\mathbf{B}_{\mathbf{\theta }}\right)
^{-1}}{\partial \theta _{p}}\mathbf{B}_{\mathbf{\theta }}^{T}\mathbf{B}_{%
\mathbf{\theta }}+\left( \mathbf{B}_{\mathbf{\theta }}^{T}\mathbf{B}_{%
\mathbf{\theta }}\right) ^{-1}\frac{\partial \left( \mathbf{B}_{\mathbf{%
\theta }}^{T}\mathbf{B}_{\mathbf{\theta }}\right) }{\partial
\theta _{p}},
\]
and $\frac{\partial }{\partial \theta _{p}}\left(
\mathbf{B}_{\mathbf{\theta
}}^{T}\mathbf{B}_{\mathbf{\theta }}\right) =\mathbf{\dot{B}}_{p}^{T}\mathbf{B%
}_{\mathbf{\theta }}+\mathbf{B}_{\mathbf{\theta
}}^{T}\mathbf{\dot{B}}_{p}$, thus
\[
\frac{\partial }{\partial \theta _{p}}\left( \mathbf{B}_{\mathbf{\theta }%
}^{T}\mathbf{B}_{\mathbf{\theta }}\right) ^{-1}=-\left( \mathbf{B}_{\mathbf{%
\theta }}^{T}\mathbf{B}_{\mathbf{\theta }}\right) ^{-1}\left( \mathbf{\dot{B}%
}_{p}^{T}\mathbf{B}_{\mathbf{\theta }}+\mathbf{B}_{\mathbf{\theta }}^{T}%
\mathbf{\dot{B}}_{p}\right) \left( \mathbf{B}_{\mathbf{\theta }}^{T}\mathbf{B%
}_{\mathbf{\theta }}\right) ^{-1}.
\]
Hence
\[
\mathbf{\dot{P}}_{p}=\left( \mathbf{\mathbf{I}}-\mathbf{P}_{\mathbf{\theta }%
}\right) \mathbf{\dot{B}}_{p}\left( \mathbf{B}_{\mathbf{\theta }}^{T}\mathbf{%
B}_{\mathbf{\theta }}\right) ^{-1}\mathbf{B}_{\mathbf{\theta }}^{T}+\mathbf{B%
}_{\mathbf{\theta }}\left( \mathbf{B}_{\mathbf{\theta }}^{T}\mathbf{B}_{%
\mathbf{\theta }}\right) ^{-1}\mathbf{\dot{B}}_{p}^{T}\left( \mathbf{\mathbf{%
I}}-\mathbf{P}_{\mathbf{\theta }}\right).
\]
Thus, (\ref{EQ:YPdotY}) follows immediately. \hfill

In practice, the estimation is implemented via the following procedure.

Step 1. \textit{Standardize the predictor vectors $\left\{\mathbf{X}%
_{i}\right\}_{i=1}^{n}$ and for each fixed $\mathbf{\theta}\in
S_{c}^{d-1}$ obtain the CDF transformed variables $\left\{U_{
\mathbf{\theta},i}\right\}_{i=1}^{n}$ of the SIP variable $\left\{X_{\mathbf{%
\theta},i}\right\}_{i=1}^{n}$ through formula (\ref{DEF:Fd}), where the
radius $a$ is taken to be the 95\% percentile of $\left\{\|\mathbf{X}%
_{i}\|\right\}_{i=1}^{n}$}.

Step 2. \textit{Compute  quadratic and cubic B-spline basis at
each value $U_{\mathbf{\theta},i}$, where the number of interior
knots $N$ is}
\begin{equation}
N=\min\left\{c_{1}\left[n^{1/5.5}\right], c_{2}\right\},
\label{EQ:Numberofknots}
\end{equation}

Step 3. \textit{Find the estimator $\hat{\mathbf{\theta}}$ of $\mathbf{\theta}_{0}$ by minimizing $%
\hat{R}^{*}$ through the port optimization routine with
$\left(0,0,...,1\right)^{T}$ as the initial value and the
empirical score vector $\hat{S}^{*}$ in (\ref{EQ:Shatstarmatrix}).
If $d<n$, one can
take the simple LSE (without the intercept) for data $\left\{Y_i,\mathbf{X}%
_{i}\right\}_{i=1}^{n}$ with its last coordinate set positive.}

Step 4. \textit{Obtain the spline estimator $\hat{g}$ of $g$ by
plugging $\hat{\mathbf{\theta}}$ obtained in Step 3 into
(\ref{DEF:ghat}).}

\vskip 0.1in \noindent \textbf{Remark 3.1.} In (\ref{EQ:Numberofknots}), $%
c_{1}$ and $c_{2}$ are positive integers and $[\nu]$ denotes the integer
part of $\nu$. The choice of the tuning parameter $c_1$ makes little
difference for a large sample and according to our asymptotic theory there
is no optimal way to set these constants. We recommend using $%
c_{1}=1$ to save computing for massive data sets. The first term ensures
Assumption A6. The addition constrain $c_{2}$ can be taken from 5 to 10 for
smooth monotonic or smooth unimodel regression and $c_{2}>10$ if has many
local minima and maxima, which is very unlikely in application.

\setcounter{chapter}{4}\setcounter{equation}{0} \renewcommand{\thetable}{{%
\arabic{table}}} \setcounter{figure}{0} \vskip .12in \noindent
\textbf{4. Simulations} \vskip 0.10in

In this section, we carry out two simulations to illustrate the
finite-sample behavior of our SIP estimation method. The number of interior
knots $N$ is computed according to (\ref{EQ:Numberofknots}) with $c_1=1,
c_2=5$. All of our codes have been written in R.

\vskip 0.1in \noindent \textbf{Example 1.} Consider the model in Xia, Li,
Tong and Zhang (2004)
\[
Y=m\left( \mathbf{X}\right) +\sigma _{0}\varepsilon,\ \sigma _{0}=0.3,0.5,\
\varepsilon \stackrel{i.i.d}{\sim}N(0,1)
\]
where $\mathbf{X}=\left( X_{1},X_{2}\right)^{T}{\sim}N(\mathbf{0},I_2)$,
truncated by $[-2.5,2.5]^2$ and
\begin{equation}
m\left( \mathbf{x}\right) =x_{1}+x_{2}+4\exp \left\{ -\left(
x_{1}+x_{2}\right) ^{2}\right\} +\delta \left( x_{1}^{2}+x_{2}^{2}\right)
^{1/2}.  \label{MODEL:Xia}
\end{equation}
If $\delta =0$, then the underlying true function $m$ is a single-index
function, i.e., $m\left( \mathbf{X}\right) =\sqrt{2}\mathbf{X}^{T}\mathbf{%
\theta } _{0}+4\exp \left\{ -2\left( \mathbf{X}^{T}\mathbf{\theta }%
_{0}\right) ^{2}\right\}$, where $\mathbf{\theta }_{0}^{T}=\left(
1,1\right) /\sqrt{2}$. While $\delta\neq 0$, then $m$ is not a
genuine single-index function. An impression of the bivariate
function $m$ for $\delta=0$ and $\delta=1$ can be gained in Figure
\ref{FIG:Xia2004} (a) and (b), respectively.

\setlength{\unitlength}{1cm}
\begin{table}[hb]
\caption{Report of Example 1 (Values out/in parentheses: $\delta=0$/$\delta=1
$)}
\label{TAB:meansdxia}\centering
\fbox{%
\begin{tabular}{c|c|c|ccc|c}
\hline
$\sigma_0$ & $n$ & $\mathbf{\theta}_{0}$ & \textbf{BIAS} & \textbf{SD} &
\textbf{MSE} & \textbf{Average MSE} \\ \hline
\multirow{8}{*}{$0.3$} & \multirow{4}{*}{$100$} & \multirow{2}{*}{$\mathbf{%
\theta}_{0,1}$} & $5e-04$ & $0.00825$ & $7e-05$ &  \\
&  &  & $(-0.00236)$ & $(0.02093)$ & $(0.00044)$ & $7e-05$ \\
&  & \multirow{2}{*}{$\mathbf{\theta}_{0,2}$} & $-6e-04$ & $0.00826$ & $7e-05
$ & $(0.00043)$ \\
&  &  & $(0.00174)$ & $(0.02083)$ & $(0.00043)$ &  \\ \cline{2-7}
& \multirow{4}{*}{$300$} & \multirow{2}{*}{$\mathbf{\theta}_{0,1}$} & $%
-0.00124$ & $0.00383$ & $2e-05$ &  \\
&  &  & $(-0.00129)$ & $(0.01172)$ & $(0.00014)$ & $2e-05$ \\
&  & \multirow{2}{*}{$\mathbf{\theta}_{0,2}$} & $-0.00124$ & $0.00383$ & $%
2e-05$ & $(0.00014)$ \\
&  &  & $(0.00110)$ & $(0.01160)$ & $(0.00013)$ &  \\ \hline
\multirow{8}{*}{$0.5$} & \multirow{4}{*}{$100$} & \multirow{2}{*}{$\mathbf{%
\theta}_{0,1}$} & $0.00121$ & $0.01346$ & $0.00018$ &  \\
&  &  & $(-0.00137 )$ & $(0.02257 )$ & $(0.00051)$ & $0.00018$ \\
&  & \multirow{2}{*}{$\mathbf{\theta}_{0,2}$} & $-0.00147$ & $0.01349$ & $%
0.00018$ & $(0.00051)$ \\
&  &  & $(0.00062)$ & $(0.02309)$ & $(0.00052)$ &  \\ \cline{2-7}
& \multirow{4}{*}{$300$} & \multirow{2}{*}{$\mathbf{\theta}_{0,1}$} & $%
-0.00204$ & $0.00639$ & $4e-05$ &  \\
&  &  & $(-0.00229 )$ & $(0.01205)$ & $(0.00015)$ & $4e-05$ \\
&  & \multirow{2}{*}{$\mathbf{\theta}_{0,2}$} & $0.00197 $ & $0.00637 $ & $%
4e-05$ & $(0.00015)$ \\
&  &  & $(0.00208)$ & $(0.01190)$ & $(0.00014)$ &  \\ \hline
\end{tabular}
}
\end{table}

For $\delta=0,1$, we draw $100$ random realizations of each sample size $%
n=50,100,300$ respectively. To demonstrate how close our SIP estimator is to
the true index parameter $\mathbf{\theta}_{0}$, Table \ref{TAB:meansdxia}
lists the sample mean (MEAN), bias (BIAS), standard deviation (SD), the mean
squared error (MSE) of the estimates of $\mathbf{\theta }_{0}$ and the
average MSE of both directions. From this table, we find that the SIP
estimators are very accurate for both cases $\delta=0$ and $\delta=1$, which
shows that our proposed method is robust against the deviation from
single-index model. As we expected, when the sample size increases, the SIP
coefficient is more accurately estimated. Moreover, for $n=100,300$, the
total average is inversely proportional to $n$.

\textbf{Example 2.} Consider the heteroscedastic regression model (\ref
{sindmodel}) with
\begin{equation}
m\left( \mathbf{X}\right) =\sin \left( \frac{\pi}{4} \mathbf{X}^{T}\mathbf{%
\theta }_{0}\right),\ \sigma \left( \mathbf{X}\right) =\sigma _{0}\frac{%
\left\{ 5-\exp \left( \left. \left\| \mathbf{X}\right\|\right/ \sqrt{d}%
\right) \right\} }{5+\exp \left( \left. \left\| \mathbf{X}\right\| \right/
\sqrt{d}\right) },  \label{MODEL:simulation}
\end{equation}
in which $\mathbf{X}_{i}=\left\{ X_{i,1},...,X_{i,d}\right\}^{T}$ and $%
\varepsilon_i$, $i=1,...,n$, are $\stackrel{i.i.d}{\sim }$ $N\left(
0,1\right) $, $\sigma _{0}=0.2$. In our simulation, the true parameter $%
\mathbf{\theta }_{0}^{T}=\left. \left( 1,1,0,...,0,1\right)
\right/ \sqrt{3}$ for different sample size $n$ and dimension $d$.
The superior performance of SIP estimators is borne out in
comparison with MAVE of Xia, Tong, Li and Zhu (2002). We also
investigate the behavior of SIP estimators in the previously
unemployed cases that sample size $n$ is smaller than or equal to
$d$, for instance, $n=100,d=100,200$ and $n=200,d=200,400$. The
average MSEs of the $d$ dimensions are listed in Table
\ref{TAB:mse-time}, from which we see that the performance of the
SIP estimators are quite reasonable and in most of the scenarios
$n\leq d$, the SIP estimators still work astonishingly well where
the MAVEs become unreliable. For $n=100$, $d=10,50,100,200$, the
estimates of the link prediction function $g$ from model
(\ref{MODEL:simulation}) are plotted in Figure
\ref{FIG:estimation}, which is rather satisfactory even when
dimension $d$ exceeds the sample size $n$.

% This is the table of computing time
\setlength{\unitlength}{1cm}
\begin{table}[ht]
\caption{Report of Example 2}
\label{TAB:mse-time}\centering
\fbox{%
\begin{tabular}{|c|c|r|r||r|r|}
\hline
\multirow{2}{*}{\textbf{Sample Size} $n$} & \multirow{2}{*}{%
\textbf{Dimension} $d$} & \multicolumn{2}{c||}{\textbf{Average MSE}} &
\multicolumn{2}{c|}{\textbf{Time}} \\ \cline{3-6}
&  & \textbf{MAVE} & \textbf{SIP} & \textbf{MAVE} & \textbf{SIP} \\
\hline\hline
\multirow{6}{*}{$50$} & $4$ & $0.00020$ & $0.00018$ & $1.91$ & $0.19$ \\
& $10$ & $0.00031$ & $0.00043$ & $2.17$ & $0.10$ \\
& $30$ & $0.00106$ & $0.00285$ & $2.77$ & $0.13$ \\
& $50$ & $0.00031$ & $0.00043$ & $3.29$ & $0.10$ \\
& $100$ & $0.00681$ & $0.00620$ & $5.94$ & $0.31$ \\
& $200$ & $0.00529$ & $0.00407$ & $27.90$ & $0.49$ \\ \hline
\multirow{6}{*}{$100$} & $4$ & $0.00008$ & $0.00008$ & $3.28$ & $0.09$ \\
& $10$ & $0.00012$ & $0.00017$ & $3.93$ & $0.13$ \\
& $30$ & $0.00017$ & $0.00058$ & $5.41$ & $0.15$ \\
& $50$ & $0.00032$ & $0.00127$ & $8.48$ & $0.16$ \\
& $100$ & --- & $0.00395$ & --- & $0.44$ \\
& $200$ & --- & $0.00324$ & --- & $0.73$ \\ \hline
\multirow{6}{*}{$200$} & $4$ & $0.00004$ & $0.00003$ & $5.32$ & $0.17$ \\
& $10$ & $0.00005$ & $0.00007$ & $7.49$ & $0.24$ \\
& $30$ & $0.00006$ & $0.00017$ & $10.08$ & $0.26$ \\
& $50$ & $0.00007$ & $0.00030$ & $15.42$ & $0.24$ \\
& $100$ & $0.00015$ & $0.00061$ & $40.81$ & $0.54$ \\
& $200$ & --- & $0.00197$ & --- & $1.44$ \\ \hline
\multirow{7}{*}{$500$} & $4$ & $0.00002$ & $0.00001$ & $14.44$ &
$0.76$
\\
& $10$ & $0.00002$ & $0.00003$ & $24.54$ & $0.79$ \\
& $30$ & $0.00002$ & $0.00008$ & $32.51$ & $0.83$ \\
& $50$ & $0.00002$ & $0.00010$ & $52.93$ & $0.89$ \\
& $100$ & $0.00003$ & $0.00012$ & $143.07$ & $0.99$ \\
& $200$ & $0.00004$ & $0.00020$ & $386.80$ & $1.96$ \\
& $400$ & --- & $0.00054$ & --- & $4.98$ \\ \hline
\multirow{7}{*}{$1000$} & $4$ & $0.00001$ & $0.00001$ & $33.57$ &
$1.95$
\\
& $10$ & $0.00001$ & $0.00001$ & $62.54$ & $3.64$ \\
& $30$ & $0.00001$ & $0.00002$ & $92.41$ & $1.95$ \\
& $50$ & $0.00001$ & $0.00003$ & $155.38$ & $2.72$ \\
& $100$ & $0.00001$ & $0.00005$ & $275.73$ & $1.81$ \\
& $200$ & $0.00008$ & $0.00006$ & $2432.56$ & $2.84$ \\
& $400$ & --- & $0.00010$ & --- & $9.35$ \\ \hline
\end{tabular}
}
\end{table}

Theorem \ref{THM:strconsistent} indicates that $\hat{\mathbf{\theta}}_{-d}$
is strongly consistent of $\mathbf{\theta }_{0,-d}$. To see the convergence,
we run $100$ replications and in each replication, the value of $\Vert \hat{%
\mathbf{\theta}}-\mathbf{\theta }_{0}\Vert/\sqrt{d}$ is computed. Figure \ref
{FIG:density} plots the kernel density estimations of the $100$ $\Vert \hat{%
\mathbf{\theta}}-\mathbf{\theta }_{0}\Vert$ in Example 2, in which
dimension $d=10,50,100,200$. There are four types of line
characteristics which correspond to the two sample sizes, the
dotted-dashed line ($n=100$), dotted line ($n=200$), dashed line
($500$) and solid line ($n=1000$). As sample sizes increasing, the
squared errors are becoming closer to $0$, with narrower spread
out, confirmative to the conclusions of Theorem \ref
{THM:strconsistent}.

Lastly, we report the average computing time of Example 2 to
generate one sample of size $n$ and perform the SIP or MAVE
procedure done on the same ordinary Pentium IV PC in Table
\ref{TAB:mse-time}. From Table \ref{TAB:mse-time}, one sees that
our proposed SIP estimator is much faster than the MAVE. The
computing time for MAVE is extremely sensitive to sample size as
we expected. For very large $d$, MAVE becomes unstable to the
point of the breaking down in four cases.

\setcounter{chapter}{5}\setcounter{equation}{0} \renewcommand{\thetable}{{%
\arabic{table}}} \setcounter{figure}{0} \vskip .12in \noindent
\textbf{5. An application} \vskip 0.10in

In this section we demonstrate the proposed SIP model through the
river flow data of J\"{o}kuls\'{a} Eystri River of Iceland, from
January 1, 1972 to December 31, 1974. There are 1096 observations,
see Tong (1990). The response variables are the daily river flow
($Y_t$), measured in meter cubed per second of J\"{o}kuls\'{a}
Eystri River. The exogenous variables are temperature ($X_t$) in
degrees Celsius and daily precipitation ($Z_t$) in millimeters
collected at the meteorological station at Hveravellir.

This data set was analyzed earlier through threshold
autoregressive (TAR) models by Tong, Thanoon and Gudmundsson
(1985), Tong (1990), and nonlinear additive autoregressive
(NAAR$X$) models by Chen and Tsay (1993). Figure
\ref{FIG:riverflow} shows the plots of the three time series, from
which some nonlinear and non-stationary features of the river flow
series are evident. To make these series stationary, we remove the
trend by a simple quadratic spline regression and these trends
(dashed lines) are shown in Figure \ref{FIG:riverflow}. By an
abuse of notation, we shall continue to use $X_t$, $Y_t$, $Z_t$ to
denote the detrended series.

In the analysis, we pre-select all the lagged values in the last 7
days (1 week), i.e., the predictor pool is $\left\{Y_{t-1},...,
Y_{t-7}, X_{t}, X_{t-1},..., X_{t-7}, Z_{t}, Z_{t-1},...,
Z_{t-7},\right\}$. Using BIC similar to Huang and Yang (2004) for
our proposed spline SIP model with 3 interior knots, the following
$9$ explanatory variables are selected from the above set
$\left\{Y_{t-1},..., Y_{t-4}, X_{t}, X_{t-1}, X_{t-2}, Z_{t},
Z_{t-1}\right\}$. Based on this selection, we fit the SIP model
again and obtain the estimate of the SIP coefficient $
\hat{\mathbf{\theta}}=\left\{-0.877, 0.382, -0.208, 0.125, -0.046,
-0.034, 0.004, -0.126, 0.079\right\}^{T}$. Figure
\ref{FIG:riverflow_fitted} (a) and (b) display the fitted river
flow series and the residuals against time.

Next we examine the forecasting performance of the SIP method. We
start with estimating the SIP estimator using only observations of
the first two years, then we perform the out-of-sample rolling
forecast of the entire third year. The observed values of the
exogenous variables are used in the forecast. Figure
\ref{FIG:riverflow_fitted} (c) shows this SIP out-of-sample
rolling forecasts. For the purpose of comparison, we also try the
MAVE method, in which the same predictor vector is selected by
using BIC. The mean squared prediction error is $60.52$ for the
SIP model, $61.25$ for MAVE, $65.62$ for NAAR$X$, $66.67$ for TAR
and $81.99$ for the linear regression model, see Chen and Tsay
(1993). Among the above five models, the SIP model produces the
best forecasts.

\setcounter{chapter}{6} \setcounter{equation}{0} \vskip .10in
\noindent \textbf{6. Conclusion} \vskip 0.05in

In this paper we propose a robust SIP model for stochastic
regression under weak dependence regardless if the underlying
function is exactly a single-index function. The proposed spline
estimator of the index coefficient possesses not only the usual
strong consistency and $\sqrt{n}$-rate asymptotically normal
distribution, but also is as efficient as if the true link
function $g$ is known. By taking advantage of the spline smoothing
method and the iterative method, the proposed procedure is much
faster than the MAVE method. This procedure is especially powerful
for large sample size $n$ and high dimension $d$ and unlike the
MAVE method, the performance of the SIP remains satisfying in the
case $d>n$.

\vskip 0.05in \noindent {\large \textbf{Acknowledgment}}

This work is part of the first author's dissertation under the
supervision of the second author, and has been supported in part
by NSF award DMS 0405330.

\setcounter{chapter}{8} \renewcommand{\thetheorem}{{A.\arabic{theorem}}}%
\renewcommand{\theproposition}{{ A.\arabic{proposition}}} %
\renewcommand{\thelemma}{{A.\arabic{lemma}}} \renewcommand{%
\thecorollary}{{A.\arabic{corollary}}} \renewcommand{\theequation}{A.%
\arabic{equation}} \renewcommand{\thesubsection}{A.\arabic{subsection}} %
\setcounter{equation}{0} \setcounter{lemma}{0} \setcounter{proposition}{0}%
\setcounter{theorem}{0} \setcounter{subsection}{0}\setcounter{corollary}{0}%
\vskip 0.10in \noindent \textbf{Appendix}

\vskip .05in \noindent \textbf{A.1. Preliminaries}

In this section, we introduce some properties of the B-spline.
\begin{lemma}
\label{LEM:normequiv} There exist constants $c>0$ such that for $%
\sum_{j=-k+1}^{N}\alpha _{j,k}B_{j,k}$ up to order $k=4$
\[
\left\{
\begin{array}{cc}
ch^{1/r}\left\| \mathbf{\alpha }\right\| _{r}\leq \left\|
\sum_{k=2}^{4}\sum_{j=-k+1}^{N}\alpha _{j,k}B_{j,k}\right\| _{r}\leq \left(
3^{r-1}h\right) ^{1/r}\left\| \mathbf{\alpha }\right\| _{r}, & 1\leq r\leq
\infty  \\
ch^{1/r}\left\| \mathbf{\alpha }\right\| _{r}\leq \left\|
\sum_{k=2}^{4}\sum_{j=-k+1}^{N}\alpha _{j,k}B_{j,k}\right\| _{r}\leq \left(
3h\right) ^{1/r}\left\| \mathbf{\alpha }\right\| _{r}, & 0<r<1
\end{array}
\right. ,
\]
where $\mathbf{\alpha :=}\left( \alpha _{-1,2},\alpha
_{0,2},...,\alpha _{N,2},...,\alpha _{N,4}\right) $. In
particular, under Assumption A2, for any fixed $\mathbf{\theta }$
\[
ch^{1/2}\left\| \mathbf{\alpha }\right\| _{2}\leq \left\|
\sum_{k=2}^{4}\sum_{j=-k+1}^{N}\alpha _{j,k}B_{j,k}\right\| _{2,\mathbf{%
\theta }}\leq Ch^{1/2}\left\| \mathbf{\alpha }\right\| _{2}.
\]
\end{lemma}

\noindent \textbf{Proof.} It follows from the B-spline property on
page 96
of de Boor (2001), $\sum_{k=2}^{4}\sum_{j=-k+1}^{N}B_{j,k}\equiv 3$ on $%
\left[ 0,1\right] $. So the right inequality follows immediate for
$r=\infty $. When $1\leq r<\infty $, we use H\"{o}lder's
inequality to find
\begin{eqnarray*}
\left| \sum_{k=2}^{4}\sum_{j=-k+1}^{N}\alpha _{j,k}B_{j,k}\right| &\leq
&\left( \sum_{k=2}^{4}\sum_{j=-k+1}^{N}\left| \alpha _{j,k}\right|
^{r}B_{j,k}\right) ^{1/r}\left(
\sum_{k=2}^{4}\sum_{j=-k+1}^{N}B_{j,k}\right) ^{1-1/r} \\
&=&3^{1-1/r}\left( \sum_{k=2}^{4}\sum_{j=-k+1}^{N}\left| \alpha
_{j,k}\right| ^{r}B_{j,k}\right) ^{1/r}.
\end{eqnarray*}
Since all the knots are equally spaced, $\int_{-\infty }^{\infty
}B_{j,k}\left( u\right) du\leq h$, the right inequality follows from
\[
\int_{0}^{1}\left| \sum_{k=2}^{4}\sum_{j=-k+1}^{N}\alpha _{j,k}B_{j,k}\left(
u\right) \right| ^{r}du\leq 3^{r-1}h\left\| \mathbf{\alpha }\right\|
_{r}^{r}.
\]
When $r<1$, we have
\[
\left| \sum_{k=2}^{4}\sum_{j=-k+1}^{N}\alpha _{j,k}B_{j,k}\right| ^{r}\leq
\sum_{k=2}^{4}\sum_{j=-k+1}^{N}\left| \alpha _{j,k}\right| ^{r}B_{j,k}^{r}.
\]
Since $\int_{-\infty }^{\infty }B_{j,k}^{r}\left( u\right) du\leq
t_{j+k}-t_{j}=kh$ and
\[
\int_{0}^{1}\left| \sum_{k=2}^{4}\sum_{j=-k+1}^{N}\alpha _{j,k}B_{j,k}\left(
u\right) \right| ^{r}du\leq \left\| \mathbf{\alpha }\right\|
_{r}^{r}\int_{-\infty }^{\infty }B_{j,k}^{r}\left( u\right) du\leq 3h\left\|
\mathbf{\alpha }\right\| _{r}^{r},
\]
the right inequality follows in this case as well. For the left
inequalities, we derive from Theorem 5.4.2, DeVore and Lorentz (1993)
\[
\left| \alpha _{j,k}\right| \leq C_{1}h^{-1/r}\int_{t_{j}}^{t_{j+1}}\left|
\sum_{j=-k+1}^{N}\alpha _{j,k}B_{j,k}\left( u\right) \right| ^{r}du
\]
for any $0<r\leq \infty ,$ so
\[
\left| \alpha _{j,k}\right| ^{r}\leq
C_{1}^{r}h^{-1}\int_{t_{j}}^{t_{j+1}}\left| \sum_{j=-k+1}^{N}\alpha
_{j,k}B_{j,k}\left( u\right) \right| ^{r}du.
\]
Since each $u\in \left[ 0,1\right] $ appears in at most $k$ intervals $%
\left( t_{j,}t_{j+k}\right) $, adding up these inequalities, we obtain that
\[
\left\| \mathbf{\alpha }\right\| _{r}^{r}\leq
C_{1}h^{-1}\sum_{k=1}^{4}\int_{t_{j}}^{t_{j+k}}\left|
\sum_{j=-k+1}^{N}\alpha _{j,k}B_{j,k}\left( u\right) \right| ^{r}du\leq
3Ch^{-1}\left\| \sum_{j=-k+1}^{N}\alpha _{j,k}B_{j,k}\right\| _{r}^{r}.
\]
The left inequality follows. \hfill

For any functions $\phi$ and $\varphi$, define the empirical inner
product and the empirical norm as
\[
\left\langle \phi ,\varphi \right\rangle _{\mathbf{\theta }}
=\int_{0}^{1}\phi \left( u\right) \varphi \left( u\right)
f_{\mathbf{\theta }}\left( u\right) du,\  \left\| \phi \right\|
_{2,n,\mathbf{\theta }}^{2}=n^{-1}\sum_{i=1}^{n}\phi ^{2}\left(
U_{\mathbf{\theta },i}\right).
\]
In addition, if functions $\phi ,\varphi$ are $L_{2}\left[ 0,1\right]$%
-integrable, define the theoretical inner product and its
corresponding theoretical $L_2$ norm as
\[
\left\| \phi \right\| _{2,\mathbf{\theta }}^{2} =\int_{0}^{1}\phi
^{2}\left( u\right) f_{\mathbf{\theta }}\left( u\right) du,\
\left\langle \phi ,\varphi \right\rangle _{n,\mathbf{\theta }}
=n^{-1}\sum_{i=1}^{n}\phi \left( U_{\mathbf{\theta ,}i}\right)
\varphi \left( U_{\mathbf{\theta ,}i}\right).
\]

\begin{lemma}
\label{LEM:th-em product} Under Assumptions A2, A5 and A6, with
probability $1$,
\[
\sup\limits_{\mathbf{\theta }\in S_{c}^{d-1}}\max_{\substack{ k,k^{\prime
}=2,3,4 \\ 1\leq j, j^{\prime }\leq N }}\left| \left\langle
B_{j,k},B_{j^{\prime },k^{\prime }}\right\rangle _{n,\mathbf{\theta }%
}-\left\langle B_{j,k},B_{j^{\prime },k^{\prime }}\right\rangle _{\mathbf{%
\theta }}\right| =O\left\{ \left( nN\right) ^{-1/2}\log n\right\} .
\]
\end{lemma}

\noindent \textbf{Proof.} We only prove the case $k=k^{\prime }=4$, all
other cases are similar. Let
\[
\zeta _{\mathbf{\theta ,}j,j^{\prime },i}=B_{j,4}\left( U_{\mathbf{\theta }%
,i}\right) B_{j^{\prime },4}\left( U_{\mathbf{\theta },i}\right)
-EB_{j,4}\left( U_{\mathbf{\theta },i}\right) B_{j^{\prime },4}\left( U_{%
\mathbf{\theta },i}\right) ,
\]
with the second moment
\[
E\zeta _{\mathbf{\theta ,}j,j^{\prime },i}^{2}=E\left[B_{j,4}^{2}\left( U_{%
\mathbf{\theta },i}\right) B_{j^{\prime },4}^{2}\left( U_{\mathbf{\theta }%
,i}\right)\right] -\left\{ EB_{j,4}\left( U_{\mathbf{\theta },i}\right)
B_{j^{\prime },4}\left( U_{\mathbf{\theta },i}\right) \right\} ^{2},
\]
where $\left\{ EB_{j,4}\left( U_{\mathbf{\theta },i}\right) B_{j^{\prime
},4}\left( U_{\mathbf{\theta },i}\right) \right\} ^{2}\sim N^{-2}$, $%
E\left[B_{j,4}^{2}\left( U_{\mathbf{\theta },i}\right) B_{j^{\prime
},4}^{2}\left( U_{\mathbf{\theta },i}\right)\right] \sim N^{-1}$ by
Assumption A2. Hence, $E\zeta _{\mathbf{\theta ,}j,j^{\prime },i}^{2}\sim
N^{-1}$. The $k$-th moment is given by
\begin{eqnarray*}
&&\left. E\left| \zeta _{\mathbf{\theta ,}j,j^{\prime },i}\right|
^{k}=E\left| B_{j,4}\left( U_{\mathbf{\theta },i}\right) B_{j^{\prime
},4}\left( U_{\mathbf{\theta },i}\right) -EB_{j,4}\left( U_{\mathbf{\theta }%
,i}\right) B_{j^{\prime },4}\left( U_{\mathbf{\theta },i}\right) \right|
^{k}\right. \\
&\leq &2^{k-1}\left\{ E\left| B_{j,4}\left( U_{\mathbf{\theta },i}\right)
B_{j^{\prime },4}\left( U_{\mathbf{\theta },i}\right) \right| ^{k}+\left|
EB_{j,4}\left( U_{\mathbf{\theta },i}\right) B_{j^{\prime },4}\left( U_{%
\mathbf{\theta },i}\right) \right| ^{k}\right\},
\end{eqnarray*}
where $\left| EB_{j,4}\left( U_{\mathbf{\theta },i}\right) B_{j^{\prime
},4}\left( U_{\mathbf{\theta },i}\right) \right| ^{k}\sim N^{-k}$, $E\left|
EB_{j,4}\left( U_{\mathbf{\theta },i}\right) B_{j^{\prime },4}\left( U_{%
\mathbf{\theta },i}\right) \right| ^{k}\sim N^{-1}$. Thus,
there exists a constant $C>0$ such that $E\left| \zeta _{\mathbf{\theta ,}%
j,j^{\prime },i}\right| ^{k}\leq C2^{k-1}k!E\zeta _{j,j^{\prime },i}^{2}$.
So the Cram\'{e}r's condition is satisfied with Cram\'{e}r's constant $c^{*}$%
. By the Bernstein's inequality (see Bosq (1998), Theorem 1.4, page 31), we
have for $k=3$
\[
P\left\{ \left| n^{-1}\sum_{i=1}^{n}\zeta _{\mathbf{\theta ,}j,j^{\prime
},i}\right| \geq \delta _{n}\right\} \leq a_{1}\exp \left( -\frac{q\delta
_{n}^{2}}{25m_{2}^{2}+5c^{*}\delta _{n}}\right) +a_{2}\left( k\right) \alpha
\left( \left[ \frac{n}{q+1}\right] \right) ^{6/7},
\]
where
\[
\delta _{n}=\delta \frac{\log n}{\sqrt{nN}},\text{ }a_{1}=2\frac{n}{q}%
+2\left( 1+\frac{\delta ^{2}\left( nN\right) ^{-1}\log ^{2}n}{%
25m_{2}^{2}+5c^{*}\delta _{n}}\right) ,\text{ }m_{2}^{2}\sim N^{-1},
\]
\[
a_{2}\left( 3\right) =11n\left( 1+\frac{5m_{3}^{6/7}}{\delta _{n}}\right) ,%
\text{ }m_{3}=\max_{1\leq i\leq n}\left\| \zeta _{\mathbf{\theta ,}%
j,j^{\prime },i}\right\| _{3}\leq cN^{1/3}.
\]
Observe that $5c\delta _{n}=o(1)$ by Assumption A6, then by taking
$q\ $such that $\left[ \frac{n}{q+1}\right] \geq c_{0}\log n$,
$q\geq c_{1}n/\log n$ for some constants $c_{0},c_{1}$, one has
$a_{1}=O(n/q)=O\left( \log n\right) $, $a_{2}\left( 3\right)
=o\left( n^{2}\right) $ via Assumption A6 again. Assumption A5
yields that
\[
\alpha \left( \left[ \frac{n}{q+1}\right] \right) ^{6/7}\leq \left\{
K_{0}\exp \left( -\lambda _{0}\left[ \frac{n}{q+1}\right] \right) \right\}
^{6/7}\leq Cn^{-6\lambda _{0}c_{0}/7}.
\]
Thus, for fixed $\mathbf{\theta }\in S_{c}^{d-1}$, when $n$ large enough
\begin{equation}
P\left\{ \frac{1}{n}\left| \sum_{i=1}^{n}\zeta _{\mathbf{\theta ,}
j,j^{\prime },i}\right| >\delta _{n}\right\} \leq c\log n\exp \left\{
-c_{2}\delta ^{2}\log n\right\} +Cn^{2-6\lambda _{0}c_{0}/7}.
\label{EQ:unifxi}
\end{equation}
We divide each range of $\theta_{p}$, $p=1,2,...,d-1$, into
$n^{6/(d-1)}\ $ equally spaced intervals with disjoint endpoints
$-1=\theta _{p,0}<\theta _{p,1}<...<\theta
_{p,M_{n}}=1$, for $p=1,...,d-1$. Projecting these small cylinders onto $%
S_{c}^{d-1}$, the radius of each patch $\Lambda _{r}$, $r=1,...,M_{n}$ is
bounded by $cM_{n}^{-1}$. Denote the projection of the $M_{n}$ points as $%
\mathbf{\theta }_{r}=\left( \mathbf{\theta }_{r,-d},\sqrt{1-\left\| \mathbf{%
\theta }_{r,-d}\right\| _{2}^{2}}\right) $, $r=0,1,...,M_{n}$. Employing the
discretization method, $\sup\limits_{\mathbf{\theta }\in
S_{c}^{d-1}}\max\limits_{\substack{ 1\leq j, j^{\prime}\leq N
}}\left| \zeta _{\mathbf{\theta },j,j^{\prime },i}\right| $ is bounded by
\begin{equation}
\sup_{0\leq r\leq M_{n}}\max_{\substack{ 1\leq j, j^{\prime}\leq N }}\left|
\zeta _{\mathbf{\theta }_{r},j,j^{\prime },i}\right| +\sup_{0\leq r\leq
M_{n}}\max_{\substack{ 1\leq j, j^{\prime}\leq N }}\sup\limits_{\mathbf{%
\theta }\in \Lambda _{r}}\left| \zeta _{\mathbf{\theta },j,j^{\prime
},i}-\zeta _{\mathbf{\theta }_{r},j,j^{\prime },i}\right|.
\label{EQ:xidecomp}
\end{equation}
By (\ref{EQ:unifxi}) and Assumption A6, there exists large enough value $%
\delta >0$ such that
\[
P\left\{ \frac{1}{n}\left| \sum_{i=1}^{n}\zeta _{\mathbf{\theta }_{r}\mathbf{%
,}j,j^{\prime },i}\right| >\delta _{n}\right\} \leq n^{-10},
\]
which implies that
\[
\sum_{n=1}^{\infty }P\left\{ \max_{\substack{ 1\leq j, j^{\prime}\leq N }%
}\left| n^{-1}\sum_{l=1}^{n}\zeta _{\mathbf{\theta }_{r},j,j^{\prime
},i}\right| \geq \delta _{n}\right\} \leq 2\sum_{n=1}^{\infty
}N^{2}M_{n}n^{-10}\leq C\sum_{n=1}^{\infty }n^{-3}<\infty .
\]
Thus, Borel-Cantelli Lemma entails that
\begin{equation}
\sup_{0\leq r\leq M_{n}}\max_{\substack{ 1\leq j, j^{\prime}\leq N }}\left|
n^{-1}\sum_{l=1}^{n}\zeta _{\mathbf{\theta }_{r},j,j^{\prime },i}\right|
=O\left( \frac{\log n}{\sqrt{nN}}\right) ,a.s..  \label{EQ:discretebound1}
\end{equation}
Employing Lipschitz continuity of the cubic B-spline, one has with
probability 1
\begin{equation}
\sup_{0\leq r\leq M_{n}}\max_{\substack{ 1\leq j, j^{\prime}\leq N }%
}\sup\limits_{\mathbf{\theta }\in \Lambda _{r}}\left|
n^{-1}\sum_{i=1}^{n}\left\{ \zeta _{\mathbf{\theta ,}j,j^{\prime },i}-\zeta
_{\mathbf{\theta }_{r}\mathbf{,}j,j^{\prime },i}\right\} \right| =O\left(
M_{n}^{-1}h^{-6}\right) .  \label{EQ:discretexidiff}
\end{equation}
Therefore Assumption A2, (\ref{EQ:xidecomp}), (\ref{EQ:discretebound1}) and (%
\ref{EQ:discretexidiff}) lead to the desired result.\hfill

Denote by $\Gamma =\Gamma ^{\left( 0\right) }\cup \Gamma ^{(1)}\cup \Gamma
^{(2)}$ the space of all linear, quadratic and cubic spline functions on $%
\left[ 0,1\right] $. We establish the uniform rate at which the empirical
inner product approximates the theoretical inner product for all B-splines $%
B_{j,k}$ with $k=2,3,4$.

\begin{lemma}
\label{LEM:Anorder} Under Assumptions A2, A5 and A6, one has
\begin{equation}
A_{n}=\sup\limits_{\mathbf{\theta }\in S_{c}^{d-1}}\sup\limits_{\gamma
_{1},\gamma _{2}\in \Gamma }\left| \frac{\left\langle \gamma _{1},\gamma
_{2}\right\rangle _{n,\mathbf{\theta }}-\left\langle \gamma _{1},\gamma
_{2}\right\rangle _{\mathbf{\theta }}}{\left\| \gamma _{1}\right\| _{2,%
\mathbf{\theta }}\left\| \gamma _{2}\right\| _{2,\mathbf{\theta }}}\right|
=O\left\{ \left( nh\right) ^{-1/2}\log n\right\} ,a.s..
\label{EQ:order of An}
\end{equation}
\end{lemma}

\noindent \textbf{Proof.} Denote without loss of generality,
\[
\gamma _{1}=\sum_{k=2}^{4}\sum_{j=-k+1}^{N}\alpha _{jk}B_{j,k}\text{, }%
\gamma _{2}=\sum_{k=2}^{4}\sum_{j=-k+1}^{N}\beta _{jk}B_{j,k}\text{,}
\]
for any two $3\left( N+3\right) $-vectors
\[
\mathbf{\alpha =}\left( \alpha _{-1,2},\alpha _{0,2},...,\alpha
_{N,2},...,\alpha _{N,4}\right) ,\mathbf{\beta =}\left( \beta _{-1,2},\beta
_{0,2},...,\beta _{N,2},...,\beta _{N,4}\right) .
\]
Then for fixed $\mathbf{\theta }$
\begin{eqnarray*}
\left\langle \gamma _{1},\gamma _{2}\right\rangle _{n,\mathbf{\theta }} &=&%
\frac{1}{n}\sum_{i=1}^{n}\left\{ \sum_{k=2}^{4}\sum_{j=-k+1}^{N}\alpha
_{j,k}B_{j,k}\left( U_{\mathbf{\theta },i}\right) \right\} \left\{
\sum_{k=2}^{4}\sum_{j=-k+1}^{N}\beta _{j,k}B_{j,k}\left( U_{\mathbf{\theta }%
,i}\right) \right\}  \\
&=&\sum_{k=2}^{4}\sum_{j=-k+1}^{N}\sum_{k^{\prime }=2}^{4}\sum_{j^{\prime
}=-k+1}^{N}\alpha _{j,k}\beta _{j^{\prime },k^{\prime }}\left\langle
B_{j,k},B_{j^{\prime },k^{\prime }}\right\rangle _{n,\mathbf{\theta }},
\end{eqnarray*}
\begin{eqnarray*}
\left\| \gamma _{1}\right\| _{2,\mathbf{\theta }}^{2}
&=&\sum_{k=2}^{4}\sum_{j=-k+1}^{N}\sum_{k^{\prime }=2}^{4}\sum_{j^{\prime
}=-k+1}^{N}\alpha _{j,k}\alpha _{j^{\prime },k^{\prime }}\left\langle
B_{j,k},B_{j^{^{\prime }},k^{\prime }}\right\rangle _{\mathbf{\theta }}, \\
\left\| \gamma _{2}\right\| _{2,\mathbf{\theta }}^{2}
&=&\sum_{k=2}^{4}\sum_{j=-k+1}^{N}\sum_{k^{\prime }=2}^{4}\sum_{j^{\prime
}=-k+1}^{N}\beta _{j,k}\beta _{j^{\prime },k^{\prime }}\left\langle
B_{j,k},B_{j^{^{\prime }},k^{\prime }}\right\rangle _{\mathbf{\theta }}.
\end{eqnarray*}
According to Lemma \ref{LEM:normequiv}, one has for any $\mathbf{\theta }\in
S_{c}^{d-1}$,
\[
c_{1}h\left\| \mathbf{\alpha }\right\| _{2}^{2}\leq \left\| \gamma
_{1}\right\| _{2,\mathbf{\theta }}^{2}\leq c_{2}h\left\| \mathbf{\alpha }%
\right\| _{2}^{2},c_{1}h\left\| \mathbf{\beta }\right\| _{2}^{2}\leq \left\|
\gamma _{2}\right\| _{2,\mathbf{\theta }}^{2}\leq c_{2}h\left\| \mathbf{%
\beta }\right\| _{2}^{2},
\]
\[
c_{1}h\left\| \mathbf{\alpha }\right\| _{2}\left\| \mathbf{\beta }\right\|
_{2}\leq \left\| \gamma _{1}\right\| _{2,\mathbf{\theta }}\left\| \gamma
_{2}\right\| _{2,\mathbf{\theta }}\leq c_{2}h\left\| \mathbf{\alpha }%
\right\| _{2}\left\| \mathbf{\beta }\right\| _{2}.
\]
Hence
\begin{eqnarray*}
A_{n} &=&\sup\limits_{\mathbf{\theta }\in S_{c}^{d-1}}\sup\limits_{\gamma
_{1}\in \gamma ,\gamma _{2}\in \Gamma }\left| \frac{\left\langle \gamma
_{1},\gamma _{2}\right\rangle _{n,\mathbf{\theta }}-\left\langle \gamma
_{1},\gamma _{2}\right\rangle _{\mathbf{\theta }}}{\left\| \gamma
_{1}\right\| _{2,\mathbf{\theta }}\left\| \gamma _{2}\right\| _{2,\mathbf{%
\theta }}}\right| \leq \frac{\left\| \mathbf{\alpha }\right\| _{\infty
}\left\| \mathbf{\beta }\right\| _{\infty }}{c_{1}h\left\| \mathbf{\alpha }%
\right\| _{2}\left\| \mathbf{\beta }\right\| _{2}} \\
&&\times \sup\limits_{\mathbf{\theta }\in S_{c}^{d-1}}\max_{\substack{
k,k^{\prime }=2,3,4 \\ 1\leq j, j^{\prime }\leq N }}\left| \frac{1}{n}%
\sum_{i=1}^{n}\left\{ \left\langle B_{j,k},B_{j^{^{\prime }},k^{\prime
}}\right\rangle _{n,\mathbf{\theta }}-\left\langle B_{j,k},B_{j^{^{\prime
}},k^{\prime }}\right\rangle _{\mathbf{\theta }}\right\} \right| ,
\end{eqnarray*}
\[
A_{n}\leq c_{0}h^{-1}\sup\limits_{\mathbf{\theta }\in S_{c}^{d-1}}\max_{%
\substack{ k,k^{\prime
}=2,3,4 \\ 1\leq j, j^{\prime }\leq N }}\left| \frac{1}{n}%
\sum_{i=1}^{n}\left\{ \left\langle B_{j,k},B_{j^{^{\prime }},k^{\prime
}}\right\rangle _{n,\mathbf{\theta }}-\left\langle B_{j,k},B_{j^{^{\prime
}},k^{\prime }}\right\rangle _{\mathbf{\theta }}\right\} \right| ,
\]
which, together with Lemma \ref{LEM:th-em product}, imply (\ref{EQ:order of
An}).\hfill

\vskip 0.10in \noindent \textbf{A.2. Proof of Proposition \ref
{PROP:ghattheta-gtheta}}

For any fixed $\mathbf{\theta }$, we write the response $\mathbf{Y}%
^{T}=\left( Y_{1},...,Y_{n}\right) $ as the sum of a signal vector \textbf{$%
\gamma $}$_{\mathbf{\theta }}$, a parametric noise vector $\mathbf{E}_{%
\mathbf{\theta }}$ and a systematic noise vector $\mathbf{E}$,
i.e.,
\[
\mathbf{Y}=\mathbf{\gamma }_{\mathbf{\theta }}+\mathbf{E}_{\mathbf{\theta }}+%
\mathbf{E,}
\]
in which the vectors \textbf{$\gamma $}$_{\mathbf{\theta
}}^{T}=\left\{
\gamma _{\mathbf{\theta }}\left( U_{\mathbf{\theta ,}1}\right) ,...,\gamma _{%
\mathbf{\theta }}\left( U_{\mathbf{\theta ,}n}\right) \right\} $, $\mathbf{E}%
^{T}=\left\{ \sigma \left( \mathbf{X}_{1}\right) \varepsilon
_{1},...,\sigma
\left( \mathbf{X}_{n}\right) \varepsilon _{n}\right\} $ and $\mathbf{E}_{%
\mathbf{\theta }}^{T}=\left\{ m\left( \mathbf{X}_{1}\right) -\gamma _{%
\mathbf{\theta }}\left( U_{\mathbf{\theta },1}\right) ,...,m\left( \mathbf{X}%
_{n}\right) -\gamma _{\mathbf{\theta }}\left( U_{\mathbf{\theta
},n}\right) \right\} $.

\vskip 0.1in \noindent \textbf{Remark A.1.} If $m$ is a genuine
single-index function, then $\mathbf{E}_{\mathbf{\theta}_0}\equiv
0$, thus the proposed SIP model is exactly the single-index model.

Let $\Gamma _{n,\mathbf{\ \theta }}^{\left( 2\right) }$ be the
cubic spline
space spanned by $\left\{ \mathbf{B}_{j,4}\left( U_{\mathbf{\theta }%
,i}\right) \right\} _{i=1}^{n}$, $-3\leq j\leq N$ for fixed
$\mathbf{\theta } $. Projecting $\mathbf{Y}$ onto $\Gamma
_{n,\mathbf{\ \theta }}^{\left( 2\right) }$ yields that
\[
\hat{\mathbf{\gamma }}_{\mathbf{\theta }}=\left\{ \hat{\gamma}_{\mathbf{%
\theta }}\left( U_{\mathbf{\theta ,}1}\right) ,...,\hat{\gamma}_{\mathbf{%
\theta }}\left( U_{\mathbf{\theta },n}\right) \right\} ^{T}=\text{Proj}%
_{\Gamma _{n,\mathbf{\theta }}^{\left( 2\right) }}\mathbf{\gamma }_{\mathbf{%
\theta }}+\text{Proj}_{\Gamma _{n,\mathbf{\theta }}^{\left( 2\right) }}%
\mathbf{E}_{\mathbf{\theta }}+\text{Proj}_{\Gamma _{n,\mathbf{\theta }%
}^{\left( 2\right) }}\mathbf{E},
\]
where $\hat{\gamma}_{\mathbf{\theta }}$ is given in
(\ref{DEF:mthetahat}). We break the cubic spline estimation error $\hat{\gamma}_{\mathbf{\theta }%
}\left( u_{\mathbf{\theta }}\right) -\gamma _{\mathbf{\theta }}\left( u_{%
\mathbf{\theta }}\right) $ into a bias term $\tilde{\gamma}_{\mathbf{\theta }%
}\left( u_{\mathbf{\theta }}\right) -\gamma _{\mathbf{\theta }}\left( u_{%
\mathbf{\theta }}\right) $ and two noise terms $\tilde{\varepsilon}_{\mathbf{%
\theta }}\left( u_{\mathbf{\theta }}\right) $ and $\hat{\varepsilon}_{%
\mathbf{\theta }}\left( u_{\mathbf{\theta }}\right) $
\begin{equation}
\hat{\gamma}_{\mathbf{\theta }}\left( u_{\mathbf{\theta }}\right) -\gamma _{%
\mathbf{\theta }}\left( u_{\mathbf{\theta }}\right) =\left\{ \tilde{\gamma}_{%
\mathbf{\theta }}\left( u_{\mathbf{\theta }}\right) -\gamma
_{\mathbf{\theta
}}\left( u_{\mathbf{\theta }}\right) \right\} +\tilde{\varepsilon}_{\mathbf{%
\theta }}\left( u_{\mathbf{\theta }}\right) +\hat{\varepsilon}_{\mathbf{%
\theta }}\left( u_{\mathbf{\theta }}\right) ,
\label{EQ:decompose1}
\end{equation}
where
\begin{equation}
\tilde{\gamma}_{\mathbf{\theta }}\left( u\right) =\left\{
B_{j,4}\left(
u\right) \right\} _{-3\leq j\leq N}^{T}\mathbf{V}_{n,\mathbf{\theta }%
}^{-1}\left\{ \left\langle \mathbf{\gamma }_{\mathbf{\theta }%
},B_{j,4}\right\rangle _{n,\mathbf{\theta }}\right\} _{j=-3}^{N},
\label{DEF:gtilde}
\end{equation}
\begin{equation}
\tilde{\varepsilon}_{\mathbf{\theta }}\left( u\right) =\left\{
B_{j,4}\left(
u\right) \right\} _{-3\leq j\leq N}^{T}\mathbf{V}_{n,\mathbf{\theta }%
}^{-1}\left\{ \left\langle \mathbf{E}_{\mathbf{\theta }},B_{j,4}\right%
\rangle _{n,\mathbf{\theta }}\right\} _{j=-3}^{N},
\label{DEF:epstilde}
\end{equation}
\begin{equation}
\hat{\varepsilon}_{\mathbf{\theta }}\left( u\right) =\left\{
B_{j,4}\left(
u\right) \right\} _{-3\leq j\leq N}^{T}\mathbf{V}_{n,\mathbf{\theta }%
}^{-1}\left\{ \left\langle \mathbf{E},B_{j,4}\right\rangle _{n,\mathbf{%
\theta }}\right\} _{j=-3}^{N}.  \label{DEF:epshat}
\end{equation}
In the above, we denote by $\mathbf{V}_{n,\mathbf{\theta }}$ the
empirical inner product matrix of the cubic B-spline basis and
similarly, the theoretical inner product matrix as
$\mathbf{V}_{\mathbf{\theta }}$
\begin{equation}
\mathbf{V}_{n,\mathbf{\theta }}=\frac{1}{n}\mathbf{B}_{\mathbf{\theta }}^{T}%
\mathbf{B}_{\mathbf{\theta }}=\left\{ \left\langle B_{j^{\prime
},4},B_{j,4}\right\rangle _{n,\mathbf{\theta }}\right\}
_{j,j^{\prime }=-3}^{N},\mathbf{V}_{\mathbf{\theta }}=\left\{
\left\langle B_{j^{\prime },4},B_{j,4}\right\rangle
_{\mathbf{\theta }}\right\} _{j,j^{\prime }=-3}^{N}.
\label{DEF:Vntheta}
\end{equation}

In Lemma \ref{LEM:ninvBB}, we provide the uniform upper bound of
$\left\| \mathbf{V}_{n,\mathbf{\theta }}^{-1}\right\|_{\infty }$
and $\left\| \mathbf{V}_{\mathbf{\theta }}^{-1}\right\|_{\infty
}$. Before that, we first describe a special case of Theorem
13.4.3 in DeVore and Lorentz (1993).
\begin{lemma}
\label{LEM:DVL13.4.3} If a bi-infinite matrix
with bandwidth $r$ has a bounded inverse $\mathbf{A}^{-1}$ on $l_{2}$ and $%
\kappa =\kappa \left( \mathbf{A}\right) :=\left\| \mathbf{A}\right\|
_{2}\left\| \mathbf{A}^{-1}\right\| _{2}$ is the condition number of $%
\mathbf{A}$, then $\left\| \mathbf{A}^{-1}\right\| _{\infty }\leq
2c_{0}\left( 1-\nu \right) ^{-1}$, with $c_{0}=\nu ^{-2r}\left\| \mathbf{A}%
^{-1}\right\| _{2}$, $\nu =\left( \kappa ^{2}-1\right) ^{1/4r}\left( \kappa
^{2}+1\right) ^{-1/4r}$.
\end{lemma}

\begin{lemma}
\label{LEM:ninvBB} Under Assumptions A2, A5 and A6, there exist
constants $0<c_{V}<C_{V}$ such that $c_{V}N^{-1}\left\|
\mathbf{w}\right\| _{2}^{2}\mathbf{\leq w}^{T}\mathbf{V}_{\mathbf{\theta }}%
\mathbf{w}\leq C_{V}N^{-1}\left\| \mathbf{w}\right\| _{2}^{2}$ and
\begin{equation}
c_{V}N^{-1}\left\| \mathbf{w}\right\| _{2}^{2}\mathbf{\leq w}^{T}\mathbf{V}%
_{n,\mathbf{\theta }}\mathbf{w}\leq C_{V}N^{-1}\left\| \mathbf{w}\right\|
_{2}^{2},a.s.,  \label{EQ:bound}
\end{equation}
with matrices $\mathbf{V}_{\mathbf{\theta }}$ and $\mathbf{V}_{n,\mathbf{%
\theta }}$ defined in (\ref{DEF:Vntheta}). In addition, there exists a
constant $C>0$ such that
\begin{equation}
\sup\limits_{\mathbf{\theta }\in S_{c}^{d-1}}\left\| \mathbf{V}_{n,\mathbf{%
\theta }}^{-1}\right\| _{\infty }\leq CN,a.s.,\sup\limits_{\mathbf{\theta }%
\in S_{c}^{d-1}}\left\| \mathbf{V}_{\mathbf{\theta }}^{-1}\right\| _{\infty
}\leq CN.  \label{EQ:ninvBB}
\end{equation}
\end{lemma}

\noindent \textbf{Proof.} First we compute the lower and upper bounds for
the eigenvalues of $\mathbf{V}_{n,\mathbf{\theta }}$. Let $\mathbf{w}$ be
any $\left( N+4\right) $-vector and denote $\gamma _{\mathbf{w}}\left(
u\right) =\sum_{j=-3}^{N}w_{j}B_{j,4}\left( u\right) $, then $\mathbf{B}_{%
\mathbf{\theta }}\mathbf{w}=\left\{ \gamma _{\mathbf{w}}\left( U_{\mathbf{%
\theta },1}\right) ,...,\gamma _{\mathbf{w}}\left( U_{\mathbf{\theta }%
,n}\right) \right\} ^{T}$ and the definition of $A_{n}$ in (\ref{EQ:order of
An}) from Lemma \ref{LEM:Anorder} entails that
\begin{equation}
\left\| \gamma _{\mathbf{w}}\right\| _{2,\mathbf{\theta }}^{2}\left(
1-A_{n}\right) \leq \mathbf{w}^{T}\mathbf{V}_{n,\mathbf{\theta }}\mathbf{w}%
=\left\| \gamma _{\mathbf{w}}\right\| _{2,n,\mathbf{\theta }}^{2}\leq
\left\| \gamma _{\mathbf{w}}\right\| _{2,\mathbf{\theta }}^{2}\left(
1+A_{n}\right) .  \label{EQ:gw-em}
\end{equation}
Using Theorem 5.4.2 of DeVore and Lorentz (1993) and Assumption A2, one
obtains that
\begin{equation}
c_{f}\frac{C}{N}\left\| \mathbf{w}\right\| _{2}^{2}\leq \Vert \gamma _{%
\mathbf{w}}\Vert _{2,\mathbf{\theta }}^{2}=\mathbf{w}^{T}\mathbf{V}_{\mathbf{%
\theta }}\mathbf{w=}\left\| \sum_{j=-3}^{N}w_{j}B_{j,4}\right\| _{2,\mathbf{%
\theta }}^{2}\leq C_{f}\frac{C}{N}\left\| \mathbf{w}\right\| _{2}^{2},
\label{EQ:gworder}
\end{equation}
which, together with (\ref{EQ:gw-em}), yield
\begin{equation}
c_{f}CN^{-1}\left\| \mathbf{w}\right\| _{2}^{2}\left( 1-A_{n}\right) \leq
\mathbf{w}^{T}\mathbf{V}_{n,\mathbf{\theta }}\mathbf{w}\leq
C_{f}CN^{-1}\left\| \mathbf{w}\right\| _{2}^{2}\left( 1+A_{n}\right) .
\label{EQ:gworder-em}
\end{equation}
Now the order of $A_{n}$ in (\ref{EQ:order of An}), together with (\ref
{EQ:gworder}) and (\ref{EQ:gworder-em}) implies (\ref{EQ:bound}), in which $%
c_{V}=c_{f}C,C_{V}=C_{f}C$. Next, denote by $\lambda _{\max }\left( \mathbf{V%
}_{n,\mathbf{\theta }}\right) $ and $\lambda _{\min }\left( \mathbf{V}_{n,%
\mathbf{\theta }}\right) $ the maximum and minimum eigenvalue of $\mathbf{V}%
_{n,\mathbf{\theta }}$, simple algebra and (\ref{EQ:bound}) entail that
\[
C_{V}N^{-1}\geq \left\| \mathbf{V}_{n,\mathbf{\theta }}\right\| _{2}=\lambda
_{\max }\left( \mathbf{V}_{n,\mathbf{\theta }}\right) ,\left\| \mathbf{V}_{n,%
\mathbf{\theta }}^{-1}\right\| _{2}=\lambda _{\min }^{-1}\left( \mathbf{V}%
_{n,\mathbf{\theta }}\right) \leq c_{V}^{-1}N,a.s.,
\]
thus
\[
\kappa :=\left\| \mathbf{V}_{n,\mathbf{\theta }}\right\| _{2}\left\| \mathbf{%
V}_{n,\mathbf{\theta }}^{-1}\right\| _{2}=\lambda _{\max }\left( \mathbf{V}%
_{n,\mathbf{\theta }}\right) \lambda _{\min }^{-1}\left( \mathbf{V}_{n,%
\mathbf{\theta }}\right) \leq C_{V}c_{V}^{-1}<\infty ,a.s..
\]
Meanwhile, let $\mathbf{w}_{j}=$ the $\left( N+4\right) $-vector with all
zeros except the $j$-th element being $1,j=-3,...,N$. Then clearly
\[
\mathbf{w}_{j}^{T}\mathbf{V}_{n,\mathbf{\theta }}\mathbf{w}_{j}=\frac{1}{n}%
\sum_{i=1}^{n}B_{j,4}^{2}\left( U_{\mathbf{\theta },i}\right) =\left\|
B_{j,4}\right\| _{n,\mathbf{\theta }}^{2},\left\| \mathbf{w}_{j}\right\|
_{2}=1,-3\leq j\leq N
\]
and in particular
\begin{eqnarray*}
\mathbf{w}_{0}^{T}\mathbf{V}_{n,\mathbf{\theta }}\mathbf{w}_{0} &\leq
&\lambda _{\max }\left( \mathbf{V}_{n,\mathbf{\theta }}\right) \left\|
\mathbf{w}_{0}\right\| _{2}=\lambda _{\max }\left( \mathbf{V}_{n,\mathbf{%
\theta }}\right) , \\
\mathbf{w}_{-3}^{T}\mathbf{V}_{n,\mathbf{\theta }}\mathbf{w}_{-3} &\geq
&\lambda _{\min }\left( \mathbf{V}_{n,\mathbf{\theta }}\right) \left\|
\mathbf{w}_{-3}\right\| _{2}=\lambda _{\min }\left( \mathbf{V}_{n,\mathbf{%
\theta }}\right) .
\end{eqnarray*}
This, together with (\ref{EQ:order of An}) yields that
\[
\kappa =\lambda _{\max }\left( \mathbf{V}_{n,\mathbf{\theta }}\right)
\lambda _{\min }^{-1}\left( \mathbf{V}_{n,\mathbf{\theta }}\right) \geq
\frac{\mathbf{w}_{0}^{T}\mathbf{V}_{n,\mathbf{\theta }}\mathbf{w}_{0}}{%
\mathbf{w}_{-3}^{T}\mathbf{V}_{n,\mathbf{\theta }}\mathbf{w}_{-3}}=\frac{%
\left\| B_{0,4}\right\| _{n,\mathbf{\theta }}^{2}}{\left\| B_{-3,4}\right\|
_{n,\mathbf{\theta }}^{2}}\geq \frac{\left\| B_{0,4}\right\| _{\mathbf{%
\theta }}^{2}}{\left\| B_{-3,4}\right\| _{\mathbf{\theta }}^{2}}\frac{1-A_{n}%
}{1+A_{n}},
\]
which leads to $\kappa \geq C>1,a.s.$ because the definition of
B-spline and
Assumption A2 ensure that $\left\| B_{0,4}\right\| _{\mathbf{\theta }%
}^{2}\geq C_{0}\left\| B_{-3,4}\right\| _{\mathbf{\theta }}^{2}$ for some
constant $C_{0}>1$. Next applying Lemma \ref{LEM:DVL13.4.3} with $\nu
=\left( \kappa ^{2}-1\right) ^{1/16}\left( \kappa ^{2}+1\right) ^{-1/16}$
and $c_{0}=\nu ^{-8}\left\| \mathbf{V}_{n,\mathbf{\theta }}^{-1}\right\| _{2}
$, one gets $\left\| \mathbf{V}_{n,\mathbf{\theta }}^{-1}\right\| _{\infty
}\leq 2\nu ^{-8}N\left( 1-\nu \right) ^{-1}=CN,a.s.$. Hence part one of (\ref
{EQ:ninvBB}) follows. Part two of (\ref{EQ:ninvBB}) is proved in the same
fashion.\hfill

In the following, we denote by $Q_{T}\left( m\right) $ the $4$-th
order quasi-interpolant of $m$ corresponding to the knots $T$, see
equation (4.12), page 146 of DeVore and Lorentz (1993). According
to Theorem 7.7.4, DeVore and Lorentz (1993), the following lemma
holds.
\begin{lemma}
\label{LEM:DVL7.7.4} There exists a constant $C>0 $, such that for
$0\leq k\leq 2$ and $\gamma \in C^{\left( 4\right) }\left[
0,1\right] $
\[
\left\| \left( \gamma -Q_{T}\left( \gamma \right) \right) ^{\left(
k\right) }\right\| _{\infty }\leq C\left\| \gamma ^{\left(
4\right) }\right\| _{\infty }h^{4-k},
\]
\end{lemma}

\begin{lemma}
\label{LEM:luniformbiasrate} Under Assumptions A2, A3, A5 and A6,
there exists an absolute constant $C>0$, such that for function
$\tilde{\gamma}_{\mathbf{\theta }}\left( u\right) $ in (\ref
{DEF:gtilde})
\begin{equation}
\sup\limits_{\mathbf{\theta }\in S_{c}^{d-1}}\left\| \frac{d^{k}}{du^{k}}%
\left( \tilde{\gamma}_{\mathbf{\theta }}-\gamma _{\mathbf{\theta
}}\right) \right\| _{\infty }\leq C\left\| m^{\left( 4\right)
}\right\| _{\infty }h^{4-k},a.s.,0\leq k\leq 2,
\label{EQ:biasbound}
\end{equation}
\end{lemma}

\noindent \textbf{Proof.} According to Theorem A.1 of Huang
(2003), there exists an absolute constant $C>0$, such that
\begin{equation}
\sup\limits_{\mathbf{\theta }\in S_{c}^{d-1}}\left\| \tilde{\gamma}_{\mathbf{%
\theta }}-\gamma _{\mathbf{\theta }}\right\| _{\infty }\leq C\sup\limits_{%
\mathbf{\theta }\in S_{c}^{d-1}}\inf\limits_{\gamma \in \Gamma
^{\left( 2\right) }}\left\| \gamma -\gamma _{\mathbf{\theta
}}\right\| _{\infty }\leq C\left\| m^{\left( 4\right) }\right\|
_{\infty }h^{4},a.s., \label{EQ:Huang2003}
\end{equation}
which proves (\ref{EQ:biasbound}) for the case $k=0$. Applying
Lemma \ref
{LEM:DVL7.7.4}, one has for $0\leq k\leq 2$%
\begin{equation}
\sup\limits_{\mathbf{\theta }\in S_{c}^{d-1}}\left\| \frac{d^{k}}{du^{k}}%
\left\{ Q_{T}\left( \gamma _{\mathbf{\theta }}\right) -\gamma _{\mathbf{%
\theta }}\right\} \right\| _{\infty }\leq
C\sup\limits_{\mathbf{\theta }\in S_{c}^{d-1}}\left\| \gamma
_{\mathbf{\theta }}^{\left( 4\right) }\right\| _{\infty
}h^{4-k}\leq C\left\| m^{\left( 4\right) }\right\| _{\infty
}h^{4-k},  \label{EQ:QTapprox}
\end{equation}
As a consequence of (\ref{EQ:Huang2003}) and (\ref{EQ:QTapprox})
for the case $k=0$, one has
\[
\sup\limits_{\mathbf{\theta }\in S_{c}^{d-1}}\left\| Q_{T}\left( \gamma _{%
\mathbf{\theta }}\right) -\tilde{\gamma}_{\mathbf{\theta
}}\right\| _{\infty }\leq C\left\| m^{\left( 4\right) }\right\|
_{\infty }h^{4},a.s.,
\]
which, according to the differentiation of B-spline given in de
Boor (2001), entails that
\begin{equation}
\sup\limits_{\mathbf{\theta }\in S_{c}^{d-1}}\left\| \frac{d^{k}}{du^{k}}%
\left\{ Q_{T}\left( \gamma _{\mathbf{\theta }}\right) -\tilde{\gamma}_{%
\mathbf{\theta }}\right\} \right\| _{\infty }\leq C\left\|
m^{\left( 4\right) }\right\| _{\infty }h^{4-k},a.s.,\mbox{\ }0\leq
k\leq 2. \label{EQ:QTapproxmtilde}
\end{equation}
Combining (\ref{EQ:QTapprox}) and (\ref{EQ:QTapproxmtilde}) proves
(\ref {EQ:biasbound}) for $k=1,2.$ \hfill

\begin{lemma}
\label{LEM:uniformbiasratethetaderiv} Under Assumptions A1, A2, A4
and A5, there exists an absolute constant $C>0 $, such that
\begin{equation}
\sup\limits_{1\leq p\leq d}\sup\limits_{\mathbf{\theta }\in
S_{c}^{d-1}}\left\| \frac{\partial }{\partial \theta _{p}}\left\{ \tilde{%
\gamma}_{\mathbf{\theta }}\left( U_{\mathbf{\theta },i}\right) -\gamma _{%
\mathbf{\theta }}\left( U_{\mathbf{\theta },i}\right) \right\}
_{i=1}^{n}\right\| _{\infty }\leq C\left\| m^{\left( 4\right)
}\right\| _{\infty }h^{3},a.s.,  \label{EQ:biasboundtthetaderiv}
\end{equation}
\begin{equation}
\sup\limits_{1\leq p,q\leq d}\sup\limits_{\mathbf{\theta }\in
S_{c}^{d-1}}\left\| \frac{\partial ^{2}}{\partial \theta
_{p}\partial \theta
_{q}}\left\{ \tilde{\gamma}_{\mathbf{\theta }}\left( U_{\mathbf{\theta }%
,i}\right) -\gamma _{\mathbf{\theta }}\left( U_{\mathbf{\theta
},i}\right) \right\} _{i=1}^{n}\right\| _{\infty }\leq C\left\|
m^{\left( 4\right) }\right\| _{\infty }h^{2},a.s..
\label{EQ:biasboundtthetaderiv2}
\end{equation}
\end{lemma}

\noindent \textbf{Proof}. According to the definition of $\tilde{\gamma}_{%
\mathbf{\theta }}$ in (\ref{DEF:gtilde}), and the fact that
$Q_{T}\left( \gamma _{\mathbf{\theta }}\right) $ is a cubic spline
on the knots $T$
\[
\left\{ \left\{ Q_{T}\left( \gamma _{\mathbf{\theta }}\right) -\tilde{\gamma}%
_{\mathbf{\theta }}\right\} \left( U_{\mathbf{\theta },i}\right)
\right\}
_{i=1}^{n}=\mathbf{P}_{\mathbf{\theta }}\left\{ \left\{ Q_{T}\left( \gamma _{%
\mathbf{\theta }}\right) -\gamma _{\mathbf{\theta }}\right\} \left( U_{%
\mathbf{\theta },i}\right) \right\} _{i=1}^{n},
\]
which entails that
\begin{eqnarray*}
&&\left. \frac{\partial }{\partial \theta _{p}}\left\{ \left\{
Q_{T}\left( \gamma _{\mathbf{\theta }}\right)
-\tilde{\gamma}_{\mathbf{\theta }}\right\}
\left( U_{\mathbf{\theta },i}\right) \right\} _{i=1}^{n}=\frac{\partial }{%
\partial \theta _{p}}\mathbf{P}_{\mathbf{\theta }}\left\{ \left\{
Q_{T}\left( \gamma _{\mathbf{\theta }}\right) -\gamma _{\mathbf{\theta }%
}\right\} \left( U_{\mathbf{\theta },i}\right) \right\} _{i=1}^{n}\right.  \\
&=&\mathbf{\dot{P}}_{p}\left\{ \left\{ Q_{T}\left( \gamma _{\mathbf{\theta }%
}\right) -\gamma _{\mathbf{\theta }}\right\} \left( U_{\mathbf{\theta }%
,i}\right) \right\} _{i=1}^{n}+\mathbf{P}_{\mathbf{\theta }}\frac{\partial }{%
\partial \theta _{p}}\left\{ \left\{ Q_{T}\left( \gamma _{\mathbf{\theta }%
}\right) -\gamma _{\mathbf{\theta }}\right\} \left( U_{\mathbf{\theta }%
,i}\right) \right\} _{i=1}^{n}.
\end{eqnarray*}
Since
\begin{eqnarray*}
&&\left. \frac{\partial }{\partial \theta _{p}}\left\{ \left\{
Q_{T}\left( \gamma _{\mathbf{\theta }}\right) -\gamma
_{\mathbf{\theta }}\right\} \left( U_{\mathbf{\theta },i}\right)
\right\} _{i=1}^{n}=\left\{ \left\{
Q_{T}\left( \frac{\partial }{\partial \theta _{p}}\gamma _{\mathbf{\theta }%
}\right) -\frac{\partial }{\partial \theta _{p}}\gamma _{\mathbf{\theta }%
}\right\} \left( U_{\mathbf{\theta },i}\right) \right\} _{i=1}^{n}\right.  \\
&&+\left\{ \frac{d}{du}\left\{ Q_{T}\left( \gamma _{\mathbf{\theta
}}\right) -\gamma _{\mathbf{\theta }}\right\} \left(
U_{\mathbf{\theta },i}\right) X_{ip}\right\} _{i=1}^{n},
\end{eqnarray*}
applying (\ref{EQ:QTapproxmtilde}) to the decomposition above
produces (\ref {EQ:biasboundtthetaderiv}). The proof of
(\ref{EQ:biasboundtthetaderiv2}) is similar. \hfill

\begin{lemma}
\label{LEM:BthetadotBpnorm} Under Assumptions A2, A5 and A6, there
exists a constant $C>0$ such that
\begin{equation}
\sup\limits_{\mathbf{\theta }\in S_{c}^{d-1}}\left\| n^{-1}\mathbf{B}_{%
\mathbf{\theta }}^{T}\right\| _{\infty }\leq Ch,a.s.,\sup\limits_{1\leq
p\leq d}\sup\limits_{\mathbf{\theta }\in S_{c}^{d-1}}\left\| n^{-1}\mathbf{%
\dot{B}}_{p}^{T}\right\| _{\infty }\leq C,a.s.,  \label{EQ:BthetadotBpnorm}
\end{equation}
\begin{equation}
\sup\limits_{\mathbf{\theta }\in S_{c}^{d-1}}\left\| \mathbf{P}_{\mathbf{%
\theta }}\right\| _{\infty }\leq C,a.s.,\sup\limits_{1\leq p\leq
d}\sup\limits_{\mathbf{\theta }\in S_{c}^{d-1}}\left\| \mathbf{\dot{P}}%
_{p}\right\| _{\infty }\leq Ch^{-1},a.s..  \label{EQ:PthetadotPpnorm}
\end{equation}
\end{lemma}

\noindent \textbf{Proof.} To prove (\ref{EQ:BthetadotBpnorm}), observe that
for any vector $\mathbf{a}\in R^{n}$, with probability $1$%
\[
\left\| n^{-1}\mathbf{B}_{\mathbf{\theta }}^{T}\mathbf{a}\right\|
_{\infty }\leq \left\| \mathbf{a}\right\| _{\infty
}\max\limits_{-3\leq j\leq N}\left|
n^{-1}\sum_{i=1}^{n}B_{j,4}\left( U_{\mathbf{\theta },i}\right)
\right| \leq Ch\left\| \mathbf{a}\right\| _{\infty },
\]
\[
\left\| n^{-1}%
\mathbf{\dot{B}}_{p}^{T}\mathbf{a}\right\| _{\infty }\leq \left\|
\mathbf{a}\right\| _{\infty }\max\limits_{-3\leq j\leq N}\left|
\frac{1}{nh}\sum_{i=1}^{n}\left\{ \left( B_{j,3}-B_{j+1,3}\right) \left( U_{%
\mathbf{\theta },i}\right) \right\} \dot{F}_{d}\left( \mathbf{X}_{\mathbf{%
\theta },i}\right) X_{i,p}\right| \leq C\left\| \mathbf{a}\right\|
_{\infty }.
\]
To prove (\ref{EQ:PthetadotPpnorm}), one only needs to use (\ref{EQ:ninvBB}%
), (\ref{EQ:BthetadotBpnorm}) and (\ref{DEF:Ptheta}). \hfill

\begin{lemma}
\label{LEM:BEsupnorm} Under Assumptions A2 and A4-A6, one has with
probability $1$
\begin{equation}
\sup\limits_{\mathbf{\theta }\in S_{c}^{d-1}}\left\| \frac{\mathbf{B}_{%
\mathbf{\theta }}^{T}\mathbf{E}}{n}\right\| _{\infty }=\max_{-3\leq j\leq
N}\left| n^{-1}\sum_{i=1}^{n}B_{j,4}\left( U_{\mathbf{\theta },i}\right)
\sigma \left( \mathbf{X}_{i}\right) \varepsilon _{i}\right| =O\left( \frac{%
\log n}{\sqrt{nN}}\right) ,  \label{EQ:BEsupnorm}
\end{equation}
\begin{equation}
\sup\limits_{1\leq p\leq d}\sup\limits_{\mathbf{\theta }\in
S_{c}^{d-1}}\left\| \frac{\partial }{\partial \theta _{p}}\left( \frac{%
\mathbf{B}_{\mathbf{\theta }}^{T}\mathbf{E}}{n}\right) \right\|
_{\infty }=\sup\limits_{1\leq p\leq d}\sup\limits_{\mathbf{\theta
}\in S_{c}^{d-1}}\left\|
\frac{\mathbf{\dot{B}}_{p}^{T}\mathbf{E}}{n}\right\| _{\infty
}=O\left( \frac{\log n}{\sqrt{nh}}\right).
\label{EQ:BEsupnormderiv}
\end{equation}
Similarly, under Assumptions A2, A4-A6, with probability 1
\begin{equation}
\sup\limits_{\mathbf{\theta }\in S_{c}^{d-1}}\left\| \frac{\mathbf{B}_{%
\mathbf{\theta }}^{T}\mathbf{E}_{\mathbf{\theta }}}{n}\right\| _{\infty
}=\sup\limits_{\mathbf{\theta }\in S_{c}^{d-1}}\max_{-3\leq j\leq N}\left|
\frac{1}{n}\sum_{i=1}^{n}B_{j,4}\left( U_{\mathbf{\theta },i}\right) \left\{
m\left( \mathbf{X}_{i}\right) -\gamma _{\mathbf{\theta }}\left( U_{\mathbf{%
\theta },i}\right) \right\} \right|=O\left( \frac{\log
n}{\sqrt{nN}}\right), \label{EQ:BEthetasupnorm}
\end{equation}
\begin{equation}
\sup\limits_{1\leq p\leq d}\sup\limits_{\mathbf{\theta }\in
S_{c}^{d-1}}\left\| \frac{\partial }{\partial \theta _{p}}\left( \frac{%
\mathbf{B}_{\mathbf{\theta }}^{T}\mathbf{E}_{\mathbf{\theta }}}{n}\right)
\right\| _{\infty }=O\left( \frac{\log n}{\sqrt{nh}}\right) ,a.s..
\label{EQ:BEthetasupnormderiv}
\end{equation}
\end{lemma}

\noindent \textbf{Proof.} We decompose the noise variable $\varepsilon _{i}$
into a truncated part and a tail part $\varepsilon _{i}=\varepsilon
_{i,1}^{D_{n}}+\varepsilon _{i,2}^{D_{n}}+m_{i}^{D_{n}}$, where $%
D_{n}=n^{\eta }\left( 1/3<\eta <2/5\right) $, $\varepsilon
_{i,1}^{D_{n}}=\varepsilon _{i}I\left\{ \left| \varepsilon _{i}\right|
>D_{n}\right\} $,
\[
\varepsilon _{i,2}^{D_{n}}=\varepsilon _{i}I\left\{ \left| \varepsilon
_{i}\right| \leq D_{n}\right\} -m_{i}^{D_{n}},m_{i}^{D_{n}}=E\left[
\varepsilon _{i}I\left\{ \left| \varepsilon _{i}\right| \leq D_{n}\right\} |%
\mathbf{X}_{i}\right] .
\]
It is straightforward to verify that the mean of the truncated part is
uniformly bounded by $D_{n}^{-2}$, so the boundedness of B spline basis and
of the function $\sigma ^{2}$ entail that
\[
\sup\limits_{\mathbf{\theta }\in S_{c}^{d-1}}\left| \frac{1}{n}%
\sum_{i=1}^{n}B_{j,4}\left( U_{\mathbf{\theta },i}\right) \sigma \left(
\mathbf{X}_{i}\right) m_{i}^{D_{n}}\right| =O\left( D_{n}^{-2}\right)
=o\left( n^{-2/3}\right) .
\]
The tail part vanishes almost surely
\[
\sum_{n=1}^{\infty }P\left\{ \left| \varepsilon _{n}\right| >D_{n}\right\}
\leq \sum_{n=1}^{\infty }D_{n}^{-3}<\infty .
\]
Borel-Cantelli Lemma implies that
\[
\left| \frac{1}{n}\sum_{i=1}^{n}B_{j,4}\left( U_{\mathbf{\theta },i}\right)
\sigma \left( \mathbf{X}_{i}\right) \varepsilon _{i,1}^{D_{n}}\right|
=O\left( n^{-k}\right) ,\text{ for any }k>0.
\]
For the truncated part, using Bernstein's inequality and
discretization as in Lemma \ref{LEM:th-em product}
\[
\sup\limits_{\mathbf{\theta }\in S_{c}^{d-1}}\sup_{1\leq j\leq N}\left|
n^{-1}\sum_{i=1}^{n}B_{j,4}\left( U_{\mathbf{\theta },i}\right) \sigma
\left( \mathbf{X}_{i}\right) \varepsilon _{i,2}^{D_{n}}\right| =O\left( \log
n/\sqrt{nN}\right) ,a.s..
\]
Therefore (\ref{EQ:BEsupnorm}) is established as with probability $1$
\[
\sup\limits_{\mathbf{\theta }\in S_{c}^{d-1}}\left\| \frac{1}{n}\mathbf{%
\mathbf{B}_{\mathbf{\theta }}^{T}\mathbf{E}}\right\| _{\infty }=o\left(
n^{-2/3}\right) +O\left( n^{-k}\right) +O\left( \log n/\sqrt{nN}\right)
=O\left( \log n/\sqrt{nN}\right) .
\]
The proofs of (\ref{EQ:BEsupnormderiv}), (\ref{EQ:BEthetasupnorm}) are
similar as $E\left\{ m\left( \mathbf{X}_{i}\right) -\gamma _{\mathbf{\theta }%
}\left( U_{\mathbf{\theta },i}\right) \left| U_{\mathbf{\theta },i}\right.
\right\} \equiv 0$, but no truncation is needed for (\ref{EQ:BEthetasupnorm}%
) as $\sup\limits_{\mathbf{\theta }\in S_{c}^{d-1}}\max\limits_{1\leq i\leq
n}\left| m\left( \mathbf{X}_{i}\right) -\gamma _{\mathbf{\theta }}\left( U_{%
\mathbf{\theta },i}\right) \right| \leq C<\infty $. Meanwhile, to prove (\ref
{EQ:BEthetasupnormderiv}), we note that for any $p=1,...,d$
\[
\frac{\partial }{\partial \theta _{p}}\left( \mathbf{B}_{\mathbf{\theta }%
}^{T}\mathbf{E}_{\mathbf{\theta }}\right) =\left\{ \sum_{i=1}^{n}\frac{%
\partial }{\partial \theta _{p}}\left[ B_{j,4}\left( U_{\mathbf{\theta }%
,i}\right) \left\{ m\left( \mathbf{X}_{i}\right) -\gamma _{\mathbf{\theta }%
}\left( U_{\mathbf{\theta },i}\right) \right\} \right] \right\} _{j=-3}^{N}.
\]
According to (\ref{EQ:mthetagtheta}), one has $\gamma _{\mathbf{\theta }%
}\left( U_{\mathbf{\theta }}\right) \equiv E\left\{ m\left( \mathbf{X}%
\right) |U_{\mathbf{\theta }}\right\} $, hence
\[
E\left[ B_{j,4}\left( U_{\mathbf{\theta }}\right) \left\{ m\left( \mathbf{X}%
\right) -\gamma _{\mathbf{\theta }}\left( U_{\mathbf{\theta }}\right)
\right\} \right] \equiv 0,-3\leq j\leq N,\mathbf{\theta }\in S_{c}^{d-1}.
\]
Applying Assumptions A2 and A3, one can differentiate through the
expectation, thus
\[
E\left\{ \frac{\partial }{\partial \theta _{p}}\left[ B_{j,4}\left( U_{%
\mathbf{\theta }}\right) \left\{ m\left( \mathbf{X}\right) -\gamma _{\mathbf{%
\theta }}\left( U_{\mathbf{\theta }}\right) \right\} \right] \right\} \equiv
0,1\leq p\leq d,-3\leq j\leq N,\mathbf{\theta }\in S_{c}^{d-1},
\]
which allows one to apply the Bernstein's inequality to obtain that with
probability $1$
\[
\left\| \left\{ \frac{1}{n}\sum_{i=1}^{n}\frac{\partial }{\partial \theta
_{p}}\left[ B_{j,4}\left( U_{\mathbf{\theta },i}\right) \left\{ m\left(
\mathbf{X}_{i}\right) -\gamma _{\mathbf{\theta }}\left( U_{\mathbf{\theta }%
,i}\right) \right\} \right] \right\} _{j=-3}^{N}\right\| _{\infty }=O\left\{
(nh)^{-1/2}\log n\right\} ,
\]
which is (\ref{EQ:BEthetasupnormderiv}). \hfill

\begin{lemma}
\label{LEM:epshatorder} Under Assumptions A2 and A4-A6, for
$\hat{\varepsilon}_{\mathbf{\theta }}\left( u\right) $ in (\ref
{DEF:epshat}), one has
\begin{equation}
\sup\limits_{\mathbf{\theta }\in S_{c}^{d-1}}\sup_{u\in \left[ 0,1\right]
}\left| \hat{\varepsilon}_{\mathbf{\theta }}\left( u\right) \right|
=O\left\{ \left( nh\right) ^{-1/2}\log n\right\} ,a.s..
\label{EQ:epshatunifsize}
\end{equation}
\end{lemma}

\noindent \textbf{Proof.} Denote $\mathbf{\hat{a}}\equiv \left( \hat{a}%
_{-3},\cdots ,\hat{a}_{N}\right) ^{T}\mathbf{=}\left( \mathbf{B}_{\mathbf{%
\theta }}^{T}\mathbf{B}_{\mathbf{\theta }}\right) ^{-1}\mathbf{B}_{\mathbf{%
\theta }}^{T}\mathbf{E}=\mathbf{V}_{n,\mathbf{\theta }}^{-1}\left( n^{-1}%
\mathbf{\mathbf{B}_{\mathbf{\theta }}^{T}\mathbf{E}}\right) $, then $\hat{%
\varepsilon}_{\mathbf{\theta }}\left( u\right) =\sum_{j=-3}^{N}\hat{a}%
_{j}B_{j,4}\left( u\right) $, so the order of $\hat{\varepsilon}_{\mathbf{%
\theta }}\left( u\right) $ is related to that of $\mathbf{\hat{a}}$. In
fact, by Theorem 5.4.2 in DeVore and Lorentz (1993)
\begin{eqnarray*}
\sup\limits_{\mathbf{\theta }\in S_{c}^{d-1}}\sup_{u\in \left[ 0,1\right]
}\left| \hat{\varepsilon}_{\mathbf{\theta }}\left( u\right) \right| &\leq
&\sup\limits_{\mathbf{\theta }\in S_{c}^{d-1}}\left\| \mathbf{\hat{a}}%
\right\| _{\infty }= \\
\sup\limits_{\mathbf{\theta }\in S_{c}^{d-1}}\left\| \mathbf{V}_{n,\mathbf{%
\theta }}^{-1}\left( n^{-1}\mathbf{\mathbf{B}_{\mathbf{\theta }}^{T}\mathbf{E%
}}\right) \right\| _{\infty } &\leq &CN\sup\limits_{\mathbf{\theta }\in
S_{c}^{d-1}}\left\| n^{-1}\mathbf{\mathbf{B}_{\mathbf{\theta }}^{T}\mathbf{E}%
}\right\| _{\infty },a.s.,
\end{eqnarray*}
where the last inequality follows from (\ref{EQ:ninvBB}) of Lemma \ref
{LEM:ninvBB}. Applying (\ref{EQ:BEsupnorm}) of Lemma \ref{LEM:BEsupnorm}, we
have established (\ref{EQ:epshatunifsize}).\hfill

\begin{lemma}
\label{LEM:epstildeorder} Under Assumptions A2
and A4-A6, for $\tilde{\varepsilon}_{\mathbf{\theta }}\left( u\right) $ in (%
\ref{DEF:epstilde}), one has
\begin{equation}
\sup\limits_{\mathbf{\theta }\in S_{c}^{d-1}}\sup_{u\in \left[ 0,1\right]
}\left| \tilde{\varepsilon}_{\mathbf{\theta }}\left( u\right) \right|
=O\left\{ \left( nh\right) ^{-1/2}\log n\right\} ,a.s..
\label{EQ:epstildeunifsize}
\end{equation}
\end{lemma}

The proof is similar to Lemma \ref{LEM:epshatorder}, thus omitted. \hfill

The next result evaluates the uniform size of the noise derivatives.

\begin{lemma}
\label{LEM:epsderivative} Under Assumptions A2-A6, one has with
probability $1$
\begin{eqnarray}
\sup\limits_{1\leq p\leq d}\sup\limits_{\mathbf{\theta }\in
S_{c}^{d-1}}\max_{1\leq i\leq n}\left| \frac{\partial }{\partial \theta _{p}}%
\hat{\varepsilon}_{\mathbf{\theta }}\left( U_{\mathbf{\theta },i}\right)
\right| =O\left\{ (nh^{3})^{-1/2}\log n\right\} ,
\label{EQ:epshatderivative} \\
\sup\limits_{1\leq p\leq d}\sup\limits_{\mathbf{\theta }\in
S_{c}^{d-1}}\max_{1\leq i\leq n}\left| \frac{\partial }{\partial \theta _{p}}%
\tilde{\varepsilon}_{\mathbf{\theta }}\left( U_{\mathbf{\theta },i}\right)
\right| =O\left\{ (nh^{3})^{-1/2}\log n\right\} ,
\label{EQ:epstildederivative} \\
\sup\limits_{1\leq p,q\leq d}\sup\limits_{\mathbf{\theta }\in
S_{c}^{d-1}}\max_{1\leq i\leq n}\left| \frac{\partial ^{2}}{\partial \theta
_{p}\partial \theta _{q}}\hat{\varepsilon}_{\mathbf{\theta }}\left( U_{%
\mathbf{\theta },i}\right) \right| =O\left\{ (nh^{5})^{-1/2}\log n\right\} ,
\label{EQ:epshatderivative2} \\
\sup\limits_{1\leq p,q\leq d}\sup\limits_{\mathbf{\theta }\in
S_{c}^{d-1}}\max_{1\leq i\leq n}\left| \frac{\partial ^{2}}{\partial \theta
_{p}\partial \theta _{q}}\tilde{\varepsilon}_{\mathbf{\theta }}\left( U_{%
\mathbf{\theta },i}\right) \right| =O\left\{ (nh^{5})^{-1/2}\log n\right\} .
\label{EQ:epstildederivative2}
\end{eqnarray}
\end{lemma}

\noindent \textbf{Proof.} Note that
\[
\left\{ \frac{\partial }{\partial \theta _{p}}\hat{\varepsilon}_{\mathbf{%
\theta }}\left( U_{\mathbf{\theta },i}\right) \right\} _{i=1}^{n}=\left(
\mathbf{I}-\mathbf{P}_{\mathbf{\theta }}\right) \mathbf{\dot{B}}_{p}\left(
\mathbf{B}_{\mathbf{\theta }}^{T}\mathbf{B}_{\mathbf{\theta }}\right) ^{-1}%
\mathbf{B}_{\mathbf{\theta }}^{T}\mathbf{E}+\mathbf{B}_{\mathbf{\theta }%
}\left( \mathbf{B}_{\mathbf{\theta }}^{T}\mathbf{B}_{\mathbf{\theta }%
}\right) ^{-1}\mathbf{\dot{B}}_{p}^{T}\left( \mathbf{I}-\mathbf{P}_{\mathbf{%
\theta }}\right) \mathbf{E}.
\]
Applying (\ref{EQ:BEsupnorm}) and (\ref{EQ:BEsupnormderiv}) of Lemma \ref
{LEM:BEsupnorm}, (\ref{EQ:ninvBB}) of Lemma \ref{LEM:ninvBB}, (\ref
{EQ:BthetadotBpnorm}) and (\ref{EQ:PthetadotPpnorm}) of Lemma \ref
{LEM:BthetadotBpnorm}, one derives (\ref{EQ:epshatderivative}). To prove (%
\ref{EQ:epstildederivative}), note that
\begin{equation}
\left\{ \frac{\partial }{\partial \theta _{p}}\tilde{\varepsilon }_{\mathbf{%
\theta }}\left( U_{\mathbf{\theta },i}\right) \right\} _{i=1}^{n}=\frac{%
\partial }{\partial \theta _{p}}\left\{ \mathbf{P}_{\mathbf{\theta }}\mathbf{%
E}_{\mathbf{\theta }}\right\} =\mathbf{\dot{P}}_{p}\mathbf{E}_{\mathbf{%
\theta }}+\mathbf{P}_{\mathbf{\theta }}\frac{\partial }{\partial \theta _{p}}%
\mathbf{E}_{\mathbf{\theta }}=T_{1}+T_{2},  \label{DEF:T1T2}
\end{equation}
in which
\begin{eqnarray*}
T_{1} &=&\left\{ \left( \mathbf{I}-\mathbf{P}_{\mathbf{\theta }}\right)
\mathbf{\dot{B}}_{p}-\mathbf{B}_{\mathbf{\theta }}\left( \mathbf{B}_{\mathbf{%
\theta }}^{T}\mathbf{B}_{\mathbf{\theta }}\right) ^{-1}\mathbf{\dot{B}}%
_{p}^{T}\mathbf{B}_{\mathbf{\theta }}\right\} \left( \mathbf{B}_{\mathbf{%
\theta }}^{T}\mathbf{B}_{\mathbf{\theta }}\right) ^{-1}\mathbf{B}_{\mathbf{%
\theta }}^{T}\mathbf{E}_{\mathbf{\theta }} \\
&=&\left\{ \left( \mathbf{I}-\mathbf{P}_{\mathbf{\theta }}\right) \mathbf{%
\dot{B}}_{p}-\mathbf{B}_{\mathbf{\theta }}\left( \frac{\mathbf{B}_{\mathbf{%
\theta }}^{T}\mathbf{B}_{\mathbf{\theta }}}{n}\right) ^{-1}\frac{\mathbf{%
\dot{B}}_{p}^{T}\mathbf{B}_{\mathbf{\theta }}}{n}\right\} \left( \frac{%
\mathbf{B}_{\mathbf{\theta }}^{T}\mathbf{B}_{\mathbf{\theta }}}{n}\right)
^{-1}\frac{\mathbf{B}_{\mathbf{\theta }}^{T}\mathbf{E}_{\mathbf{\theta }}}{n}%
,
\end{eqnarray*}
\[
T_{2}=\mathbf{B}_{\mathbf{\theta }}\left( \frac{\mathbf{B}_{\mathbf{\theta }%
}^{T}\mathbf{B}_{\mathbf{\theta }}}{n}\right) ^{-1}\frac{\partial }{\partial
\theta _{p}}\left( \frac{\mathbf{B}_{\mathbf{\theta }}^{T}\mathbf{E}_{%
\mathbf{\theta }}}{n}\right) .
\]
By (\ref{EQ:BEsupnorm}), (\ref{EQ:ninvBB}), (\ref{EQ:BthetadotBpnorm}) and (%
\ref{EQ:PthetadotPpnorm}), one derives
\begin{equation}
\sup\limits_{\mathbf{\theta }\in S_{c}^{d-1}}\left\| T_{1}\right\| _{\infty
}=O\left( n^{-1/2}N^{3/2}\log n\right) ,a.s.,  \label{EQ:T1bound}
\end{equation}
while (\ref{EQ:BEthetasupnormderiv}) of Lemma \ref{LEM:BEsupnorm}, (\ref
{EQ:ninvBB}) of Lemma \ref{LEM:ninvBB}
\begin{equation}
\sup\limits_{\mathbf{\theta }\in S_{c}^{d-1}}\left\| T_{2}\right\| _{\infty
}=N\times O\left( n^{-1/2}h^{-1/2}\log n\right) =O\left(
n^{-1/2}h^{-3/2}\log n\right) ,a.s..  \label{EQ:T2bound}
\end{equation}
Now, putting together (\ref{DEF:T1T2}), (\ref{EQ:T1bound}) and (\ref
{EQ:T2bound}), we have established (\ref{EQ:epstildederivative}). The proof
for (\ref{EQ:epshatderivative2}) and (\ref{EQ:epstildederivative2}) are
similar. \hfill

\noindent \textbf{Proof of Proposition \ref{PROP:ghattheta-gtheta}.}
According to the decomposition (\ref{EQ:decompose1})
\[
\left| \hat{\gamma}_{\mathbf{\theta }}\left( u\right) -\gamma _{\mathbf{%
\theta }}\left( u\right) \right| =\left| \left\{ \tilde{\gamma}_{\mathbf{%
\theta }}\left( u\right) -\gamma _{\mathbf{\theta }}\left( u\right) \right\}
+\tilde{\varepsilon}_{\mathbf{\theta }}\left( u\right) +\hat{\varepsilon}_{%
\mathbf{\theta }}\left( u\right) \right| .
\]
Then (\ref{EQ:ghattheta-gtheta}) follows directly from (\ref{EQ:biasbound})
of Lemma \ref{LEM:luniformbiasrate}, (\ref{EQ:epshatunifsize}) of Lemma \ref
{LEM:epshatorder} and (\ref{EQ:epstildeunifsize}) of Lemma \ref
{LEM:epstildeorder}. Again by definitions (\ref{DEF:epstilde}) and (\ref
{DEF:epshat}), we write
\[
\frac{\partial }{\partial \theta _{p}}\left\{ \left( \hat{\gamma}_{\mathbf{%
\theta }}-\gamma _{\mathbf{\theta }}\right) \left( U_{\mathbf{\theta }%
,i}\right) \right\} =\frac{\partial }{\partial \theta _{p}}\left( \tilde{%
\gamma}_{\mathbf{\theta }}-\gamma _{\mathbf{\theta }}\right) \left( U_{%
\mathbf{\theta },i}\right) +\frac{\partial }{\partial \theta _{p}}\tilde{%
\gamma}_{\mathbf{\theta }}\left( U_{\mathbf{\theta },i}\right) +\frac{%
\partial }{\partial \theta _{p}}\hat{\varepsilon}_{\mathbf{\theta }}\left(
U_{\mathbf{\theta },i}\right) .
\]
It is clear from (\ref{EQ:biasboundtthetaderiv}), (\ref{EQ:epshatderivative}%
) and (\ref{EQ:epstildederivative}) that with probability $1$
\[
\sup\limits_{1\leq p\leq d}\sup\limits_{\mathbf{\theta }\in
S_{c}^{d-1}}\max\limits_{1\leq i\leq n}\left| \frac{\partial }{\partial
\theta _{p}}\left( \tilde{\gamma}_{\mathbf{\theta }}-\gamma _{\mathbf{\theta
}}\right) \left( U_{\mathbf{\theta },i}\right) \right| =O\left( h^{3}\right)
,
\]
\[
\sup\limits_{1\leq p\leq d}\sup\limits_{\mathbf{\theta }\in
S_{c}^{d-1}}\max\limits_{1\leq i\leq n}\left\{ \left| \frac{\partial }{%
\partial \theta _{p}}\tilde{\varepsilon}_{\mathbf{\theta }}\left( U_{\mathbf{%
\theta },i}\right) \right| +\left| \frac{\partial }{\partial \theta _{p}}%
\hat{\varepsilon}_{\mathbf{\theta }}\left( U_{\mathbf{\theta },i}\right)
\right| \right\} =O\left\{ \left( nh^{3}\right) ^{-1/2}\log n\right\} .
\]
Putting together all the above yields (\ref{EQ:ghatthetaderiv-gthetaderv}).
The proof of (\ref{EQ:ghatthetaderiv2-gthetaderv2}) is similar. \hfill

\vskip 0.10in \noindent \textbf{A.3. Proof of Proposition \ref{PROP:unif}}

\begin{lemma}
\label{LEM:k0} Under Assumptions A2-A6, one has
\[
\sup\limits_{\mathbf{\theta }\in S_{c}^{d-1}}\left| \hat{R}\left( \mathbf{%
\theta }\right) -R\left( \mathbf{\theta }\right) \right| =o(1),a.s..
\]
\end{lemma}

\noindent \textbf{Proof.} For the empirical risk function $\hat{R}\left(
\mathbf{\theta }\right) $ in (\ref{DEF:Rhat}), one has
\[
\hat{R}\left( \mathbf{\theta }\right) =n^{-1}\sum_{i=1}^{n}\left\{ \hat{%
\gamma}_{\mathbf{\theta }}\left( U_{\mathbf{\theta },i}\right) -m\left(
\mathbf{X}_{i}\right) -\sigma \left( \mathbf{X}_{i}\right) \varepsilon
_{i}\right\} ^{2}
\]
\[
=n^{-1}\sum_{i=1}^{n}\left\{ \hat{\gamma}_{\mathbf{\theta }}\left( U_{%
\mathbf{\theta },i}\right) -\gamma _{\mathbf{\theta }}\left( U_{\mathbf{%
\theta },i}\right) +\gamma _{\mathbf{\theta }}\left( U_{\mathbf{\theta }%
,i}\right) -m\left( \mathbf{X}_{i}\right) -\sigma \left( \mathbf{X}%
_{i}\right) \varepsilon _{i}\right\} ^{2},
\]
hence
\[
\hat{R}\left( \mathbf{\theta }\right) =n^{-1}\sum_{i=1}^{n}\left\{ \hat{%
\gamma}_{\mathbf{\theta }}\left( U_{\mathbf{\theta },i}\right) -\gamma _{%
\mathbf{\theta }}\left( U_{\mathbf{\theta },i}\right) \right\}
^{2}+n^{-1}\sum_{i=1}^{n}\sigma ^{2}\left( \mathbf{X}_{i}\right) \varepsilon
_{i}^{2}
\]
\[
+2n^{-1}\sum_{i=1}^{n}\left\{ \hat{\gamma}_{\mathbf{\theta }}\left( U_{%
\mathbf{\theta },i}\right) -\gamma _{\mathbf{\theta }}\left( U_{\mathbf{%
\theta },i}\right) \right\} \left\{ \gamma _{\mathbf{\theta }}\left( U_{%
\mathbf{\theta },i}\right) -m\left( \mathbf{X}_{i}\right) -\sigma \left(
\mathbf{X}_{i}\right) \varepsilon _{i}\right\}
\]
\[
+n^{-1}\sum_{i=1}^{n}\left\{ \gamma _{\mathbf{\theta }}\left( U_{\mathbf{%
\theta },i}\right) -m\left( \mathbf{X}_{i}\right) \right\}
^{2}+2n^{-1}\sum_{i=1}^{n}\left\{ \gamma _{\mathbf{\theta }}\left( U_{%
\mathbf{\theta },i}\right) -m\left( \mathbf{X}_{i}\right) \right\} \sigma
\left( \mathbf{X}_{i}\right) \varepsilon _{i},
\]
where $\hat{\gamma}_{\mathbf{\theta }}\left( x\right) $ is defined in (\ref
{DEF:mthetahat}). Using the expression of $R\left( \mathbf{\theta }\right) $
in (\ref{EQ:Rthetagtheta}), one has
\[
\sup\limits_{\mathbf{\theta }\in S_{c}^{d-1}}\left| \hat{R}\left( \mathbf{%
\theta }\right) -R\left( \mathbf{\theta }\right) \right| \leq
I_{1}+I_{2}+I_{3}+I_{4},
\]
with
\[
I_{1}=\sup\limits_{\mathbf{\theta }\in S_{c}^{d-1}}\left|
n^{-1}\sum_{i=1}^{n}\left\{ \hat{\gamma}_{\mathbf{\theta }}\left( U_{\mathbf{%
\theta },i}\right) -\gamma _{\mathbf{\theta }}\left( U_{\mathbf{\theta }%
,i}\right) \right\} ^{2}\right| ,
\]
\[
I_{2}=\sup\limits_{\mathbf{\theta }\in S_{c}^{d-1}}\left|
2n^{-1}\sum_{i=1}^{n}\left\{ \hat{\gamma}_{\mathbf{\theta }}\left( U_{%
\mathbf{\theta },i}\right) -\gamma _{\mathbf{\theta }}\left( U_{\mathbf{%
\theta },i}\right) \right\} \left\{ \gamma _{\mathbf{\theta }}\left( U_{%
\mathbf{\theta },i}\right) -m\left( \mathbf{X}_{i}\right) -\sigma \left(
\mathbf{X}_{i}\right) \varepsilon _{i}\right\} \right| ,
\]
\[
I_{3}=\sup\limits_{\mathbf{\theta }\in S_{c}^{d-1}}\left|
n^{-1}\sum_{i=1}^{n}\left\{ \gamma _{\mathbf{\theta }}\left( U_{\mathbf{%
\theta },i}\right) -m\left( \mathbf{X}_{i}\right) \right\} ^{2}-E\left\{
\gamma _{\mathbf{\theta }}\left( U_{\mathbf{\theta }}\right) -m\left(
\mathbf{X}\right) \right\} ^{2}\right| ,
\]
\[
I_{4}=\sup\limits_{\mathbf{\theta }\in S_{c}^{d-1}}\left\{ \left| \frac{1}{n}%
\sum_{i=1}^{n}\sigma ^{2}\left( \mathbf{X}_{i}\right) \varepsilon
_{i}^{2}-E\sigma ^{2}\left( \mathbf{X}\right) \right| +\left| \frac{2}{n}%
\sum_{i=1}^{n}\left\{ \gamma _{\mathbf{\theta }}\left( U_{\mathbf{\theta }%
,i}\right) -m\left( \mathbf{X}_{i}\right) \right\} \sigma \left( \mathbf{X}%
_{i}\right) \varepsilon _{i}\right| \right\} .
\]
Bernstein inequality and strong law of large number for $\alpha $ mixing
sequence imply that
\begin{equation}
I_{3}+I_{4}=o(1),a.s..  \label{EQ:I3I4order}
\end{equation}
Now (\ref{EQ:ghattheta-gtheta}) of Proposition \ref{PROP:ghattheta-gtheta}
provides that
\[
\sup\limits_{\mathbf{\theta }\in S_{c}^{d-1}}\sup\limits_{u\in \left[
0,1\right] }\left| \hat{\gamma}_{\mathbf{\theta }}\left( u\right) -\gamma _{%
\mathbf{\theta }}\left( u\right) \right| =O\left( n^{-1/2}h^{-1/2}\log
n+h^{4}\right) ,a.s.,
\]
which entail that
\begin{equation}
I_{1}=O\left\{ \left( n^{-1/2}h^{-1/2}\log n\right) ^{2}+\left( h^{4}\right)
^{2}\right\} ,a.s.,  \label{EQ:I1order}
\end{equation}
\[
I_{2}\leq O\left\{ (nh)^{-1/2}\log n+h^{4}\right\} \times \sup\limits_{%
\mathbf{\theta }\in S_{c}^{d-1}}2n^{-1}\sum_{i=1}^{n}\left| \gamma _{\mathbf{%
\theta }}\left( U_{\mathbf{\theta },i}\right) -m\left( \mathbf{X}_{i}\right)
-\sigma \left( \mathbf{X}_{i}\right) \varepsilon _{i}\right| .
\]
Hence
\begin{equation}
I_{2}\leq O\left( n^{-1/2}h^{-1/2}\log n+h^{4}\right) ,a.s..
\label{EQ:I2order}
\end{equation}
The lemma now follows from (\ref{EQ:I3I4order}), (\ref{EQ:I1order}) and (\ref
{EQ:I2order}) and Assumption A6.\hfill

\begin{lemma}
\label{LEM:Rhat-Lderiv} Under Assumptions A2 - A6, one has
\begin{equation}
\sup\limits_{\mathbf{\theta }\in S_{c}^{d-1}}\sup\limits_{1\leq p\leq
d}\left| \frac{\partial }{\partial \theta _{p}}\left\{ \hat{R}\left( \mathbf{%
\theta }\right) -R\left( \mathbf{\theta }\right) \right\}
-n^{-1}\sum_{i=1}^{n}\xi _{\mathbf{\theta },i,p}\right| =o\left(
n^{-1/2}\right) ,\text{ }a.s.,  \label{EQ:Rhat-Lderivlinearization}
\end{equation}
in which
\begin{equation}
\xi _{\mathbf{\theta },i,p}=2\left\{ \gamma _{\mathbf{\theta }}\left( U_{%
\mathbf{\theta },i}\right) -Y_{i}\right\} \frac{\partial }{\partial \theta
_{p}}\gamma _{\mathbf{\theta }}\left( U_{\mathbf{\theta },i}\right) -\frac{%
\partial }{\partial \theta _{p}}R\left( \mathbf{\theta }\right) ,\ E\left(
\xi _{\mathbf{\theta },i,p}\right) =0.  \label{DEF:xithetaip}
\end{equation}
Furthermore for $k=1,2$
\begin{equation}
\sup\limits_{\mathbf{\theta }\in S_{c}^{d-1}}\left| \frac{\partial ^{k}}{%
\partial \mathbf{\theta }^{k}}\left\{ \hat{R}\left( \mathbf{\theta }\right)
-R\left( \mathbf{\theta }\right) \right\} \right| =O\left(
n^{-1/2}h^{-1/2-k}\log n+h^{4-k}\right) ,a.s..
\label{EQ:Rhat-Lderivunifrate}
\end{equation}
\end{lemma}

\noindent \textbf{Proof.} Note that for any $p=1,2,...,d$%
\[
\frac{1}{2}\frac{\partial }{\partial \theta _{p}}\hat{R}\left( \mathbf{%
\theta }\right) =n^{-1}\sum_{i=1}^{n}\left\{ \hat{\gamma}_{\mathbf{\theta }%
}\left( U_{\mathbf{\theta },i}\right) -Y_{i}\right\} \frac{\partial }{%
\partial \theta _{p}}\hat{\gamma}_{\mathbf{\theta }}\left( U_{\mathbf{\theta
},i}\right),
\]
\begin{eqnarray*}
\frac{1}{2}\frac{\partial }{\partial \theta _{p}}R\left( \mathbf{\theta }%
\right)  &=&E\left[ \left\{ \gamma _{\mathbf{\theta }}\left( U_{\mathbf{%
\theta }}\right) -m\left( \mathbf{X}\right) \right\} \frac{\partial }{%
\partial \theta _{p}}\gamma _{\mathbf{\theta }}\left( U_{\mathbf{\theta }%
}\right) \right]  \\
&=&E\left[ \left\{ \gamma _{\mathbf{\theta }}\left( U_{\mathbf{\theta }%
}\right) -m\left( \mathbf{X}\right) -\sigma \left( \mathbf{X}\right)
\varepsilon \right\} \frac{\partial }{\partial \theta _{p}}\gamma _{\mathbf{%
\theta }}\left( U_{\mathbf{\theta }}\right) \right] .
\end{eqnarray*}
Thus $E\left( \xi _{\mathbf{\theta },i,p}\right) =2E\left[ \left\{ \gamma _{%
\mathbf{\theta }}\left( U_{\mathbf{\theta },i}\right) -Y_{i}\right\} \frac{%
\partial }{\partial \theta _{p}}\gamma _{\mathbf{\theta }}\left( U_{\mathbf{%
\theta },i}\right) \right] -\frac{\partial }{\partial \theta _{p}}R\left(
\mathbf{\theta }\right) =0$ and
\begin{equation}
\frac{1}{2}\frac{\partial }{\partial \theta _{p}}\left\{ \hat{R}\left(
\mathbf{\theta }\right) -R\left( \mathbf{\theta }\right) \right\} =\left(
2n\right) ^{-1}\sum_{i=1}^{n}\xi _{\mathbf{\theta },i,p}+J_{1,\mathbf{\theta
},p}+J_{2,\mathbf{\theta },p}+J_{3,\mathbf{\theta },p},
\label{EQ:J1J2J3decomp}
\end{equation}
with
\begin{eqnarray*}
J_{1,\mathbf{\theta },p} &=&n^{-1}\sum_{i=1}^{n}\left\{ \hat{\gamma}_{%
\mathbf{\theta }}\left( U_{\mathbf{\theta },i}\right) -\gamma _{\mathbf{%
\theta }}\left( U_{\mathbf{\theta },i}\right) \right\} \frac{\partial }{%
\partial \theta _{p}}\left( \hat{\gamma}_{\mathbf{\theta }}-\gamma _{\mathbf{%
\theta }}\right) \left( U_{\mathbf{\theta },i}\right) , \\
J_{2,\mathbf{\theta },p} &=&n^{-1}\sum_{i=1}^{n}\left\{ \gamma _{\mathbf{%
\theta }}\left( U_{\mathbf{\theta },i}\right) -m\left( \mathbf{X}_{i}\right)
-\sigma \left( \mathbf{X}_{i}\right) \varepsilon _{i}\right\} \frac{\partial
}{\partial \theta _{p}}\left( \hat{\gamma}_{\mathbf{\theta }}-\gamma _{%
\mathbf{\theta }}\right) \left( U_{\mathbf{\theta },i}\right) , \\
J_{3,\mathbf{\theta },p} &=&n^{-1}\sum_{i=1}^{n}\left\{ \hat{\gamma}_{%
\mathbf{\theta }}\left( U_{\mathbf{\theta },i}\right) -\gamma _{\mathbf{%
\theta }}\left( U_{\mathbf{\theta },i}\right) \right\} \frac{\partial }{%
\partial \theta _{p}}\gamma _{\mathbf{\theta }}\left( U_{\mathbf{\theta }%
,i}\right) .
\end{eqnarray*}
Bernstein inequality implies that
\begin{equation}
\sup\limits_{\mathbf{\theta }\in S_{c}^{d-1}}\sup\limits_{1\leq p\leq
d}\left| n^{-1}\sum_{i=1}^{n}\xi _{\mathbf{\theta },i,p}\right| =O\left(
n^{-1/2}\log n\right) ,a.s..  \label{EQ:xiavgorder}
\end{equation}
Meanwhile, applying (\ref{EQ:ghattheta-gtheta}) and (\ref
{EQ:ghatthetaderiv-gthetaderv}) of Proposition \ref{PROP:ghattheta-gtheta},
one obtains that
\begin{eqnarray}
\sup\limits_{\mathbf{\theta }\in S_{c}^{d-1}}\sup\limits_{1\leq p\leq
d}\left| J_{1,\mathbf{\theta },p}\right|  &=&O\left\{ \left( nh\right)
^{-1/2}\log n+h^{4}\right\} \times O\left\{ \left( nh^{3}\right) ^{-1/2}\log
n+h^{3}\right\}   \nonumber \\
&=&O\left( n^{-1}h^{-2}\log ^{2}n+h^{7}\right) ,a.s..  \label{EQ:J1order}
\end{eqnarray}
Note that
\begin{eqnarray*}
J_{2,\mathbf{\theta },p} &=&n^{-1}\sum_{i=1}^{n}\left\{ \gamma _{\mathbf{%
\theta }}\left( U_{\mathbf{\theta },i}\right) -m\left( \mathbf{X}_{i}\right)
-\sigma \left( \mathbf{X}_{i}\right) \varepsilon _{i}\right\} \frac{\partial
}{\partial \theta _{p}}\left( \tilde{\gamma}_{\mathbf{\theta }}-\gamma _{%
\mathbf{\theta }}\right) \left( U_{\mathbf{\theta },i}\right)  \\
&&-n^{-1}\left( \mathbf{E+E}_{\mathbf{\theta }}\right) ^{T}\frac{\partial }{%
\partial \theta _{p}}\left\{ \mathbf{P}_{\mathbf{\theta }}\left( \mathbf{E+E}%
_{\mathbf{\theta }}\right) \right\} .
\end{eqnarray*}
Applying (\ref{EQ:ghattheta-gtheta}), one gets
\[
\sup\limits_{\mathbf{\theta }\in S_{c}^{d-1}}\sup\limits_{1\leq p\leq
d}\left| J_{2,\mathbf{\theta },p}+n^{-1}\left( \mathbf{E+E}_{\mathbf{\theta }%
}\right) ^{T}\frac{\partial }{\partial \theta _{p}}\left\{ \mathbf{P}_{%
\mathbf{\theta }}\left( \mathbf{E+E}_{\mathbf{\theta }}\right) \right\}
\right| =O\left( h^{3}\right) ,a.s.,
\]
while (\ref{EQ:BEsupnorm}), (\ref{EQ:BEthetasupnorm}) and (\ref{EQ:ninvBB})
entail that with probability $1$%
\[
\sup\limits_{\mathbf{\theta }\in S_{c}^{d-1}}\sup\limits_{1\leq p\leq
d}\left| n^{-1}\left( \mathbf{E+E}_{\mathbf{\theta }}\right) ^{T}\frac{%
\partial }{\partial \theta _{p}}\left\{ \mathbf{P}_{\mathbf{\theta }}\left(
\mathbf{E+E}_{\mathbf{\theta }}\right) \right\} \right|
\]
\[
=O\left\{ \left( nN\right) ^{-1/2}\log n\right\} \times N\times N\times
O\left\{ \left( nN\right) ^{-1/2}\log n\right\} =O\left\{ n^{-1}N\log
^{2}n\right\} ,
\]
thus
\begin{equation}
\sup\limits_{\mathbf{\theta }\in S_{c}^{d-1}}\sup\limits_{1\leq p\leq
d}\left| J_{2,\mathbf{\theta },p}\right| =O\left( h^{3}+n^{-1}N\log
^{2}n\right) ,a.s..  \label{EQ:J2order}
\end{equation}
Lastly
\[
J_{3,\mathbf{\theta },p}-n^{-1}\sum_{i=1}^{n}\left( \tilde{\gamma}_{\mathbf{%
\theta }}-\gamma _{\mathbf{\theta }}\right) \frac{\partial }{\partial \theta
_{p}}\gamma _{\mathbf{\theta }}\left( U_{\mathbf{\theta },i}\right)
=n^{-1}\left( \mathbf{E+E}_{\mathbf{\theta }}\right) ^{T}\mathbf{B}_{\mathbf{%
\theta }}\left( \frac{\mathbf{B}_{\mathbf{\theta }}^{T}\mathbf{B}_{\mathbf{%
\theta }}}{n}\right) ^{-1}\frac{\mathbf{B}_{\mathbf{\theta }}^{T}}{n}\frac{%
\partial }{\partial \theta _{p}}\mathbf{\gamma }_{\mathbf{\theta }}.
\]
By applying (\ref{EQ:BEsupnorm}), (\ref{EQ:BEthetasupnorm}), and (\ref
{EQ:ninvBB}), it is clear that with probability $1$%
\begin{eqnarray*}
&&\sup\limits_{\mathbf{\theta }\in S_{c}^{d-1}}\sup\limits_{1\leq p\leq
d}\left| \left( n^{-1}\mathbf{B}_{\mathbf{\theta }}^{T}\mathbf{E+}n^{-1}%
\mathbf{\mathbf{B}_{\mathbf{\theta }}^{T}E}_{\mathbf{\theta }}\right)
^{T}\left( \frac{\mathbf{B}_{\mathbf{\theta }}^{T}\mathbf{B}_{\mathbf{\theta
}}}{n}\right) ^{-1}\frac{\mathbf{B}_{\mathbf{\theta }}^{T}}{n}\frac{\partial
}{\partial \theta _{p}}\mathbf{\gamma }_{\mathbf{\theta }}\right|  \\
&=&O\left\{ \left( nN\right) ^{-1/2}\log n\right\} \times N\times O\left\{
h+\left( nN\right) ^{-1/2}\log n\right\}  \\
&=&O\left\{ n^{-1}\log ^{2}n+\left( nN\right) ^{-1/2}\log n\right\} ,
\end{eqnarray*}
while by applying (\ref{EQ:biasbound}) of Lemma \ref{LEM:luniformbiasrate},
one has
\[
\sup\limits_{\mathbf{\theta }\in S_{c}^{d-1}}\sup\limits_{1\leq p\leq
d}\left| n^{-1}\sum_{i=1}^{n}\left( \tilde{\gamma}_{\mathbf{\theta }}-\gamma
_{\mathbf{\theta }}\right) \frac{\partial }{\partial \theta _{p}}\gamma _{%
\mathbf{\theta }}\left( U_{\mathbf{\theta },i}\right) \right| =O\left(
h^{4}\right) ,a.s.,
\]
together, the above entail that
\begin{equation}
\sup\limits_{\mathbf{\theta }\in S_{c}^{d-1}}\sup\limits_{1\leq p\leq
d}\left| J_{3,\mathbf{\theta },p}\right| =O\left\{ h^{4}+n^{-1}\log
^{2}n+\left( nN\right) ^{-1/2}\log n\right\} ,a.s..  \label{EQ:J3order}
\end{equation}
Therefore, (\ref{EQ:J1J2J3decomp}), (\ref{EQ:J1order}), (\ref{EQ:J2order}), (%
\ref{EQ:J3order}) and Assumption A6 lead to (\ref
{EQ:Rhat-Lderivlinearization}), which, together with (\ref{EQ:xiavgorder}),
establish (\ref{EQ:Rhat-Lderivunifrate}) for $k=1$.

Note that the second order derivative of $\hat{R}\left( \mathbf{\theta }%
\right) $ and $R\left( \mathbf{\theta }\right) $ with respect to $\theta _{p}
$, $\theta _{q}$ are
\[
2n^{-1}\left[ \sum_{i=1}^{n}\left\{ \hat{\gamma}_{\mathbf{\theta }}\left( U_{%
\mathbf{\theta },i}\right) -Y_{i}\right\} \frac{\partial ^{2}}{\partial
\theta _{p}\partial \theta _{q}}\hat{\gamma}_{\mathbf{\theta }}\left( U_{%
\mathbf{\theta },i}\right) +\sum_{i=1}^{n}\frac{\partial }{\partial \theta
_{q}}\hat{\gamma}_{\mathbf{\theta }}\left( U_{\mathbf{\theta },i}\right)
\frac{\partial }{\partial \theta _{p}}\hat{\gamma}_{\mathbf{\theta }}\left(
U_{\mathbf{\theta },i}\right) \right] ,
\]
\[
2\left[ E\left\{ \gamma _{\mathbf{\theta }}\left( U_{\mathbf{\theta }%
}\right) -m\left( \mathbf{X}\right) \right\} \frac{\partial ^{2}}{\partial
\theta _{p}\partial \theta _{q}}\gamma _{\mathbf{\theta }}\left( U_{\mathbf{%
\theta }}\right) +E\left\{ \frac{\partial }{\partial \theta _{q}}\gamma _{%
\mathbf{\theta }}\left( U_{\mathbf{\theta }}\right) \frac{\partial }{%
\partial \theta _{p}}\gamma _{\mathbf{\theta }}\left( U_{\mathbf{\theta }%
}\right) \right\} \right] .
\]
The proof of (\ref{EQ:Rhat-Lderivunifrate}) for $k=2$ follows from (\ref
{EQ:ghattheta-gtheta}), (\ref{EQ:ghatthetaderiv-gthetaderv}) and (\ref
{EQ:ghatthetaderiv2-gthetaderv2}). \hfill

\vskip .05in \noindent \textbf{Proof of Proposition \ref{PROP:unif}.} The
result follows from Lemma \ref{LEM:k0}, Lemma \ref{LEM:Rhat-Lderiv},
equations (\ref{EQ:Ldot}) and (\ref{EQ:Lddot}). \hfill

\vskip .10in \noindent \textbf{A.4. Proof of the Theorem 2}

Let $\hat{S}_{p}^{*}\left( \mathbf{\theta }_{-d}\right) $ be the $p$-th
element of $\hat{S}^{*}\left( \mathbf{\theta }_{-d}\right) $ and for $\gamma
_{\mathbf{\theta }}$ in (\ref{EQ:mthetagtheta}), denote
\begin{equation}
\eta _{i,p}:=2\left\{ \dot{\gamma}_{p}-\theta _{0,p}\theta _{0,d}^{-1}\dot{%
\gamma}_{d}\right\} \left( U_{\mathbf{\theta }_{0},i}\right) \left\{ \gamma
_{\mathbf{\theta }_{0}}\left( U_{\mathbf{\theta }_{0},i}\right)
-Y_{i}\right\} ,  \label{DEF:etaip}
\end{equation}
where $\dot{\gamma}_{p}$ is value of $\frac{\partial }{\partial \theta _{p}}%
\gamma _{\mathbf{\theta }}$ taking at $\mathbf{\theta =\theta }_{0}$, for
any $p,q=1,2,...,d-1$.

\begin{lemma}
\label{LEM:sptheta0} Under Assumptions A2-A6, one has
\begin{equation}
\sup\limits_{1\leq p\leq d-1}\left| \hat{S}_{p}^{*}\left( \mathbf{\theta }%
_{0,-d}\right) -n^{-1}\sum_{i=1}^{n}\eta _{i,p}\right| =o\left(
n^{-1/2}\right) ,a.s..  \label{EQ:sptheta0}
\end{equation}
\end{lemma}

\noindent \textbf{Proof.} For any $p=1,...,d-1$%
\[
\hat{S}_{p}^{*}\left( \mathbf{\theta }_{-d}\right) -S_{p}^{*}\left( \mathbf{%
\theta }_{-d}\right) =\left( \frac{\partial }{\partial \theta
_{p}}-\theta _{p}\theta _{d}^{-1}\frac{\partial }{\partial \theta
_{d}}\right) \left\{ \hat{R}\left( \mathbf{\theta }\right)
-R\left( \mathbf{\theta }\right) \right\}.
\]
Therefore, according to (\ref{EQ:Rhat-Lderivlinearization}), (\ref
{DEF:xithetaip}) and (\ref{DEF:etaip})
\[
\eta _{i,p}=n^{-1}\sum_{i=1}^{n}\xi _{\mathbf{\theta }_{0},i,p}-\theta
_{0,p}\theta _{0,d}^{-1}n^{-1}\sum_{i=1}^{n}\xi _{\mathbf{\theta }_{0},i,d},%
\text{ }E\left( \eta _{i,p}\right) =0,
\]
\[
\sup\limits_{1\leq p\leq d-1}\left| \hat{S}_{p}^{*}\left( \mathbf{\theta }%
_{0,-d}\right) -S_{p}^{*}\left( \mathbf{\theta }_{0,-d}\right)
-n^{-1}\sum_{i=1}^{n}\eta _{i,p}\right| =o\left( n^{-1/2}\right) ,a.s..
\]
Since $S^{*}\left( \mathbf{\theta }_{-d}\right) $ attains its minimum at $%
\mathbf{\theta }_{0,-d}$, for $p=1,...,d-1$
\[
S_{p}^{*}\left( \mathbf{\theta }_{0,-d}\right) \equiv \left. \left( \frac{%
\partial }{\partial \theta _{p}}-\theta _{p}\theta _{d}^{-1}\frac{\partial }{%
\partial \theta _{d}}\right) R\left( \mathbf{\theta }\right) \right| _{%
\mathbf{\theta =\theta }_{0}}\equiv 0,
\]
which yields (\ref{EQ:sptheta0}).\hfill

\begin{lemma}
\label{LEM:lpq} The $\left( p,q\right) $-th entry of the Hessian
matrix $H^{*}\left( \mathbf{\theta }_{0,-d}\right) $ equals
$l_{p,q}$ given in Theorem \ref{THM:normality}.
\end{lemma}

\noindent \textbf{Proof.} It is easy to show that for any $p,q=1,2,...,d$,
\[
\frac{\partial }{\partial \theta _{p}}R\left( \mathbf{\theta }\right) =\frac{%
\partial }{\partial \theta _{p}}E\left\{ m\left( \mathbf{X}\right) -\gamma _{%
\mathbf{\theta }}\left( U_{\mathbf{\theta }}\right) \right\} ^{2}=-2E\left[
\gamma _{\mathbf{\theta }}\left( U_{\mathbf{\theta }}\right) \frac{\partial
}{\partial \theta _{p}}\gamma _{\mathbf{\theta }}\left( U_{\mathbf{\theta }%
}\right) \right] ,
\]
\[
\frac{\partial ^{2}}{\partial \theta _{p}\partial \theta _{q}}R\left(
\mathbf{\theta }\right) =-2E\left[ \frac{\partial }{\partial \theta _{p}}%
\gamma _{\mathbf{\theta }}\left( U_{\mathbf{\theta }}\right) \frac{\partial
}{\partial \theta _{q}}\gamma _{\mathbf{\theta }}\left( U_{\mathbf{\theta }%
}\right) +\gamma _{\mathbf{\theta }}\left( U_{\mathbf{\theta }}\right) \frac{%
\partial ^{2}}{\partial \theta _{p}\partial \theta _{q}}\gamma _{\mathbf{%
\theta }}\left( U_{\mathbf{\theta }}\right) \right] .
\]
Note that
\begin{equation}
\frac{\partial }{\partial \theta _{p}}R^{*}\left( \mathbf{\theta }%
_{-d}\right) =\frac{\partial }{\partial \theta _{p}}R\left( \mathbf{\theta }%
\right) -\frac{\theta _{p}}{\theta _{d}}\frac{\partial }{\partial \theta _{d}%
}R\left( \mathbf{\theta }\right) ,  \label{EQ:Ldot}
\end{equation}
\begin{eqnarray}
&&\frac{\partial ^{2}}{\partial \theta _{p}\partial \theta _{q}}R^{*}\left(
\mathbf{\theta }_{-d}\right) =\frac{\partial ^{2}}{\partial \theta
_{p}\partial \theta _{q}}R\left( \mathbf{\theta }\right) -\frac{\theta _{q}}{%
\theta _{d}}\frac{\partial ^{2}}{\partial \theta _{p}\partial \theta _{d}}%
R\left( \mathbf{\theta }\right) -\frac{\theta _{p}}{\theta _{d}}\frac{%
\partial ^{2}}{\partial \theta _{d}\partial \theta _{q}}R\left( \mathbf{%
\theta }\right)   \nonumber \\
&&-\frac{\partial }{\partial \theta _{q}}\left( \frac{\theta _{p}}{\sqrt{%
1-\left\| \mathbf{\theta }_{-d}\right\| _{2}^{2}}}\right) \frac{\partial }{%
\partial \theta _{d}}R\left( \mathbf{\theta }\right) +\frac{\theta
_{p}\theta _{q}}{\theta _{d}^{2}}\frac{\partial ^{2}}{\partial \theta
_{d}\partial \theta _{d}}R\left( \mathbf{\theta }\right) .  \label{EQ:Lddot}
\end{eqnarray}
Thus
\[
\frac{\partial }{\partial \theta _{p}}R^{*}\left( \mathbf{\theta }%
_{-d}\right) =-2E\left[ \gamma _{\mathbf{\theta }}\left( U_{\mathbf{\theta }%
}\right) \frac{\partial }{\partial \theta _{p}}\gamma _{\mathbf{\theta }%
}\left( U_{\mathbf{\theta }}\right) \right] +2\theta _{d}^{-1}\theta
_{p}E\left[ \gamma _{\mathbf{\theta }}\left( U_{\mathbf{\theta }}\right)
\frac{\partial }{\partial \theta _{d}}\gamma _{\mathbf{\theta }}\left( U_{%
\mathbf{\theta }}\right) \right] ,
\]
\[
\frac{\partial ^{2}}{\partial \theta _{p}\partial \theta _{q}}R^{*}\left(
\mathbf{\theta }_{-d}\right) =-2E\left\{ \frac{\partial }{\partial \theta
_{p}}\gamma _{\mathbf{\theta }}\left( U_{\mathbf{\theta }}\right) \frac{%
\partial }{\partial \theta _{q}}\gamma _{\mathbf{\theta }}\left( U_{\mathbf{%
\theta }}\right) +\gamma _{\mathbf{\theta }}\left( U_{\mathbf{\theta }%
}\right) \frac{\partial ^{2}}{\partial \theta _{p}\partial \theta _{q}}%
\gamma _{\mathbf{\theta }}\left( U_{\mathbf{\theta }}\right) \right\}
\]
\begin{eqnarray*}
&&+2\theta _{q}\theta _{d}^{-1}E\left\{ \frac{\partial }{\partial \theta _{d}%
}\gamma _{\mathbf{\theta }}\left( U_{\mathbf{\theta }}\right) \frac{\partial
}{\partial \theta _{p}}\gamma _{\mathbf{\theta }}\left( U_{\mathbf{\theta }%
}\right) +\gamma _{\mathbf{\theta }}\left( U_{\mathbf{\theta }}\right) \frac{%
\partial ^{2}}{\partial \theta _{p}\partial \theta _{d}}\gamma _{\mathbf{%
\theta }}\left( U_{\mathbf{\theta }}\right) \right\}  \\
&&+2\frac{\partial }{\partial \theta _{q}}\left( \frac{\theta _{p}}{\sqrt{%
1-\left\| \mathbf{\theta }_{-d}\right\| _{2}^{2}}}\right) E\left\{ \gamma _{%
\mathbf{\theta }}\left( U_{\mathbf{\theta }}\right) \frac{\partial }{%
\partial \theta _{d}}\gamma _{\mathbf{\theta }}\left( U_{\mathbf{\theta }%
}\right) \right\}  \\
&&+2\theta _{p}\theta _{d}^{-1}E\left\{ \frac{\partial }{\partial \theta _{p}%
}\gamma _{\mathbf{\theta }}\left( U_{\mathbf{\theta }}\right) \frac{\partial
}{\partial \theta _{q}}\gamma _{\mathbf{\theta }}\left( U_{\mathbf{\theta }%
}\right) +\gamma _{\mathbf{\theta }}\left( U_{\mathbf{\theta }}\right) \frac{%
\partial ^{2}}{\partial \theta _{p}\partial \theta _{q}}\gamma _{\mathbf{%
\theta }}\left( U_{\mathbf{\theta }}\right) \right\}  \\
&&-2\theta _{p}\theta _{q}\theta _{d}^{-2}E\left[ \left\{ \frac{\partial }{%
\partial \theta _{d}}\gamma _{\mathbf{\theta }}\left( U_{\mathbf{\theta }%
}\right) \right\} ^{2}+\gamma _{\mathbf{\theta }}\left( U_{\mathbf{\theta }%
}\right) \frac{\partial ^{2}}{\partial \theta _{d}\partial \theta _{d}}%
\gamma _{\mathbf{\theta }}\left( U_{\mathbf{\theta }}\right) \right] .
\end{eqnarray*}
Therefore we obtained the desired result.\hfill

\vskip .05in \noindent \textbf{Proof of Theorem \ref{THM:normality}.} For
any $p=1,2,...,d-1$, let
\[
f_{p}\left( t\right) =\hat{S}_{p}^{*}\left( t\hat{\mathbf{\theta}}%
_{-d}+\left( 1-t\right) \mathbf{\theta }_{0,-d}\right) ,t\in [0,1],
\]
then
\[
\frac{d}{dt}f_{p}\left( t\right) =\sum_{q=1}^{d-1}\frac{\partial }{\partial
\theta _{q}}\hat{S}_{p}^{*}\left( t\hat{\mathbf{\theta}}_{-d}\mathbf{+}%
\left( 1-t\right) \mathbf{\theta}_{0,-d}\right) \left( \hat{\theta}%
_{q}-\theta _{0,q}\right) .
\]
Note that $\hat{S}^{*}\left( \mathbf{\theta }_{-d}\right) $ attains its
minimum at $\hat{\mathbf{\theta}}_{-d}$, i.e., $\hat{S}_{p}^{*}\left( \hat{%
\mathbf{\theta}}_{-d}\right) \equiv 0$. Thus, for any $p=1,2,...,d-1$, $%
t_{p}\in \left[ 0,1\right] $, one has
\begin{eqnarray*}
&&\left. -\hat{S}_{p}^{*}\left( \mathbf{\theta }_{0,-d}\right) =\hat{S}%
_{p}^{*}\left( \hat{\mathbf{\theta}}_{-d}\right) -\hat{S}_{p}^{*}\left(
\mathbf{\theta }_{0,-d}\right) =f_{p}\left( 1\right) -f_{p}\left( 0\right)
\right. \\
&=&\left\{ \frac{\partial ^{2}}{\partial \theta _{q}\theta _{p}}\hat{R}%
^{*}\left( t_{p}\hat{\mathbf{\theta}}_{-d}+\left( 1-t_{p}\right) \mathbf{%
\theta }_{0,-d}\right) \right\} _{q=1,...,d-1}^{T}\left( \hat{\mathbf{\theta}%
}_{-d}\mathbf{-\theta }_{0,-d}\right) ,
\end{eqnarray*}
then
\[
-\hat{S}^{*}\left( \mathbf{\theta }_{0,-d}\right) =\left\{ \frac{\partial
^{2}}{\partial \theta _{q}\partial \theta _{p}}\hat{R}^{*}\left( t_{p}\hat{%
\mathbf{\theta}}_{-d}+\left( 1-t_{p}\right) \mathbf{\theta }_{0,-d}\right)
\right\} _{\substack{p,q=1,...,d-1}}\left( \hat{\mathbf{\theta}}_{-d}-%
\mathbf{\theta }_{0,-d}\right) .
\]
Now (\ref{EQ:strongconsist}) of Theorem \ref{THM:strconsistent} and
Proposition \ref{PROP:unif} with $k=2$ imply that uniformly in $%
p,q=1,2,...,d-1$
\begin{equation}
\frac{\partial ^{2}}{\partial \theta _{q}\partial \theta _{p}}\hat{R}%
^{*}\left( t_{p}\hat{\mathbf{\theta}}_{-d}\mathbf{+}\left( 1-t_{p}\right)
\mathbf{\theta }_{0,-d}\right) \longrightarrow l_{q,p},a.s.,
\label{EQ:Hessianconverge}
\end{equation}
where $l_{p,q}$ is given in Theorem \ref{THM:normality}. Noting that $\sqrt{n%
}\left( \hat{\mathbf{\theta}}_{-d}\mathbf{-\theta }_{0,-d}\right) $ is
represented as
\[
-\left[ \left\{ \frac{\partial ^{2}}{\partial \theta _{q}\partial \theta _{p}%
}\hat{R}^{*}\left( t_{p}\hat{\mathbf{\theta}}_{-d}+\left( 1-t_{p}\right)
\mathbf{\theta }_{0,-d}\right) \right\} _{\substack{p,q=1,...,d-1}}\right]
^{-1}\sqrt{n}\hat{S}^{*}\left( \mathbf{\theta }_{0,-d}\right),
\]
where $\hat{S}^{*}\left( \mathbf{\theta }_{0,-d}\right) =\left\{ \hat{S}%
_{p}^{*}\left( \mathbf{\theta }_{0,-d}\right) \right\} _{p=1}^{d-1}$ and
according to (\ref{DEF:etaip}) and Lemma \ref{LEM:sptheta0}
\[
\hat{S}_{p}^{*}\left( \mathbf{\theta }_{0,-d}\right)
=n^{-1}\sum_{i=1}^{n}\eta _{p,i}+o\left( n^{-1/2}\right) ,a.s.,\text{ }%
E\left( \eta _{p,i}\right) =0.
\]
Let $\Psi \left( \mathbf{\theta }_{0}\right) =\left( \psi _{pq}\right)
_{p,q=1}^{d-1}$ be the covariance matrix of $\sqrt{n}\left\{ \hat{S}%
_{p}^{*}\left( \mathbf{\theta }_{0,-d}\right) \right\} _{p=1}^{d-1}$ with $%
\psi _{pq}$ given in Theorem \ref{THM:normality}. Cram\'{e}r-Wold device and
central limit theorem for $\alpha $ mixing sequences entail that
\[
\sqrt{n}\hat{S}^{*}\left( \mathbf{\theta }_{0,-d}\right) \stackrel{d}{%
\longrightarrow }N\left\{ \mathbf{0},\Psi \left( \mathbf{\theta }_{0}\right)
\right\} .
\]
Let $\Sigma \left( \mathbf{\theta }_{0}\right) =\left\{ H^{*}\left( \mathbf{%
\theta }_{0,-d}\right) \right\} ^{-1}\Psi \left( \mathbf{\theta }_{0}\right)
\left[ \left\{ H^{*}\left( \mathbf{\theta }_{0,-d}\right) \right\}
^{T}\right] ^{-1}$, with $H^{*}\left( \mathbf{\theta }_{0,-d}\right) $ being
the Hessian matrix defined in (\ref{DEF:SHstarmatrices}). The above limiting
distribution of $\sqrt{n}\hat{S}^{*}\left( \mathbf{\theta }_{0,-d}\right)$, (%
\ref{EQ:Hessianconverge}) and Slutsky's theorem imply that
\[
\hspace{3.2cm}\sqrt{n}\left( \hat{\mathbf{\theta}}_{-d}\mathbf{-\theta }%
_{0,-d}\right) \stackrel{d}{\longrightarrow }N\left\{ \mathbf{0},\Sigma
\left( \mathbf{\theta }_{0}\right) \right\}.\hspace{3.5cm}
\]

\noindent{\large \textbf{References}}

\begin{description}
\item  Bosq, D. (1998). \emph{Nonparametric Statistics for Stochastic
Processes}. Springer-Verlag, New York.

\item  Carroll, R., Fan, J., Gijbles, I. and Wand, M. P. (1997). Generalized
partially linear single-index models. \textit{J. Amer. Statist. Assoc.}
\textbf{92} 477-489.

\item  Chen, H. (1991). Estimation of a projection -persuit type regression
model. \textit{Ann. Statist.} \textbf{19} 142-157.

\item  de Boor, C. (2001). \emph{A Practical Guide to Splines}.
Springer-Verlag, New York.

\item  DeVore, R. A. and Lorentz, G. G. (1993). \emph{Constructive
Approximation: Polynomials and Splines Approximation}. Springer-Verlag,
Berlin.

\item  Fan, J. and Gijbels, I. (1996). \emph{Local Polynomial Modelling and
Its Applications}. Chapman and Hall, London.

\item  Friedman, J. H. and Stuetzle, W. (1981). Projection pursuit
regression. \textit{J. Amer. Statist. Assoc.} \textbf{76} 817-823.

\item  H\"{a}rdle, W. (1990). \emph{Applied Nonparametric Regression}.
Cambridge University Press, Cambridge.

\item  H\"{a}rdle, W. and Hall, P. and Ichimura, H. (1993). Optimal
smoothing in single-index models. \textit{Ann. Statist.} \textbf{21} 157-178.

\item  H\"{a}rdle, W. and Stoker, T. M. (1989). Investigating smooth
multiple regression by the method of average derivatives. \textit{J. Amer.
Statist. Assoc.} \textbf{84} 986-995.

\item  Hall, P. (1989). On projection pursuit regression. \textit{Ann.
Statist.} \textbf{17} 573-588.

\item  Hastie, T. J. and Tibshirani, R. J. (1990). \emph{Generalized
Additive Models}. Chapman and Hall, London.

\item  Horowitz, J. L. and H\"{a}rdle, W. (1996). Direct semiparametric
estimation of single-index models with discrete covariates. \textit{J. Amer.
Statist. Assoc.} \textbf{91} 1632-1640.

\item  Hristache, M., Juditski, A. and Spokoiny, V. (2001). Direct
estimation of the index coefficients in a single-index model.
\textit{Ann. Statist.} \textbf{29} 595-623.

\item  Huang, J. Z. (2003). Local asymptotics for polynomial spline
regression. \textit{Ann. Statist.} \textbf{31} 1600-1635.

\item  Huang, J. and Yang, L. (2004). Identification of nonlinear additive
autoregressive models. \textit{J. R. Stat. Soc. Ser. B Stat. Methodol.}
\textbf{66} 463-477.

\item  Huber, P. J. (1985). Projection pursuit (with discussion). \textit{%
Ann. Statist.} \textbf{13} 435-525.

\item  Ichimura, H. (1993). Semiparametric least squares (SLS) and weighted
SLS estimation of single-index models \textit{Journal of Econometrics}
\textbf{58} 71-120.

\item  Klein, R. W. and Spady. R. H. (1993). An efficient semiparametric
estimator for binary response models. \textit{Econometrica} \textbf{61}
387-421.

\item  Mammen, E., Linton, O. and Nielsen, J. (1999). The existence and
asymptotic properties of a backfitting projection algorithm under weak
conditions. \textit{Ann. Statist.} \textbf{27} 1443-1490.

\item  Pham, D. T. (1986). The mixing properties of bilinear and generalized
random coefficient autoregressive models. \textit{Stochastic Anal. Appl.}
\textbf{23} 291-300.

\item  Powell, J. L., Stock, J. H. and Stoker, T. M. (1989). Semiparametric
estimation of index coefficients. \textit{Econometrica.} \textbf{57}
1403-1430.

\item  Tong, H. (1990) \emph{Nonlinear Time Series: A Dynamical System
Approach}. Oxford, U.K.: Oxford University Press.

\item  Tong, H., Thanoon, B. and Gudmundsson, G. (1985) Threshold time
series modeling of two icelandic riverflow systems. \textit{Time Series
Analysis in Water Resources}. ed. K. W. Hipel, American Water Research
Association.

\item  Wang, L. and Yang, L. (2007). Spline-backfitted kernel
smoothing of nonlinear additive autoregression model. \textit{Ann.
Statist.} Forthcoming.

\item  Xia, Y. and Li, W. K. (1999). On single-index coefficient regression
models. \textit{J. Amer. Statist. Assoc.} \textbf{94} 1275-1285.

\item  Xia, Y., Li, W. K., Tong, H. and Zhang, D. (2004). A goodness-of-fit
test for single-index models. \textit{Statist. Sinica.} \textbf{14} 1-39.

\item  Xia, Y., Tong, H., Li, W. K. and Zhu, L. (2002). An adaptive
estimation of dimension reduction space. \textit{J. R. Stat. Soc. Ser. B
Stat. Methodol.} \textbf{64} 363-410.

\item Xue, L. and Yang, L. (2006 a). Estimation of semiparametric
additive coefficient model. \textit{J. Statist. Plann. Inference}
\textbf{136}, 2506-2534.

\item Xue, L. and Yang, L. (2006 b). Additive coefficient modeling
via polynomial spline. \textit{Statistica Sinica} \textbf{16}
1423-1446.

\end{description}

%\vskip .65cm \noindent Department of Statistics, University of
%Georgia,
%Athens, GA 30602, USA \vskip 2pt \noindent E-mail: (wangli4@stt.msu.edu) %
%\vskip 2pt \noindent Department of Statistics and Probability,
%Michigan State University, East Lansing, MI 48824, USA \vskip 2pt
%\noindent E-mail: (wangli4@stt.msu.edu) \vskip .3cm

%\centerline{(Received xxx 200?; accepted xxx 200?)}\par

\setcounter{chapter}{8} \renewcommand{\thefigure}{{\arabic{figure}}} %
\setcounter{figure}{0}

% This is the figure of Xia 2004 estimators
\newpage \setlength{\unitlength}{1cm} \pagestyle{empty}
% This is the figure of spline estimators
\newpage \setlength{\unitlength}{1cm} \pagestyle{empty}
\begin{figure}[th]
\begin{center}
\begin{picture}(-19.5,23)
\put(-19.2,11.4){\includegraphics{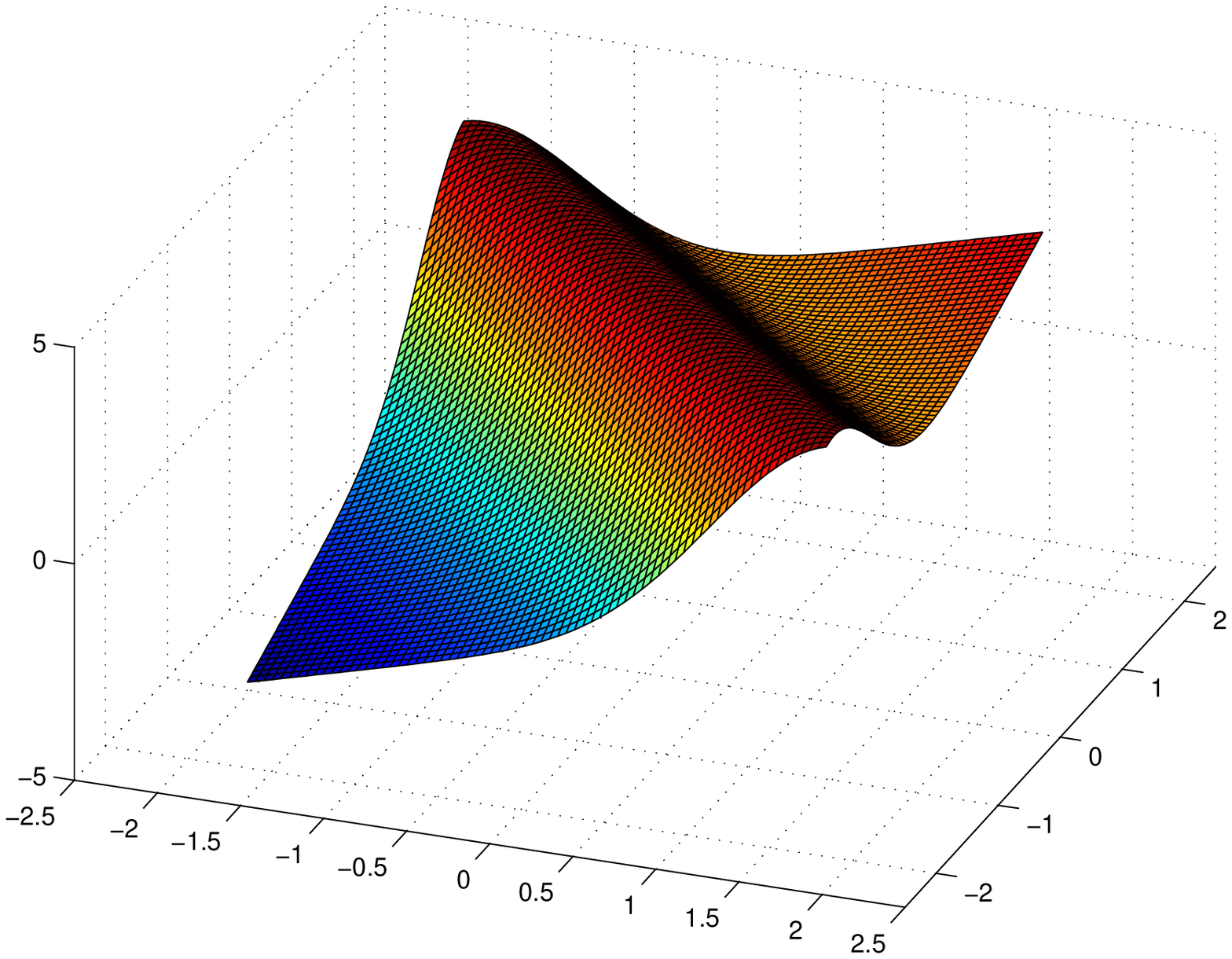}}%
\put(-14.2,14.0){(a)}%
\put(-10.5,11.4){\includegraphics{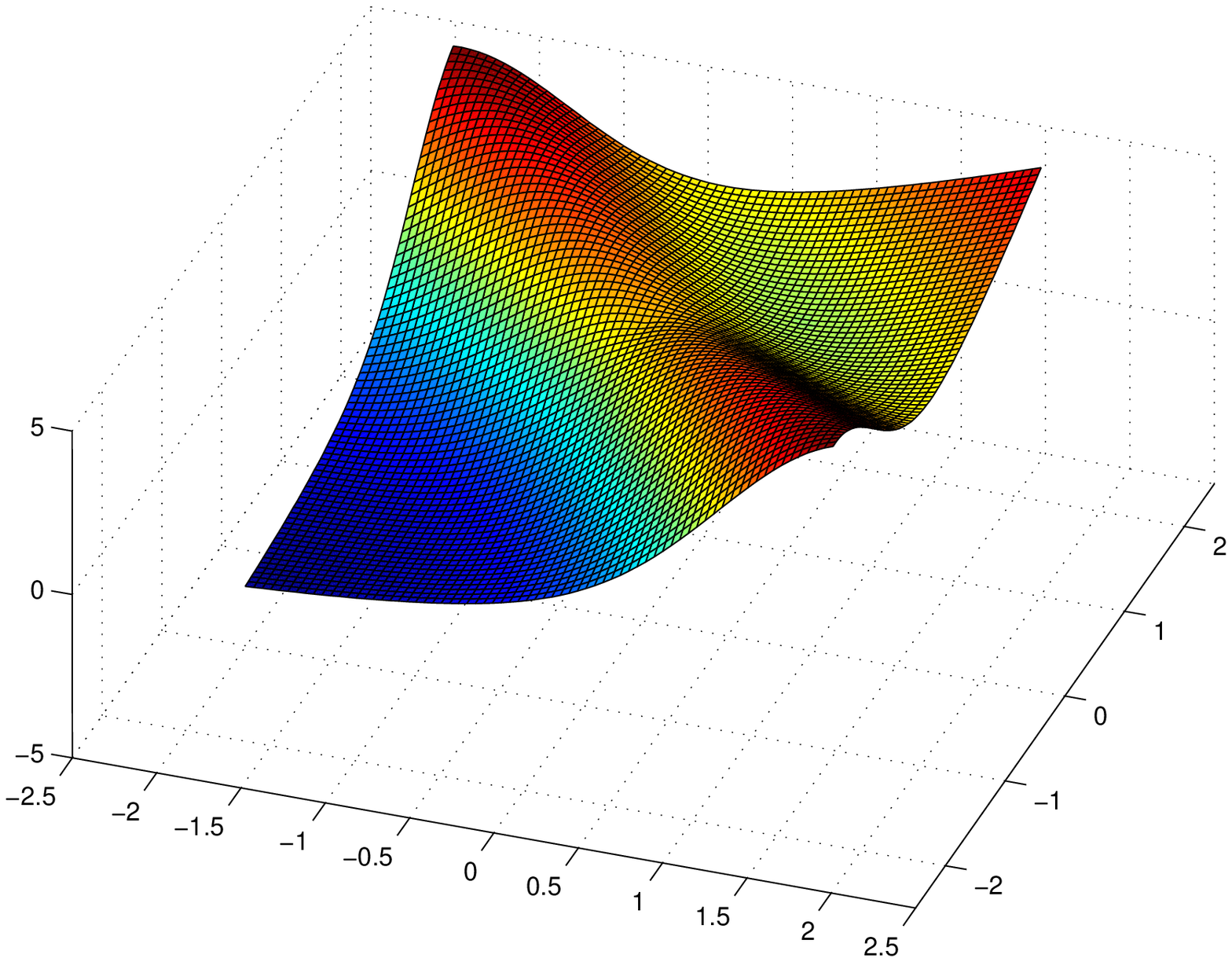}}%
\put(-5.6,14.0){(b)}%
\put(-21,0.4){\includegraphics{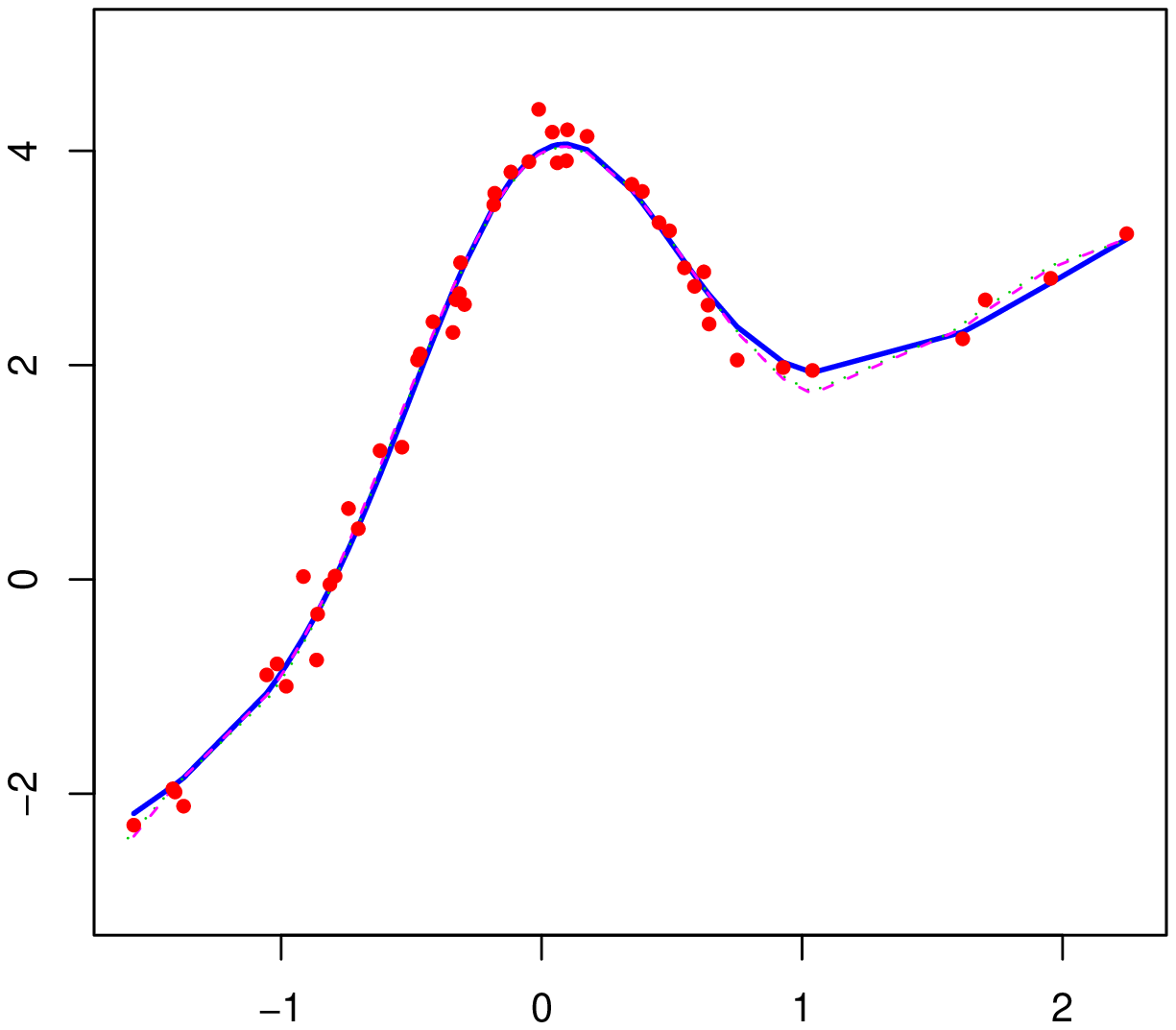}}%
\put(-14.2,5.2){(c)}%
\put(-12.5,0.4){\includegraphics{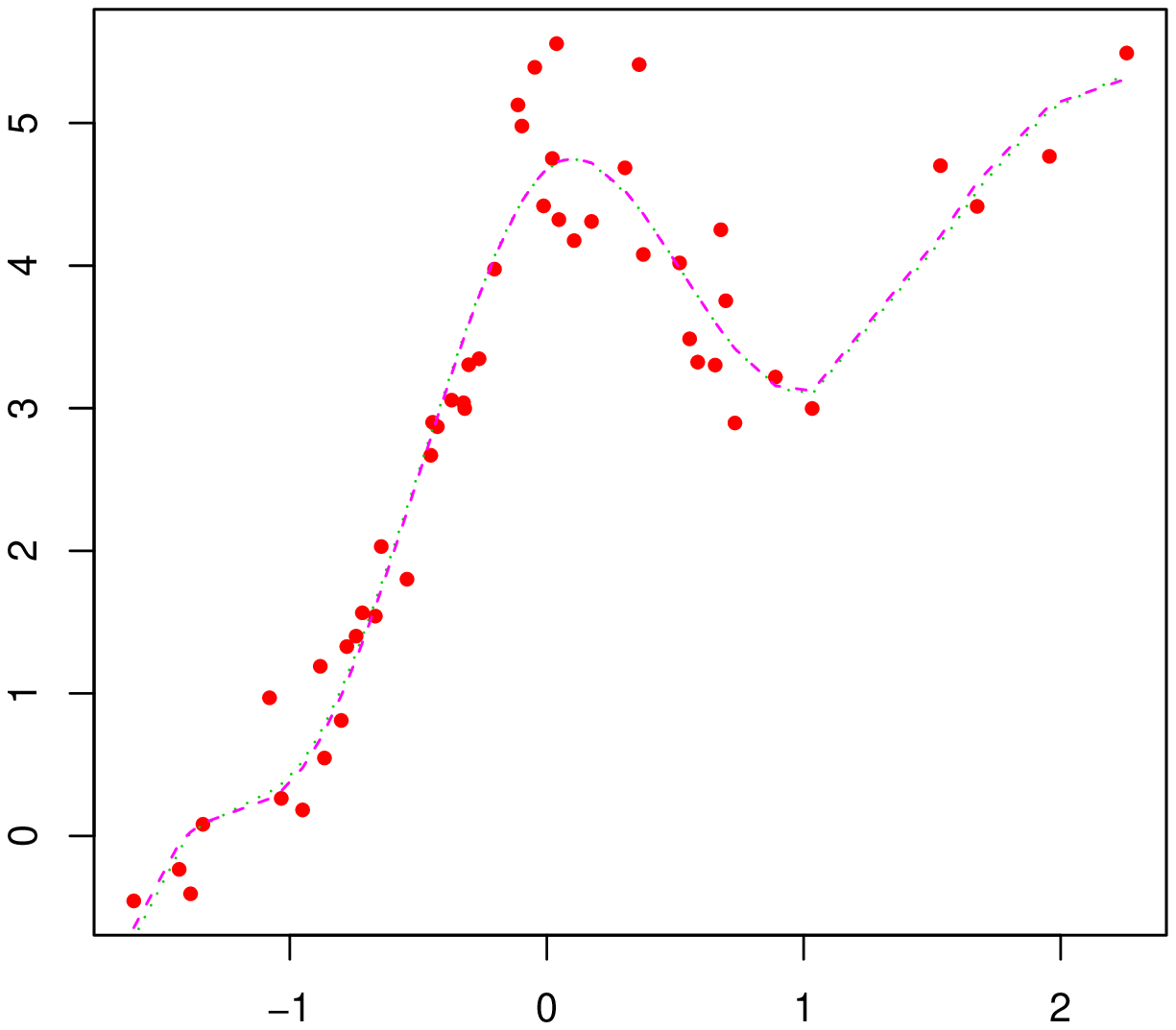}}%
\put(-5.6,5.2){(d)}%
\end{picture}
\end{center}
\par
\vskip -5.5cm \caption{Example 1. (a) and (b) Plots of the actual
surface $m$ in model (\ref {MODEL:Xia}) with respect to
$\delta=0,1$; (c) and (d) Plots of various univariate functions
with respect to $\delta=0,1$: $\left\{
\mathbf{X}_i^T\hat{\mathbf{\theta}},Y_i\right\}, {1 \leq i \leq
50}$ (dots); the univariate function $g$ (solid line); the
estimated
function of $g$ by plugging in the true index coefficient $\mathbf{%
\theta}_0$ (dotted line); the estimated function of $g$ by
plugging in the estimated index coefficient (dashed line)
$\hat{\mathbf{\theta}}=(0.69016, 0.72365)^T$ for $\delta=0$ and
$(0.72186, 0.69204)^T$ for $\delta=1$.} \label{FIG:Xia2004}
\end{figure}

% This is the figure of spline estimators
\newpage \setlength{\unitlength}{1cm} \pagestyle{empty}
\begin{figure}[th]
\begin{center}
\begin{picture}(-19.5,23)
\put(-21,9.4){\includegraphics{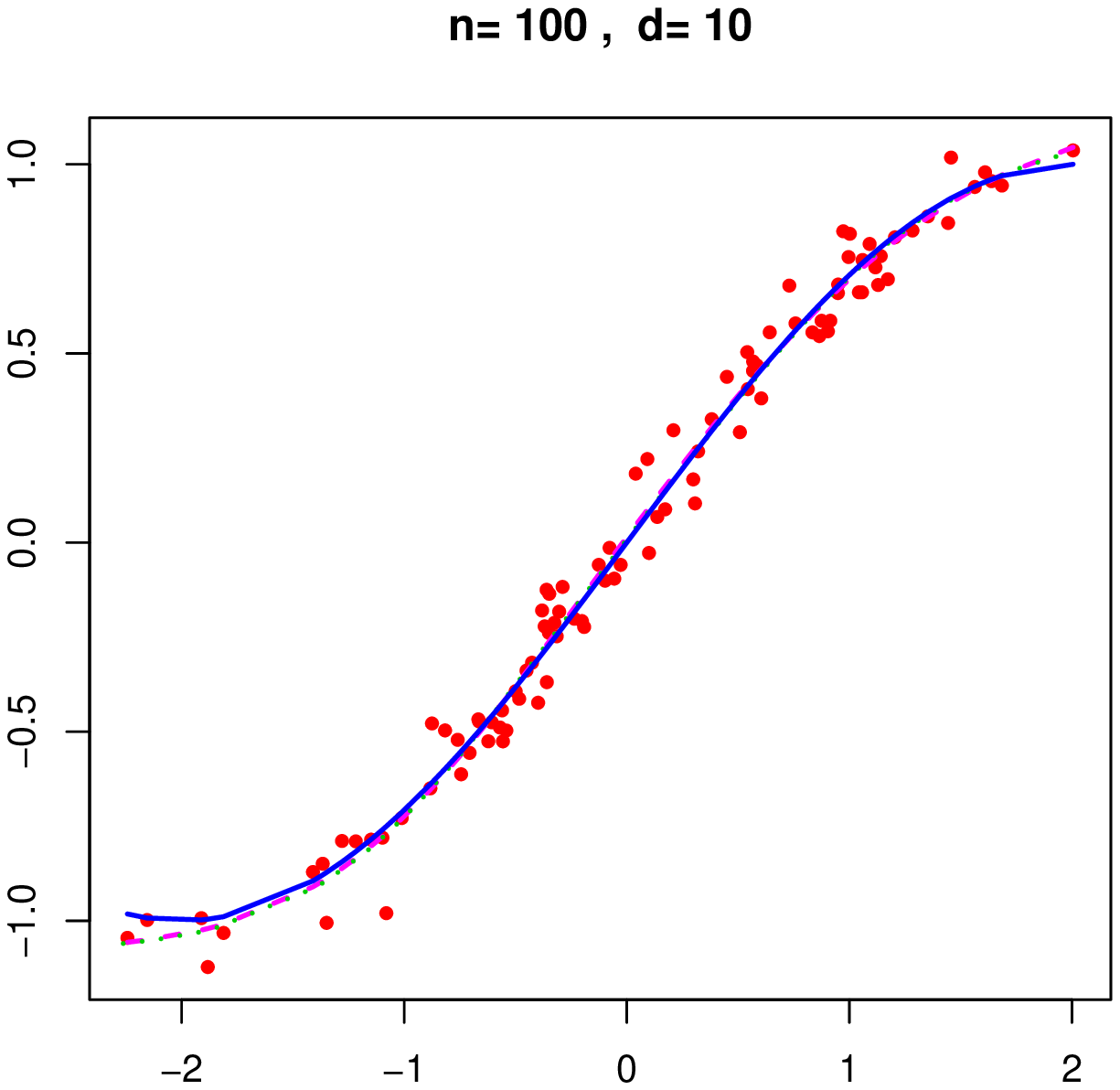}}%
\put(-14.2,14.0){(a)}%
\put(-12.5,9.4){\includegraphics{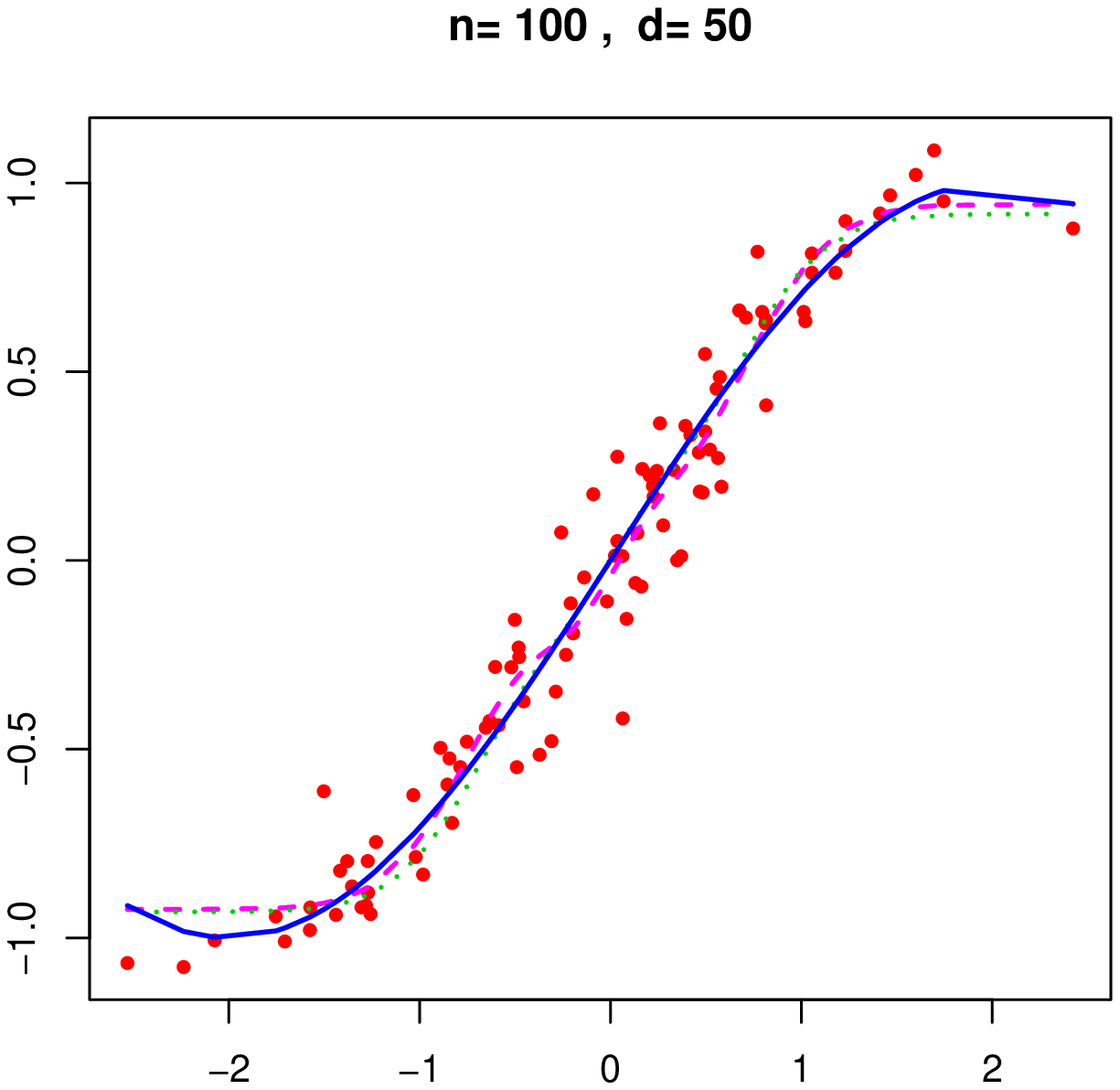}}%
\put(-5.6,14.0){(b)}%
\put(-21,-0.8){\includegraphics{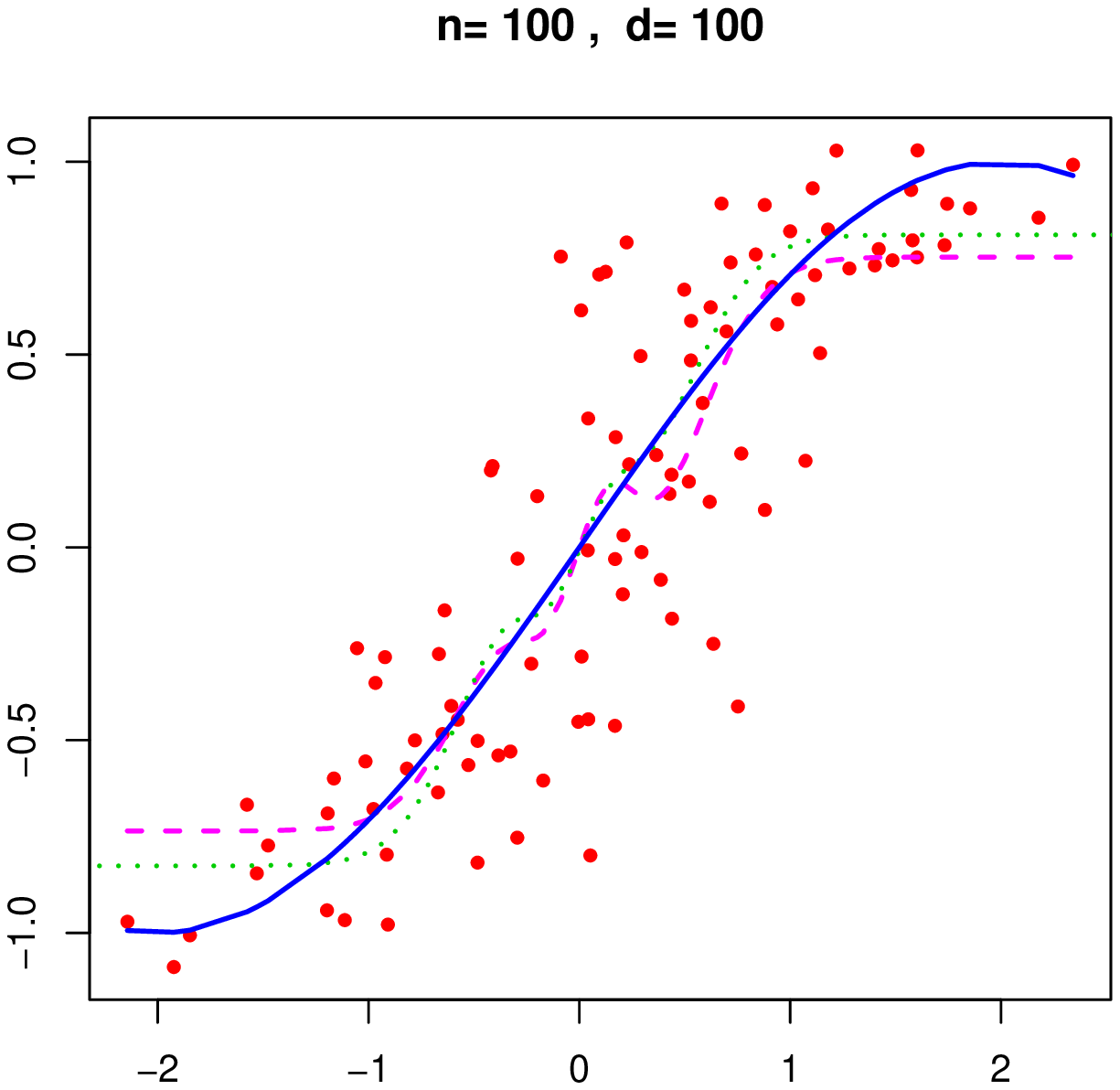}}%
\put(-14.2,3.7){(c)}%
\put(-12.5,-0.8){\includegraphics{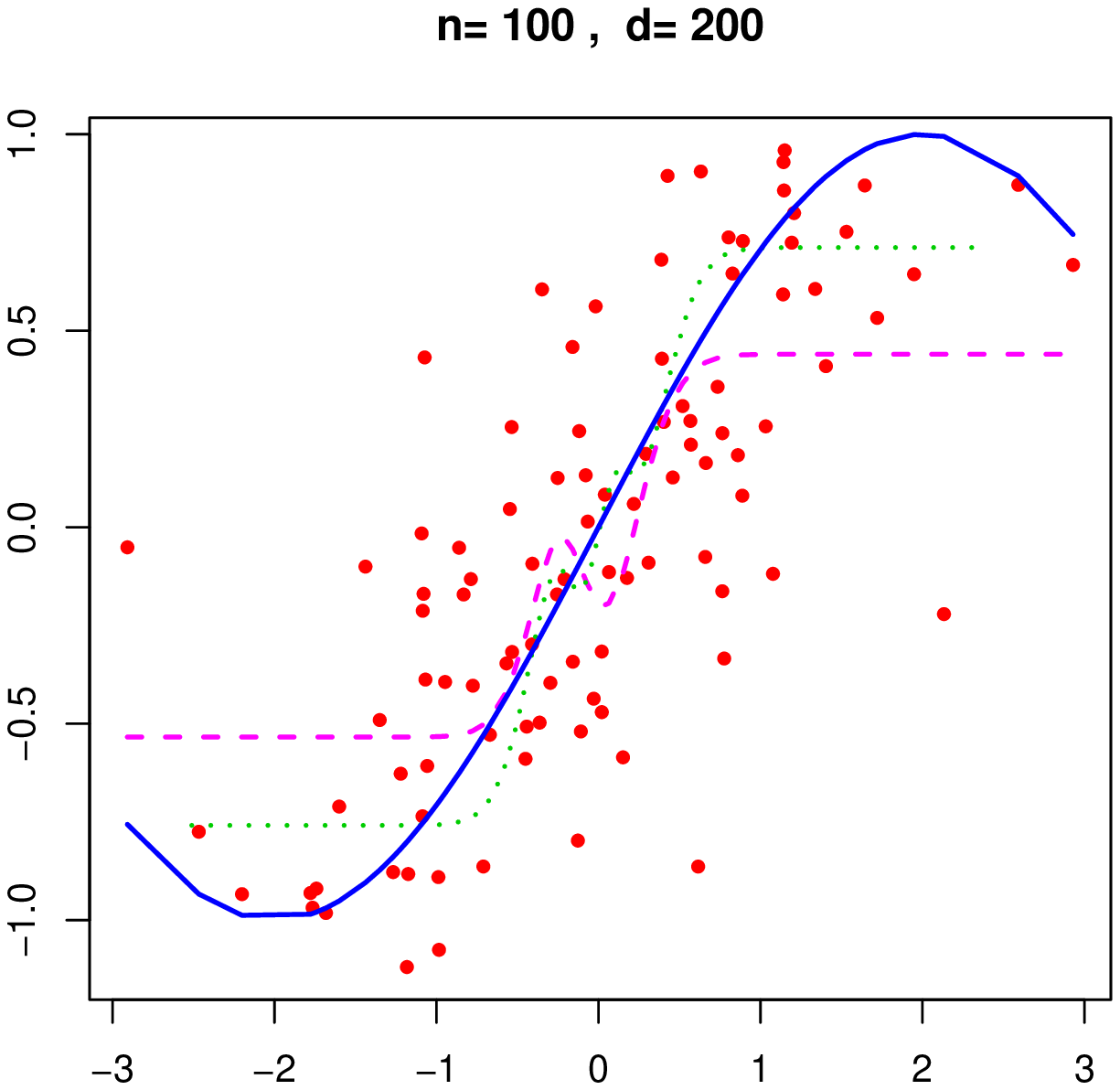}}%
\put(-5.6,3.7){(d)}%
\end{picture}
\end{center}
\par
\vskip -3.8cm \caption{Example 2. Plots of the spline estimator of
$g$ with the estimated index parameter $\hat{\mathbf{\theta}}$
(dotted curve), cubic spline estimator of
$g$ with the true index parameter $\mathbf{\theta}_0$ (dashed curves), the true function $m\left(\mathbf{x}\right)$ in (%
\ref{MODEL:simulation}) (solid curve), and the data scatter plots (dots).}
\label{FIG:estimation}
\end{figure}

\newpage % This is the figure to see the convergence behavior
\setlength{\unitlength}{1cm} \pagestyle{empty}
\begin{figure}[th]
\begin{center}
\begin{picture}(-19.5,24)
\put(-21,9.4){\includegraphics{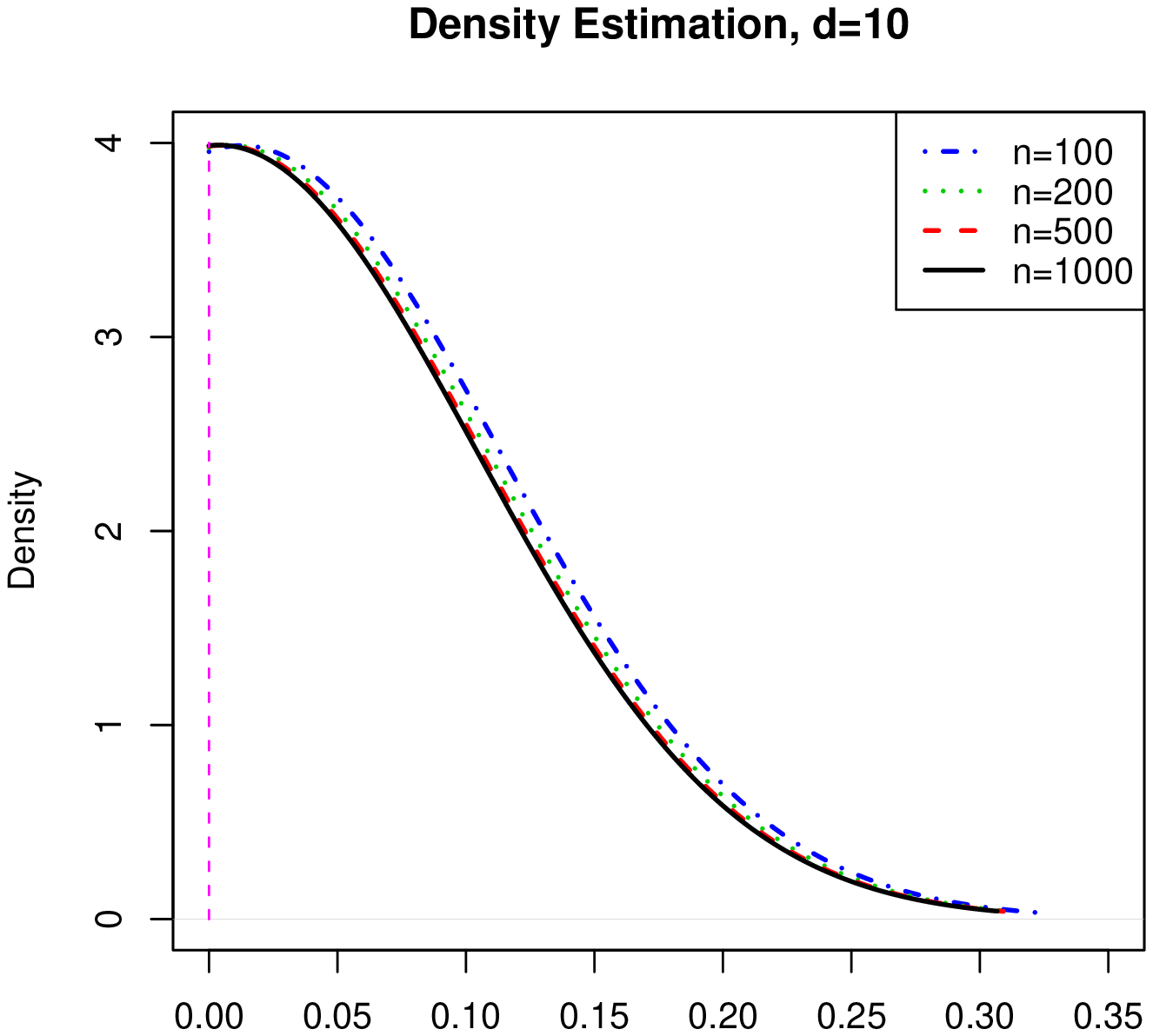}}%
\put(-14.2,14.2){(a)}%
\put(-12.5,9.4){\includegraphics{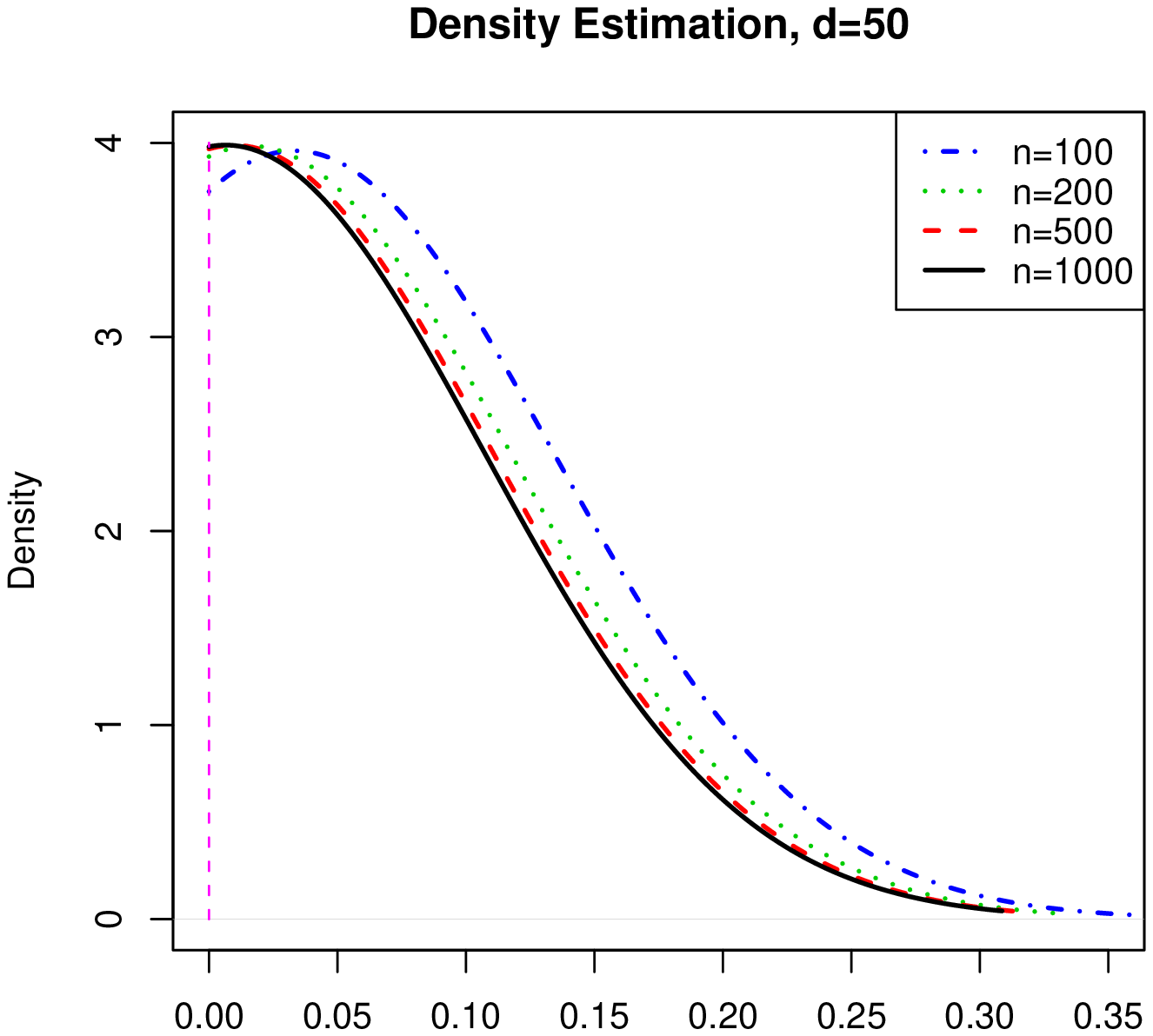}}%
\put(-5.6,14.2){(b)}%
\put(-21,-0.5){\includegraphics{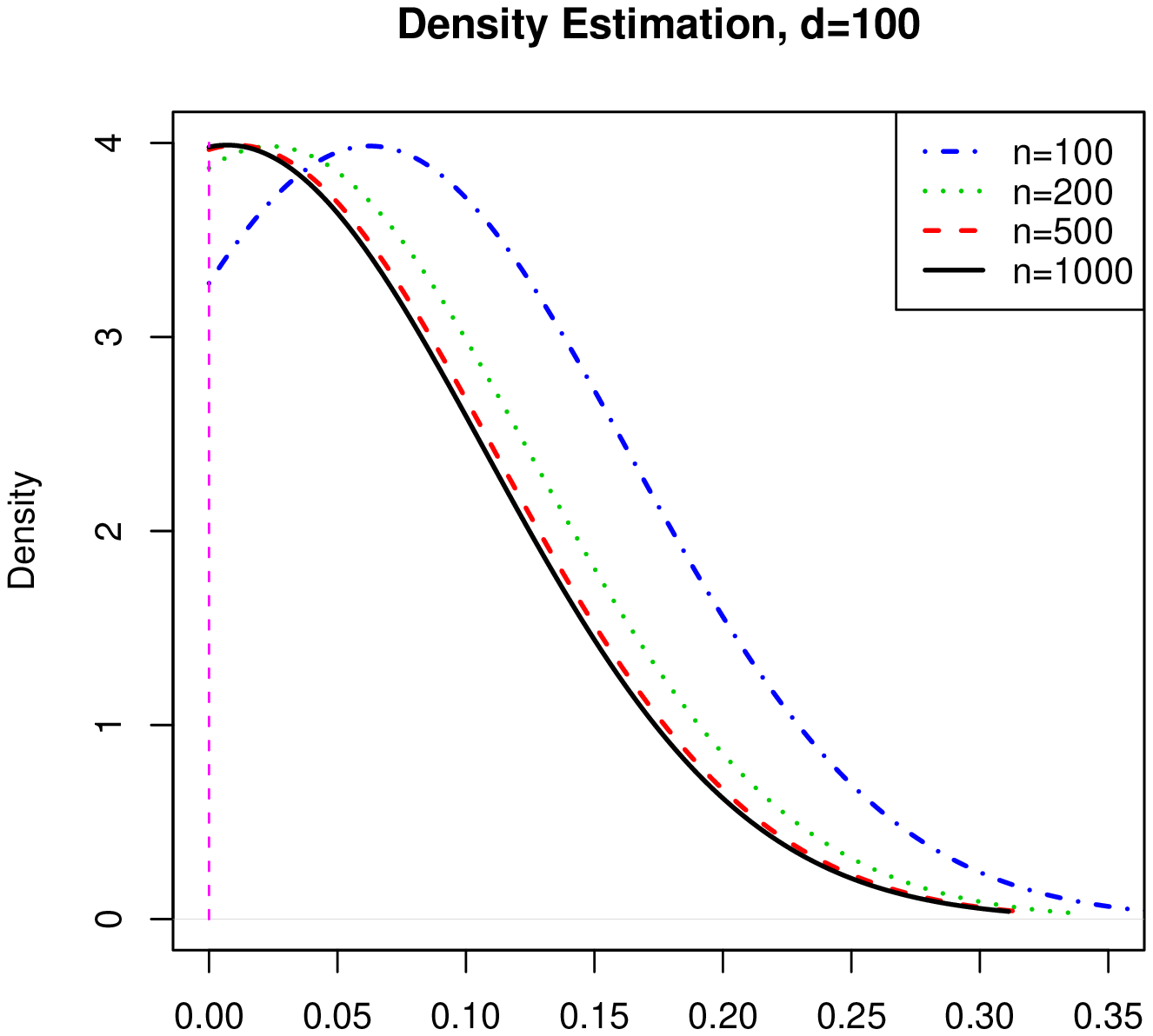}}%
\put(-14.2,4.0){(c)}%
\put(-12.5,-0.5){\includegraphics{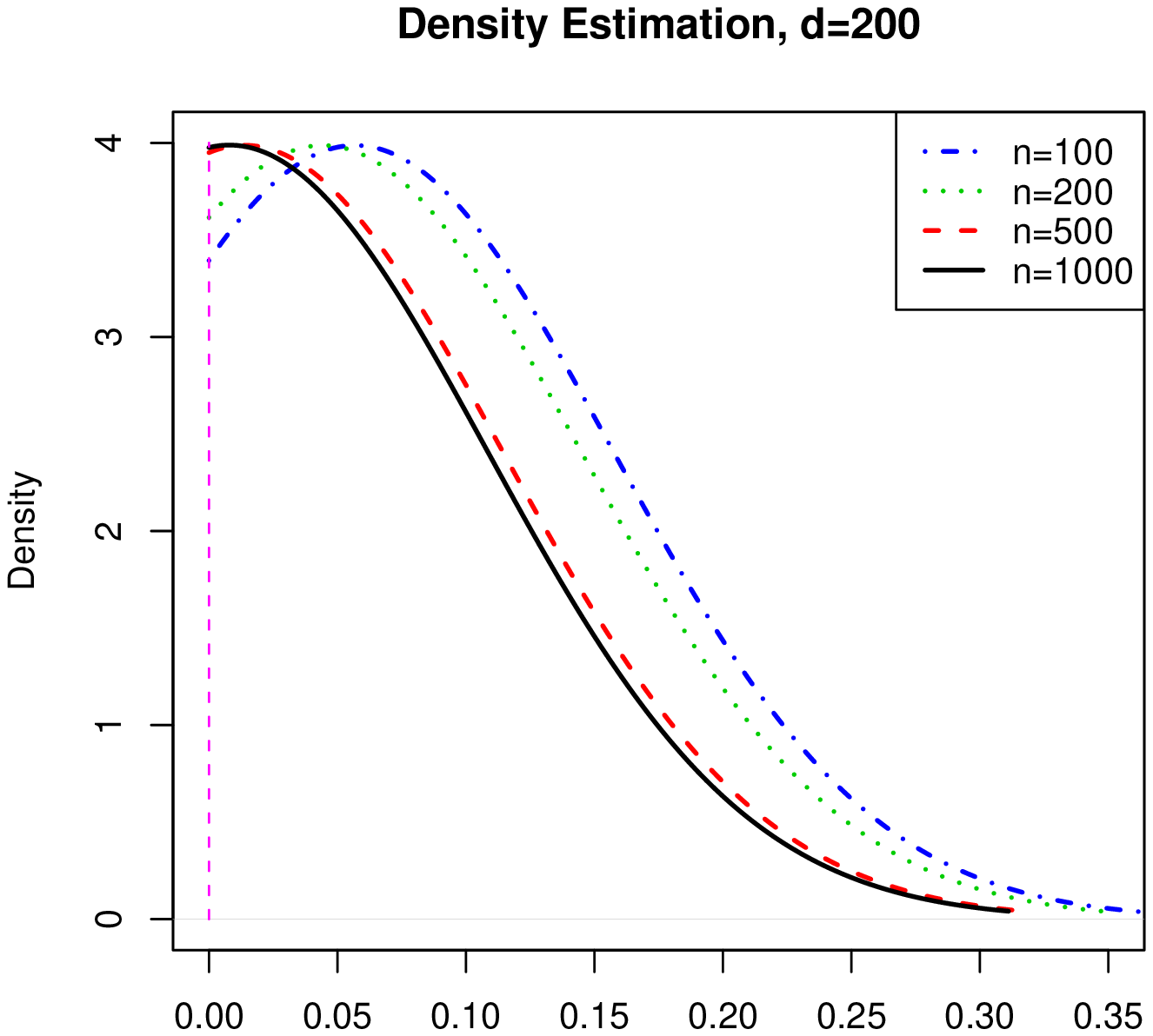}}%
\put(-5.6,4.0){(d)}%
\end{picture}
\end{center}
\par
\vskip -4.0cm
\caption{Example 2. Kernel density estimators of the $100$ $\|\hat{\mathbf{%
\theta}}-\mathbf{\theta}_0\|/\protect\sqrt{d}$.}
\label{FIG:density}
\end{figure}

\newpage % This is the figure to see the time plots with trends
\setlength{\unitlength}{1cm} \pagestyle{empty}
\begin{figure}[th]
\begin{center}
\begin{picture}(-19.5,20)
\put(-16.5,9.4){\includegraphics{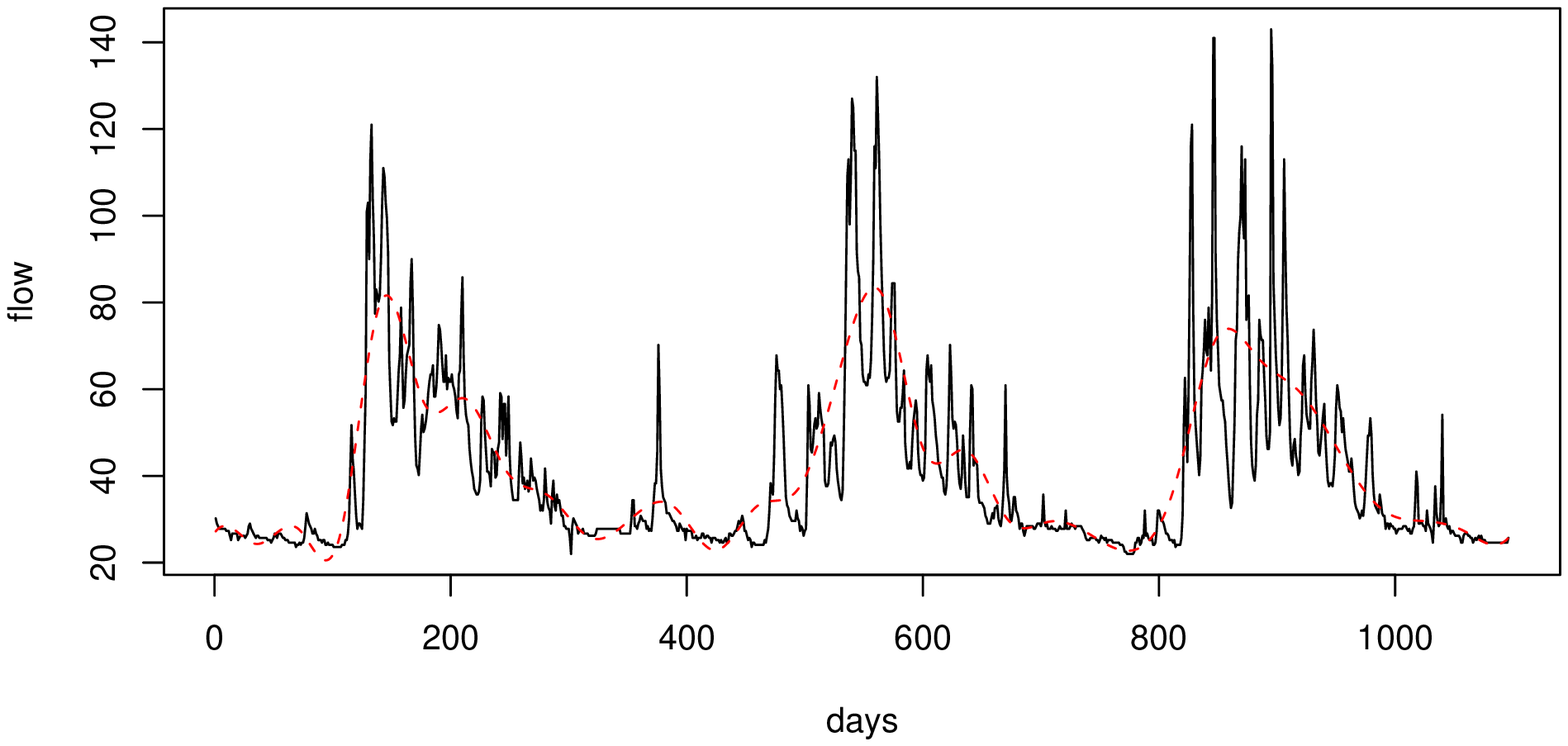}}%
\put(-9.4,13.8){(a)}%
\put(-16.5,3.4){\includegraphics{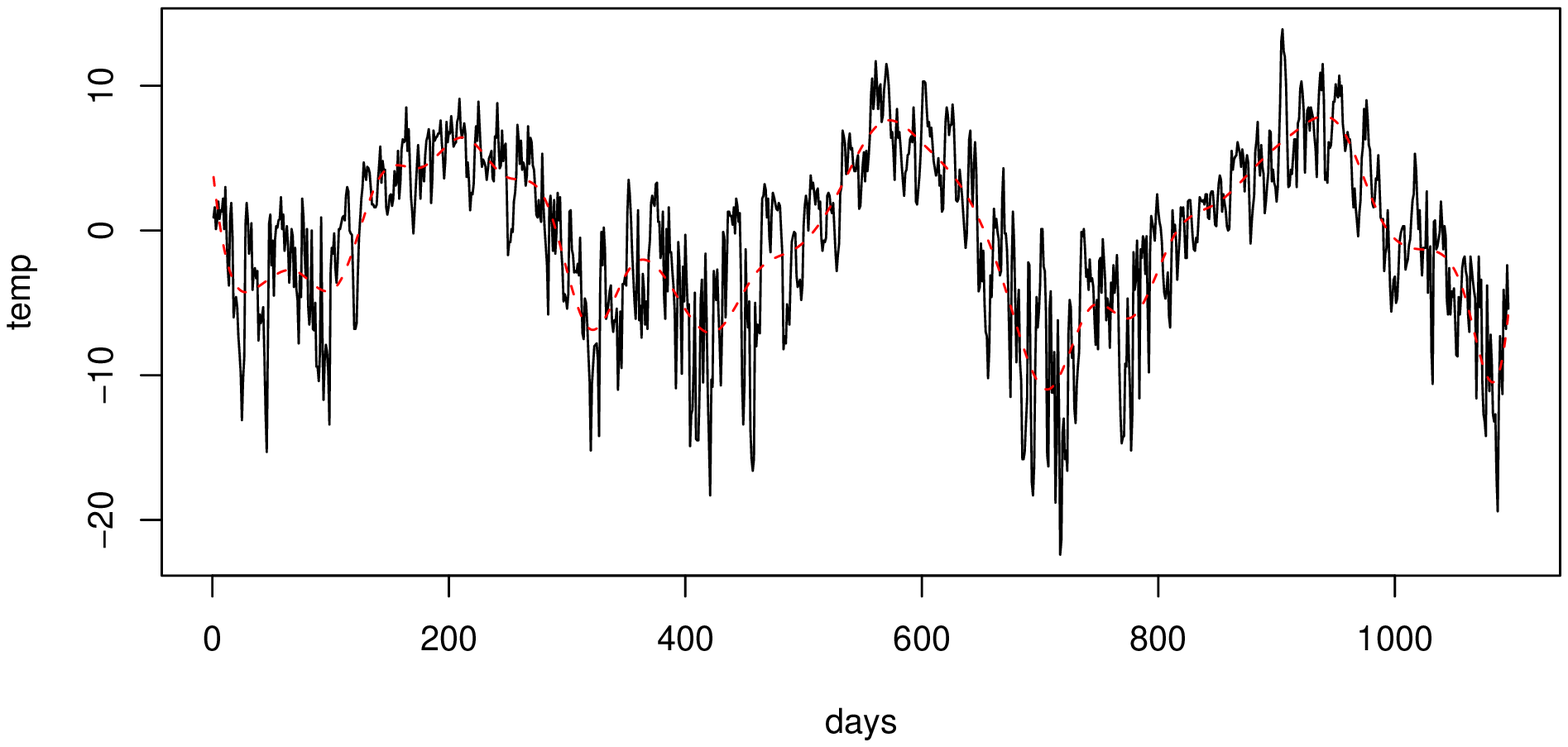}}%
\put(-9.4,7.8){(b)}%
\put(-16.5,-2.6){\includegraphics{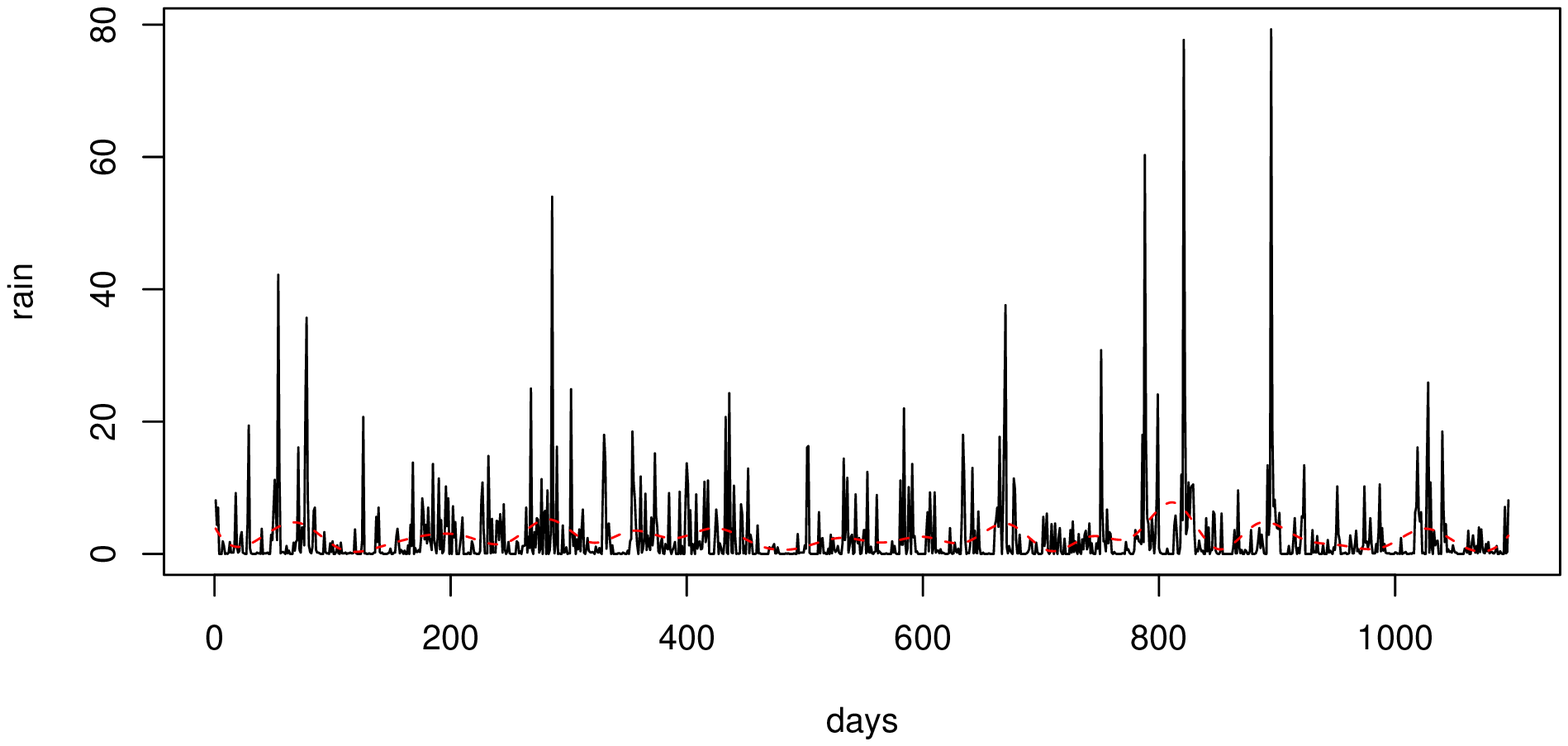}}%
\put(-9.4,1.8){(c)}%
\end{picture}
\end{center}
\par
\vskip -2.0cm \caption{Time plots of the daily J\"{o}kuls\'{a}
Eystri River data (a) river flow $Y_t$ (solid line) with its trend
(dashed line) (b) temperature $X_t$ (solid line) with its trend
(dashed line) (c) precipitation $Z_t$ (solid line) with its trend
(dashed line).} \label{FIG:riverflow}
\end{figure}

\newpage % This is the figure to see the convergence behavior
\setlength{\unitlength}{1cm} \pagestyle{empty}
\begin{figure}[th]
\begin{center}
\begin{picture}(-19.5,20)
\put(-16.5,9.4){\includegraphics{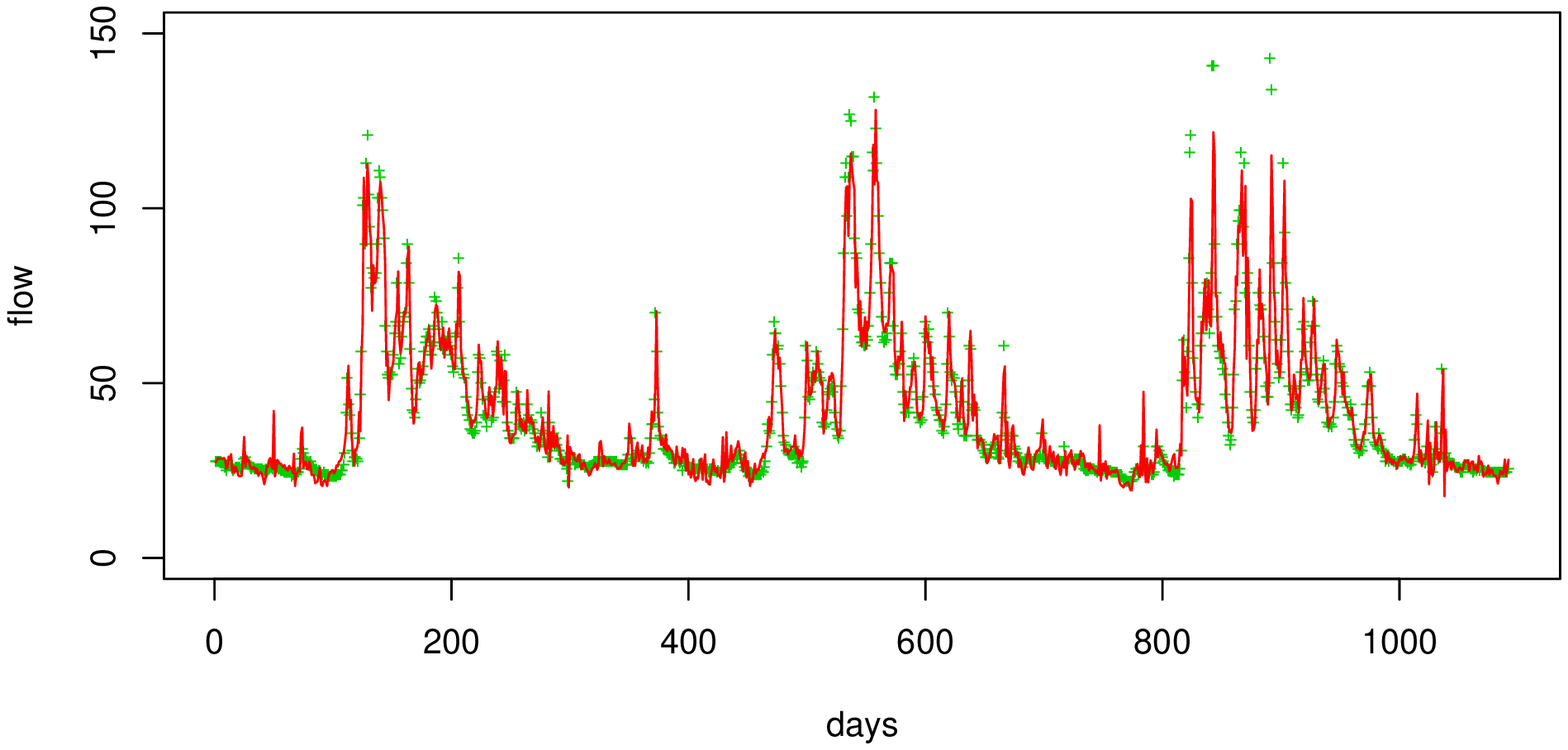}}%
\put(-9.4,13.8){(a)}%
\put(-16.5,3.4){\includegraphics{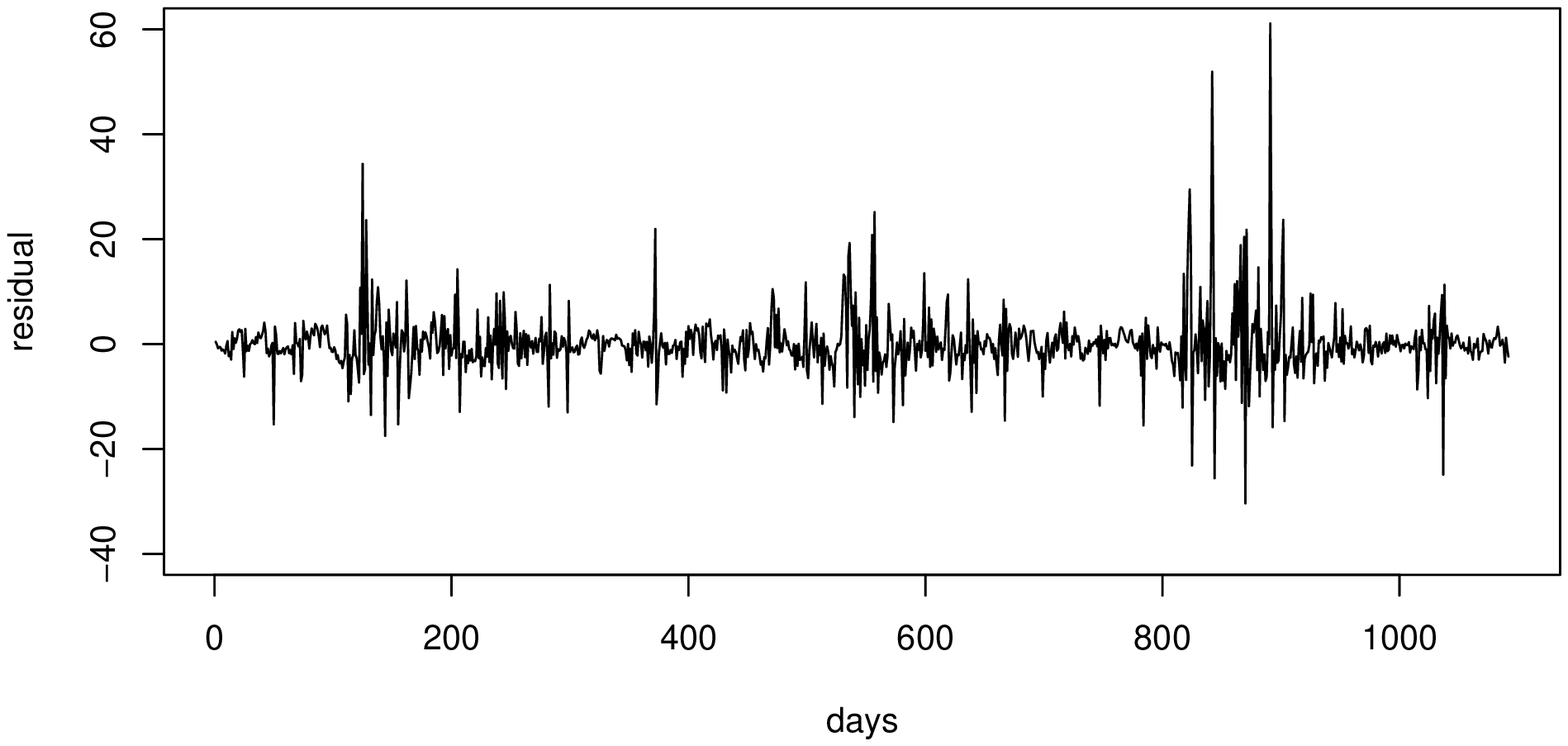}}%
\put(-9.4,7.8){(b)}%
\put(-16.5,-2.6){\includegraphics{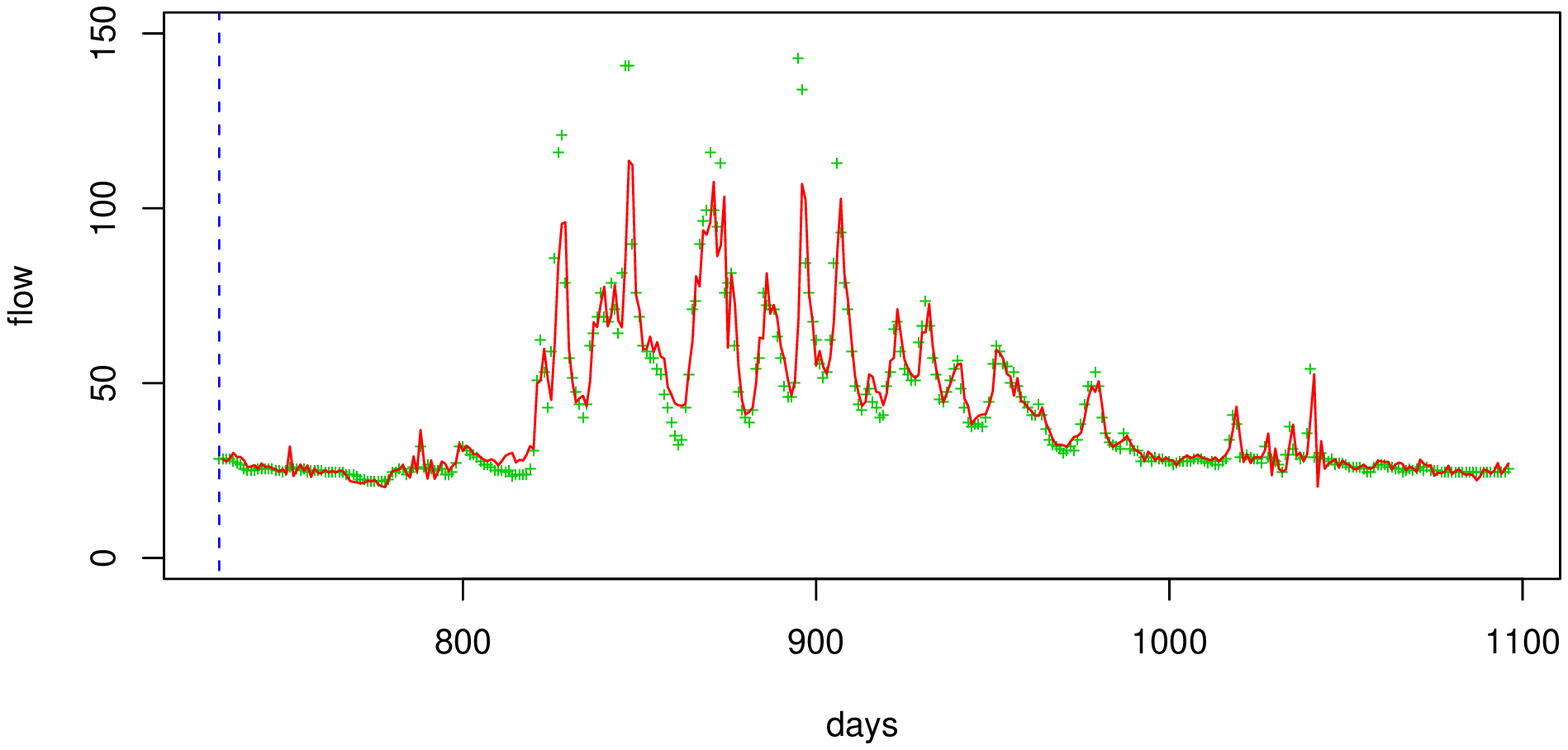}}%
\put(-9.4,1.8){(c)}%
\end{picture}
\end{center}
\par
\vskip -2.0cm \caption{(a) The scatter plot of the river flow
(``+") and the fitted plot of the river flow (line) and (b)
Residuals of the fitted SIP model (c) Out-of-sample rolling
forecasts (line) of the river flow for the entire third year
(``+") based on the first two years' river flow.}
\label{FIG:riverflow_fitted}
\end{figure}
\end{document}